\documentclass[final]{siamltex}
\usepackage{amsmath}
\usepackage{amssymb}
\usepackage{amstext}

\usepackage{graphicx}
\usepackage{hhline}
\usepackage{psfrag}
\usepackage{float}
\usepackage[dvips]{epsfig}
\usepackage{subfig}
\usepackage{tikz}
\usetikzlibrary{decorations.pathreplacing}
\usepackage{epstopdf} 
\usepackage{color}

\newtheorem{rem}{Remark}

\begin{document}
%	\nocite{*}
\title{A meshfree arbitrary Lagrangian-Eulerian  method for the BGK model of the Boltzmann equation with moving boundaries}

\author{ S. Tiwari  \footnotemark[1] , 
A. Klar \footnotemark[1] \footnotemark[2]  \and G. Russo \footnotemark[3] }
\footnotetext[1]{Technische Universit\"at Kaiserslautern, Department of Mathematics, Erwin-Schr\"odinger-Stra{\ss}e, 67663 Kaiserslautern, Germany 
  (\{klar, tiwari\}@mathematik.uni-kl.de)}
\footnotetext[2]{Fraunhofer ITWM, Fraunhoferplatz 1, 67663 Kaiserslautern, Germany} 
\footnotetext[3]{Department of Mathematics and Computer Science, University of Catania, Italy (russo@dmi.unict.it)}

%%%%%%%%%%%%%%%%%%%%%%%%%%%%%%%%%%%%%%%%%%%%%%%%%%%%%%%%%%%

\maketitle
In this paper we present a novel technique for the simulation of moving boundaries and moving rigid bodies immersed in a rarefied gas using an Eulerian-Lagrangian formulation based on least square method. The rarefied gas is simulated by solving the Bhatnagar-Gross-Krook (BGK) model for the Boltzmann equation of rarefied gas dynamics.  
  The BGK model is solved by an Arbitrary Lagrangian-Eulerian (ALE)  method, where grid-points/particles are moved with the mean velocity of the gas.  
 The computational domain  for the  rarefied gas  changes with time due to the motion of the boundaries. To allow a simpler handling of the interface motion  we have used a meshfree method based on a  least-square approximation  for the reconstruction procedures required for the scheme.
 We have considered a one way, as well as a two-way coupling of 
 boundaries/rigid bodies and gas flow.    The  numerical results are compared with analytical
 as well as with Direct Simulation Monte Carlo (DSMC) solutions of the Boltzmann equation. Convergence studies are performed for one-dimensional and two-dimensional  test-cases.  Several further test problems and 
 applications illustrate the versatility of the approach. 
\begin{abstract}
 %%%

%%%%
 \end{abstract}

{\bf MSC2020:} 65C05, 65M99,  70E99, 76P05, 76T20\\
{\bf Keywords:} rarefied gas, kinetic equation, BGK model, meshfree method, ALE method, semi-implicit method,  least squares method, gas rigid body interactions

\section{Introduction}
In  recent years moving boundary problems for rarefied gas dynamics have been extensively investigated  in the connection with Micro-Electro-Mechanical-Systems (MEMS), see \cite{BTSKH, DM,DM1, FFL, KBA, RF, STK, STKH, TKHD, TA,TA1}. 
%\giovanni{For the order of bibliography we should use either alphabetical order or appearing order. Please, correct}. 
In micro scale geometries the mean free path is often of the order  or larger than the characteristic length of the geometry, even at standard condition of temperature and pressure, thus requiring the
physical system to be described by kinetic  equations. Usually, these  flows have low Mach numbers, therefore, stochastic methods like DSMC  are not the optimal choice, since  statistical noise dominates the flow quantities. Moreover, when one considers  moving rigid body, the gas domain will change in time and one has  to encounter  unsteady flow problems, so that  averages over long runs cannot be taken. Instead, one has to perform many independent runs in order to get smooth solutions. 
Although some attempts have been made to reduce the statistical noise of DSMC type methods, see, for example, \cite{DDP}, or to adopt efficient solvers for the Boltzmann equation, such as those based on the Fourier-spectral method (see for example the review paper \cite{DP}), many works rather employ deterministic approaches  for  simplified models of  the Boltzmann equation, like the Bhatnagar-Gross-Krook (BGK) model, see \cite{DM, RF, TA,TKR,TKR2}. 
In the above mentioned works either Finite-Difference schemes or Semi-Lagrangian methods are used to solve the moving boundary problems, see \cite{DM} 
%\giovanni{are we sure this is a review paper . Yes. It is not a review paper but contains a nice survey of the literature in the introduction} 
for an overview  of methods used for the BGK equation. 
Since the rigid body moves in time,  classical 
interpolation procedures near the rigid body become complicated and possibly inaccurate because of the arbitrary intersection of cells by the rigid body. Thus,  a  Cartesian cut cell method has been introduced in 
\cite{DM1} to handle the moving object in the rarefied gas. A different technique has been used in \cite{CCKR}, where  the authors have used  ghost point methods in a finite difference framework to treat moving boundaries. 
For immersed boundary type approaches  applied to kinetic equations to simulate the fluid-rigid body interactions
see \cite{AKF, DM,TKR2}.

In the present  paper we use a deterministic Arbitrary Lagrangian-Eulerian approach  for the   BGK model. First and second order versions of the scheme and associated  upwinding procedure are described and numerically tested.
This approach, based on moving grid points, is simple, well suited and very efficient for the treatment of problems with moving boundaries. While the interior grid points are moved with the mean velocity of the gas,
the moving  boundaries  are as well  approximated by a discrete set of boundary points moving with the boundaries.
This leads to a very flexible scheme also suited  for complicated geometries and flows.

The paper is organised as follows. In section 2 we present the BGK model for the Boltzmann equation,   the Newton-Euler equations for rigid body motions and the Chu reduction procedure. In section 3 we introduce the numerical  scheme for the BGK model, in particular the spatial and temporal discretization with first and second order accuracy.  Section 4 illustrates various numerical results in one and two space dimensions including a convergence study in 1D and 2D and comparisons with DSMC results. Finally, in section 5 some conclusions and an outlook  are presented.

\section{The BGK model for rarefied gas dynamics}
 \label{sec:model}
    
We consider the BGK  model of the Boltzmann equation for  rarefied gas dynamics, where the collision term is modeled by a relaxation of the distribution function $f(t,x,v)$ to the Maxwellian equilibrium distribution. The evolution equation for the distribution function $f(t,x,v)$  is given by the following initial boundary value problem
\begin{equation}
\frac{\partial f}{\partial t} + v\cdot\nabla_x f = \frac{1}{\tau}(M - f)
\label{bgk_eqn}
\end{equation}
with $t \ge 0, x \in \Omega \subset \mathbb R^{d_x}, (d_x=1,2,3), \;  v  \in \mathbb R^{d_v}, (d_v=1,2,3)$ and initial condition $f(0,x, v) = f_0(x, v)$. Additionally, suitable  boundary conditions are  described, see   the next section. 
 Here $\tau$ is the relaxation time, which may depend on local density and temperature, and $M$ is the local Maxwellian given by 
\begin{equation}
M = \frac{\rho}{(2\pi R T)^{d_v/2}} \exp \left(\frac{| v -  U|^2}{2RT}\right), 
\label{maxwellian}
\end{equation}
where the parameters $\rho (x,t) \in \mathbb R , U (x,t) \in \mathbb R^{d_v} , T (x,t)\in \mathbb R$ are the macroscopic quantities mass density, mean velocity and temperature, respectively.  $R$ is the universal gas constant divided by the molecular mass of the gas.  $\rho, U ,T $ are computed from $f$ as  follows.
Let the moments of $f$ be  defined by 
\begin{equation}
(\rho, \rho U, E) = \int_{\mathbb R^{3}} \phi(v) f(t,  x,  v) d v.
\label{moments}
\end{equation}
where  $ { \phi } (v)=\left (1,  v ,\frac{| v|^2}{2} \right )$ denotes  the vector of collision invariants.
$E$ is the total energy density which is related to the temperature through the internal energy 
\begin{equation}
e(t, { x}) = \frac{3}{2}R T, \quad \quad \rho e = E - \frac{1}{2}\rho |U|^2. 
\label{internal_energy}
\end{equation}
The relaxation time $\tau=\tau(x,t)$ and the mean free path $\lambda$ are related according to  \cite{CC}
\begin{equation}
\tau = \frac{4 \lambda}{\pi \bar C},
\label{tau}
\end{equation}
where $\bar C = \sqrt{{8RT}/{\pi}}$ and the mean free path $\lambda$ is given by 

\begin{equation}
\lambda = \frac{k_B}{\sqrt{2}\pi\rho R d^2},
\label{lambda}
\end{equation}
where $k_B$ is the Boltzmann constant and $d$ is the diameter of the gas molecules.

\subsection{Newton-Euler equations for rigid body motion}
The  motion of a rigid body $S\subset \mathbb{R}^{d_x}$ is given by the Newton-Euler equations, compare \cite{TKR},
\begin{equation}
M\frac{d { V} }{dt} = {\mathcal{ { F} } }, \; \;
  [I]\cdot \frac{d  { \omega} }  {dt} + {{ \omega} }\times ({[I]\cdot}{ { \omega} }) =  {\mathcal{ { T} } },
  \label{euler_newton}
  \end{equation}
where $M$ is the total mass of the body with center of mass ${ X}_c$, ${ V} $ is the velocity of the center of mass ${ X}_c$ and $ {\omega}$ is the angular velocity of the rigid body. $\mathcal{ { F} }$ is the translation force, $ \mathcal{ { T} } $ is the torque  and $[I]$ is the moment of inertia.  
%We note that,  in $2D$ $\omega = (0, 0, \Omega)$ the nonlinear term in the second equation of (\ref{euler_newton}) vanishes and one obtains a scalar equation for the 
%angular velocity $\Omega$.
The center of mass of the rigid body is obtained by  
\begin{equation}
\frac{d { X}_c}{dt} = { V}.   %, \quad \frac{d{\bf \Theta}}{dt} = {\bf \omega}
\label{pos_velo_rigid}
\end{equation}
Finally, the velocity of a point on the surface of the rigid body is given by   $ {{ U}_w} ={{ V}}+{ { \omega} } \times({ { x} } -{ {X}_c}) ,~{ { x} }\in \partial S$. 
  
The force $\mathcal{ { F} }$ and torque ${  \mathcal{ { T} } }$, that the gas exerts on the  rigid body, 
are
%are computed in DSMC \cite{STKH} according to
%  \begin {equation}
%  \mathcal{ {\bf F } } = \frac{\sum_{i} m  {\bf v} ^{in}_i- \sum_{i}m {\bf v}^{out}_i}  {\Delta t},\; \;
%   \mbox{and} \;\;   \mathcal{ T} = ( {\bf x} - {\bf X}_c)\times \mathcal{ {\bf F} },
%  \label{force_torque_dsmc}
%  \end{equation}
%  where $m$ is the mass of the gas molecules. 
 computed according to 
  \begin{equation}
  \mathcal{ { F }} = \int_{\partial S} (-\varphi \cdot{ { n} }_{s}) dA, \; \;
 \mathcal{ { T }}= \int_{\partial S} ( { x} - {X}_c)\times(-\varphi \cdot{ { n} }_{s}) dA,
 \label{force_torque_bgk}
 \end{equation}
 where $\varphi\in\mathbb{R}^{d_x\times d_x}$ is the stress  tensor and is given by 
 \begin{equation}
 \varphi=\int_{\mathbb{R}^3}( { v} - { U}_w)\otimes( { v} - { U}_w) f(t, { x},  { v}) d { v}.
 \label{pressure_tensor}
 \end{equation}

\subsection{Chu-reduction}
In one and two 
physical space  dimensions $d_x =1,2$ one might consider mathematically a one or two dimensional  velocity space
$d_v = 1,2$,
respectively. 
However,   it is physically correct to consider in these situations  still three  velocity  dimensions.
To resolve the three-dimensional  velocity space numerically  requires unnecessary memory and computational time. In these cases, for the BGK model, 
the  3D velocity  space can be reduced as  suggested by Chu \cite{Chu}. This reduction yields a considerable  savings in memory allocation and  computational time. 
For example, in a  physically one-dimensional situation, in which all variables depend on $x \in \mathbb R$ and $t$ (slab geometry),  the 
velocity space is reduced from three dimensions  to one dimension defining  the following reduced distributions \cite{GR}. Considering  $v  = (v_1,v_2,v_3) \in \mathbb R^3$ we define
\begin{equation}
g_{\rm 1}(t,x,v_1) = \int_{\mathbb R^2} f(t,x,v_1,v_2,v_3)dv_2dv_3, \quad
g_{\rm 2} (t,x,v_1) = \int_{\mathbb R^2} (v_2^2+v_3^2)f(t,x,v_1,v_2,v_3)dv_2dv_3.
\label{reduced_g}
\end{equation}
Multiplying (\ref{bgk_eqn}) by $1$ and $v_2^2+v_2^2$ and integrating with respect to $(v_2,v_3) \in \mathbb R^2$, we obtain the following 
system of two equations 
\begin{equation}
\frac{\partial g_1}{\partial t} + v \frac{\partial g_1}{\partial x} = \frac{1}{\tau}(G_1 - g_1),\quad
\frac{\partial g_2}{\partial t} + v \frac{\partial g_2}{\partial x} = \frac{1}{\tau}(G_2 - g_2),
\label{eqduced_bgk}
\end{equation}
where we denoted $v_1$ by $v$, and  
\begin{equation}
G_1 =  \int_{\mathbb R^2} M dv_2 dv_3 = \frac{\rho}{\sqrt{2RT}} 
\exp\left({-\frac{(v-U)^2}{2RT}}\right),\quad
G_2 =  \int_{\mathbb R^2} (v_2^2+v_3^2) M dv_2 dv_3 = (2RT) G_1.  
\label{Maxw_G}
\end{equation}
Assuming the initial condition is a local equilibrium, the initial distributions are defined via the parameters $(\rho_0, U_0, T_0)\in\mathbb{R}^3$ and are given as 
\begin{eqnarray}
g_1(0,x,v) =  \frac{\rho_0}{\sqrt{2RT_0}} \exp\left({-\frac{(v-U_0)^2}{2RT_0}}\right), \quad 
g_2(0,x,v) = (2RT_0) g_1(0,x,v). 
\end{eqnarray}
The macroscopic quantities are given through the reduced distributions as 
\begin{equation}
\rho = \int_{\mathbb R} g_1dv, \; \rho U = \int_{\mathbb R}v g_1 dv, \; 3\rho RT = \int_{\mathbb R}(v-U)^2 g_1dv + \int_{\mathbb R} g_2 dv.  
\label{macroeq}
\end{equation}

Similarly,  in two spatial dimensions $x \in \mathbb R^2$,  the reduction from a three dimensional to a two dimensional velocity space  is obtained by multiplying the BGK model (\ref{bgk_eqn}) by $1$ and $v_3^2$ and 
integrating wrt $dv_3$ over $\mathbb R$. The reduced equations are two-dimensional versions of (\ref{eqduced_bgk})
with  $v=(v_1,v_2) \in \mathbb R^2$, but the reduced Maxwellians $G_1$ and $G_2$ are 
given as 
\begin{equation}
G_1 = \frac{\rho}{2RT}\exp\left({-\frac{|v-U|^2}{2RT}}\right), \quad G_2 = (RT) G_1
\end{equation}
 with $U= (U_1,U_2) \in \mathbb R^2 $. The distribution functions are 
\[
g_1(t,x,v_1,v_2) = \int_{\mathbb R} f(t,x,v_1,v_2,v_3)dv_3, \quad
 g_2(t,x,v_1,v_2)= \int_{\mathbb R} v_3^2 f(t,x,v_1,v_2,v_3)dv_3.
\]

 %%%%%%%%%%%%%%%%%%%%%%%%%%%%%%
\section{Numerical schemes}

We solve the original  equation (\ref{bgk_eqn}) and the  reduced system of equations (\ref{eqduced_bgk}) by the ALE method described below. We use a time splitting, where the 
advection step is solved explicitly and the relaxation part is solved  implicitly.
%\sout{and a meshfree particle method to solve this system of equations.}  
%\sout{Here, by particle we actually mean} 
%\giovanni{I was not satisfied by the sentence, since these are not particles, since they are not in phase space, but only in physical space, and their velocity grid does not change. I correct as follows} 
Using  a discrete velocity approximation of the distribution function (see Section \ref{sec:velocity_discretization}) the information is stored on grid points in physical space  moving  with the mean velocity $U$ of the gas. The spatial derivatives 
of the distribution function  at an arbitrary particle position are approximated using values at  the point-cloud surrounding the particle and   a weighted least squares method. 

In the following. we present first and second order schemes in time as well as in space. 

\subsection{ALE formulation}
We consider original and reduced model.

\subsubsection{ALE formulation for the original model}
We rewrite  the equations (\ref{bgk_eqn}) in  Lagrangian form as 
\begin{eqnarray}
\frac{dx}{dt} &=& U \label{originalALE1}
\\
\frac{d f}{d t} &=& -(v-U) \cdot \nabla_x f+ \frac{1}{\tau}(M-f)
\label{originalALE2}
\end{eqnarray}
where ${d }/{dt} = {\partial }/{\partial t} + U\cdot \nabla_x$. The first equation describes motion  with the macroscopic 
mean velocity $U$ of the gas determined by (\ref{macroeq}). The second equation includes the remaining advection with the difference between microscopic and macroscopic velocity. 

\subsubsection{ALE for reduced  model}
In this case  the equations (\ref{reduced_g}) are reformulated in  Lagrangian form as 
\begin{eqnarray}
\frac{dx}{dt} &=& U \label{ALE1}
\\
\frac{d g_1}{d t} &=& -(v-U) \cdot \nabla_x  g_1+ \frac{1}{\tau}(G_1 - g_1)
\label{ALE2}
\\
\frac{d g_2}{d t} &=& - (v-U) \cdot \nabla_x g_2 + \frac{1}{\tau}(G_2 - g_2).
\label{ALE3}
\end{eqnarray}

\subsection{Time discretization}

\subsubsection{First order time splitting scheme for the original model}
\label{1st_order_step}
Time is discretized as  $t^n = n \Delta t, n= 0, 1, \cdots, N_t$. We denote the numerical approximation of $f$ at $t_n$ by $f^n = f(t^n, x, v)$. 
We use a time splitting scheme for equation (\ref{originalALE2}), where the advection term is solved explicitly and the collision term is solved implicitly. 
In the first step of the splitting scheme we obtain the intermediate distribution $\tilde{f}^n$ by solving
\begin{eqnarray}
\tilde{f}^n &=& f^n - \Delta t (v-U^n) \cdot \nabla_x f^n .
\end{eqnarray}
In the second step we obtain the new distribution by solving 
\begin{eqnarray}
f^{n+1} &=& \tilde{f}^n+ \frac{\Delta t}{\tau}(M^{n+1} - f^{n+1} )
\label{implicitfull_1}
\end{eqnarray}
and the new positions of the grids are updated by 
\begin{equation}
x^{n+1} = x^n + \Delta t U^n. 
\end{equation}

In the first step, we have to approximate the spatial derivatives of $f$ at every grid point. This is described in the following section.

Following \cite{GRS,GR,XRQ} 
we  obtain $f^{n+1}$ in the second step by  first  determining  the parameters $\rho^{n+1}, U^{n+1}$ and $T^{n+1}$ for $M^{n+1}$. 
Multiplying (\ref{implicitfull_1}) by $1$, $v$ and $(v-U)^2$ and  integrating  over velocity space,  we get 
\begin{equation}
\rho^{n+1} = \int_{\mathbb R} \tilde{f}^n dv, \quad (\rho U)^{n+1} = \int_{\mathbb R} v\tilde{f}^n dv, \quad 
3\rho RT^{n+1} = \int_{\mathbb R^3} |v-U|^2 \tilde{f}^n dv 
\label{fullmoment1a}
\end{equation}
where we have used the conservation of mass, momentum and energy of the original BGK model.

	Now, the parameters $\rho^{n+1}, U^{n+1}$ and $T^{n+1}$ of $M^{n+1}$  are given in terms of $\tilde{f}$  from (\ref{fullmoment1a}) . Hence the implicit step (\ref{implicitfull_1}) can be explicitly solved  as 
	\begin{eqnarray}
	f^{n+1} = \frac{\tau \tilde{f}^n + \Delta t M^{n+1}}{\tau + \Delta t}.
	\end{eqnarray}
	
	\subsubsection{Second order splitting scheme for the original model (ARS(2,2,2))}
\label{ARS}
For the second order  splitting scheme we use the stiffly accurate ARS(2,2,2) scheme, \cite{ARS}, and compare the results with a slightly simpler scheme ARS(2,2,1). The Butcher tabueau of both schemes are reported below, in the usual form expressed in Table \ref{tab:Butcher}

\begin{table}
\begin{center}
\begin{tabular} { l | c }
$c$ & $A$ \\[1mm] \hline\\[-3mm]
& $b^\intercal$
\end{tabular} \quad
\begin{tabular} { l | c }
$\tilde{c}$ & $\tilde{A}$ \\[1mm] \hline \\[-3mm]
& $\tilde{b}^\intercal$ 
\end{tabular}
\end{center}
\caption{Classical form of the double Butcher tableau of an IMEX scheme: matrix $A$ and vectors $b$ and $c$ are relative to the implicit scheme, while $\tilde{A}$, $\tilde{b}$, $\tilde{c}$ denote the RK coefficients of the implicit scheme.}
\label{tab:Butcher}
\end{table}

\begin{table}
\begin{center}
\[\quad\quad
\begin{array} { c | ccc }
0      & 0       & 0       & 0 \\ 
\beta  & 0       & \beta   & 0 \\ 
1      & 0       & 1-\beta & \beta \\ 
\hline
       & 0       & 1-\beta & \beta \\ 
\end{array}
\quad
\begin{array} { c | ccc }
0      & 0       & 0       & 0 \\ 
\beta  & \beta   & 0       & 0 \\ 
1      & \beta-1 & 2-\beta & 0 \\ 
\hline
       & \beta-1 & 2-\beta & 0 \\ 
\end{array}
\quad
\quad
\quad
\begin{array} { c | ccc }
0      & 0       & 0    & 0 \\ 
1/2    & 0       & 1/2  & 0 \\ 
1      & 0       & 0    & 1 \\ 
\hline
       & 0       & 0    & 1 \\ 
\end{array}
\quad
\begin{array} { c | ccc }
0      & 0       & 0       & 0 \\ 
1/2    & 1/2     & 0       & 0 \\ 
1      & 0       & 1       & 0 \\ 
\hline
       & 0       & 1       & 0 \\ 
\end{array}
\]
\end{center}
\caption{Tableau of ARS(2,2,2) scheme (left) and of ARS(2,2,1) scheme (right). $\beta = 1-1/\sqrt{2}$.} 
\label{tab:ARS(2,2,2)}
\end{table}

For  equation (\ref{originalALE1}-\ref{originalALE2}) this leads to

{\bf Step 1:}\\
\begin{eqnarray}
x^{n+\frac{1}{2}} &= &x^n + \beta \Delta t U^n. \\
\tilde{f}^{n+\frac{1}{2}} &=& f^n - \beta \Delta t (v-U^n)   \cdot \nabla_x  f^n.
\end{eqnarray}
The intermediate distributions $f^{n+\frac{1}{2}}$ are then obtained by solving 
\begin{eqnarray*}
	f^{n+\frac{1}{2}} &=& \tilde{f}^{n+\frac{1}{2}}+ \beta  \frac{\Delta t}{\tau}(M^{n+\frac{1}{2}} - f^{n+\frac{1}{2}} )
\end{eqnarray*}
or
\begin{eqnarray}
f^{n+\frac{1}{2}} = \frac{  \tau \tilde{f}^{n+\frac{1}{2}} + \beta \Delta t M^{n+\frac{1}{2}}}{  \tau + \beta  \Delta t}.
\end{eqnarray}

{\bf Step 2:}\\
\begin{eqnarray}
x^{n+1} &= &x^n + \Delta t ((\beta-1)U^n + (2-\beta)U^{n+\frac{1}{2}}). \\
\tilde{f}^{n+1 }&=& f^n - (2-\beta)\Delta t (v-U^{n+\frac{1}{2}})  \cdot \nabla_x   f^{n+\frac{1}{2}} \\
&&- (\beta-1)\Delta t (v-U^{n})  \cdot \nabla_x   f^{n} \nonumber\\
&&+ (1-\beta) \frac{\Delta t}{\tau}(M^{n+\frac{1}{2}} - f^{n+\frac{1}{2}} ).\nonumber
\end{eqnarray}
The new distributions are obtained by solving 
\begin{eqnarray*}
	f^{n+1} &=& \tilde{f}^{n+1} + \beta \frac{\Delta t}{\tau}(M^{n+1} - f^{n+1} ).
\end{eqnarray*}
or
\begin{eqnarray}
f^{n+1} = \frac{\tau \tilde{f}^{n+1} + \beta \Delta t M^{n+1}}{\tau + \beta \Delta t}.
\end{eqnarray}

with $\beta = 1-1/{\sqrt{2}}$.
We note that the implicit computations of $M^{n+\frac{1}{2}}$ and of $M^n$ 
are similar to the implicit computation of $M^{n+1}$ 
in the first order scheme as described above.

	\subsubsection{Partial second order time splitting scheme for the original model ARS(2,2,1))}
	\label{RK2}
For later use we also describe a simplified scheme with an explicit second order solution of the advection equation  and an implicit first order solution of the collision term. 
	For the second order  scheme we use a two step Runge-Kutta scheme. For  equation (\ref{originalALE1}-\ref{originalALE2}) the scheme  is given by 
	\\
	{\bf Step 1:}\\
	\begin{eqnarray}
	x^{n+\frac{1}{2}} &= &x^n + \frac{\Delta t}{2} U^n. \\
	\tilde{f}^{n+\frac{1}{2}} &=& f^n - \frac{\Delta t}{2} (v-U^n)   \cdot \nabla_x  f^n.
	\end{eqnarray}
	The intermediate distributions $f^{n+\frac{1}{2}}$ are then obtained by solving 
	\begin{eqnarray*}
	f^{n+\frac{1}{2}} &=& \tilde{f}^{n+\frac{1}{2}}+ \frac{1}{2} \frac{\Delta t}{\tau}(M^{n+\frac{1}{2}} - f^{n+\frac{1}{2}} )
	\end{eqnarray*}
	i.e.
	\begin{eqnarray}
	f^{n+\frac{1}{2}} = \frac{ 2 \tau \tilde{f}^{n+\frac{1}{2}} + \Delta t M^{n+\frac{1}{2}}}{ 2 \tau + \Delta t}.
	\end{eqnarray}
	
	{\bf Step 2:}\\
	\begin{eqnarray}
	x^{n+1} &= &x^n + \Delta t U^{n+\frac{1}{2}}. \\
	\tilde{f}^{n+1 }&=& f^n - \Delta t (v-U^{n+\frac{1}{2}})  \cdot \nabla_x   f^{n+\frac{1}{2}} .
	\end{eqnarray}
	The new distributions are obtained by solving 
	\begin{eqnarray}
	f^{n+1} &=& \tilde{f}^{n+1} + \frac{\Delta t}{\tau}(M^{n+1} - f^{n+1} ).
	\end{eqnarray}
	or
	\begin{eqnarray}
	f^{n+1} = \frac{\tau \tilde{f}^{n+1} + \Delta t M^{n+1}}{\tau + \Delta t}.
	\end{eqnarray}
	
\begin{rem}
Note that this scheme is not the 
Midpoint rule, which is A-stable, but not L-stable.  It is not second order, but it is L-stable, therefore it can be adopted with arbitrarily small values of the relaxation time $\tau$.
However, the scheme is simpler and less costly than the ARS scheme and  in the examples considered here, we obtain numerically 
 second  order of convergence. 
\end{rem}

\subsubsection{Time splitting scheme for the reduced model}
\label{1st_order_step_red}

We use again a  time splitting scheme. For the first order scheme with one-dimensional  physical space $x \in \mathbb R$,
we proceed as follows.
In the first step we obtain the intermediate distributions $\tilde{g}^n_1$  and $\tilde{g}^n_2$ by solving
for $v \in \mathbb R$ and $U \in \mathbb R$
\[
\tilde{g}_1^n = g_1^n - \Delta t (v-U^n) \partial_x g_1^n\quad
\tilde{g}_2^n = g_2^n - \Delta t (v-U^n) \partial_x g_2^n .
\]
In the second step we obtain the new distributions by solving 
\begin{eqnarray}
g_1^{n+1} &=& \tilde{g}_1^n + \frac{\Delta t}{\tau}(G_1^{n+1} - g_1^{n+1} )
\label{implicit_1}
\\
g_2^{n+1} &=& \tilde{g}_2^n + \frac{\Delta t}{\tau}(G_2^{n+1} - g_2^{n+1} )
\label{implicit_2}
\end{eqnarray}
and the new positions of the grids are updated by 
\begin{equation}
x^{n+1} = x^n + \Delta t U^n. 
\end{equation}

For  the second step  we have to determine first the parameters $\rho^{n+1}, U^{n+1}$ and $T^{n+1}$ for $G_1^{n+1}$ and $G_2^{n+1}$. 
Multiplying (\ref{implicit_1}) by $1$ and $v$ and integrating with respect to $v$ over $\mathbb R$ we get 
\begin{equation}
\rho^{n+1} = \int_{\mathbb R} \tilde{g}_1^n dv, \quad (\rho U)^{n+1} = \int_{\mathbb R} v\tilde{g}_1^n dv, 
\label{moment1a}
\end{equation}
where we have used  the conservation of mass and momentum of the original BGK model.  In order to compute $T^{n+1}$ 
we note that the following identity is valid
\begin{equation}
\int_{\mathbb R} (v-U)^2(G_1-g_1)dv + \int_{\mathbb R} (G_2 - g_2) dv = 0. 
\label{id3}
\end{equation}

Multiplying the equation (\ref{implicit_1}) by $(v-U)^2$ and integrate with respect to $v$ over $\mathbb R$ we get 
\begin{equation}
\int_{\mathbb R} (v-U)^2 g_1^{n+1} dv = \int_{\mathbb R} (v-U)^2 \tilde{g}_1^n dv + \frac{\Delta t}{\tau} \int_{\mathbb R} (v-U)^2(G_1^{n+1} - g_1^{n+1} ) dv.
\label{id4}
\end{equation}

Next,  integrate both sides of {\ref{implicit_2}) with respect to $v$ over $\mathbb R$ we get 
\begin{equation}
\int_{\mathbb R} g_2^{n+1}dv = \int_{\mathbb R}\tilde{g}_2^n dv + \frac{\Delta t}{\tau} \int_{\mathbb R} (G_2^{n+1} - g_2^{n+1}) dv. 
\label{id5}
\end{equation}

Adding (\ref{id4}) and (\ref{id5}) and making use of the identity (\ref{id3}) we get 
\begin{equation}
3 \rho^{n+1} R T^{n+1} = \int_{\mathbb R} (v-U)^2 \tilde{g}_1^n dv + \int_{\mathbb R} \tilde{g}_2^n dv. 
\label{moment1b}
\end{equation} 

Now, the parameters $\rho^{n+1}, U^{n+1}$ and $T^{n+1}$ of $G_1^{n+1}$ and $G_2^{n+1}$ are given in terms of $\tilde{g}_1^n$ and $\tilde{g}_2^n$ from (\ref{moment1a}) and 
(\ref{moment1b}). Hence the implicit steps (\ref{implicit_1}) and (\ref{implicit_2}) can be rewritten as 
\begin{eqnarray}
g_1^{n+1} = \frac{\tau \tilde{g}_1^n + \Delta t G_1^{n+1}}{\tau + \Delta t}\\
g_2^{n+1} = \frac{\tau \tilde{g}_2^n + \Delta t G_2^{n+1}}{\tau + \Delta t}.
\end{eqnarray}

The second order time splitting for the reduced model follows the lines of the second order splitting procedure for the original model.

\subsection{Velocity discretization} 
\label{sec:velocity_discretization}
For the sake of simplicity we consider a one-dimensional velocity domain. Consider $N_v+1$ velocity grid points and a uniform velocity grid of size $\Delta v = 2 v_{\rm max}/ N_v$ 
We assume that the distribution function is negligible for $|v| > v_{\rm max} $ and  discretize $[-v_{\rm max},v_{\rm max}]$.  
That means for each velocity direction  we have the discretization points 
$v_j = -v_{\rm max} + (j-1)\Delta {v}, j = 1,\ldots , N_v+1$.  
Note that the performance of the method could be improved by using a grid adapted to the mean velocity $U$, see, for example,
\cite{DM}.

\subsection{Spatial discretization}
 
 We discuss the spatial discretization and upwinding procedures for first and second order schemes.
 
\subsubsection{Approximation of spatial derivatives}

In the above  numerical schemes an approximation of  the spatial derivatives of $g_1$ and $g_2$ is required. In this subsection, we describe a least squares approximation of the derivatives on the moving point cloud based on  so called generalized finite differences, see \cite{LO, SKT} and  references there in.  A stabilizing procedure using  upwinding and a WENO
type discretization for the higher order schemes  will 
be described in the following. 

For the sake of simplicity we consider a one-dimensional spatial domain $\Omega$. We first approximate 
the boundary of the domain by a set of discrete points called  boundary particles.   In the second step we approximate the 
interior of the computational domain using another set of interior points or interior particles. The sum of boundary and interior points gives the total number of points. 
We note that the boundary conditions are applied on the boundary points. The boundary points move together with the boundaries. The initial generation of grid points 
can be regular as well as arbitrary.  When the points 
move they can form a cluster or can scatter away from each other. In these cases, either some grid points have to be removed or new grid points have to be added. We will describe this particle management in the next subsection. 

Let  ${x}_i  \in \Omega, i = 1 , \ldots , N_x$, where $N_x$ is the total number of grid points with initial average spacing $\Delta x$.  Let $f(x)$ be a scalar function and
$f_i$  its discrete values in $x_i$. 
 Our main task is to approximate the spatial derivatives of $f_i$ at an arbitrary position $x_i$ from its neighboring particles.  We call $x_i$ a central point. 
 We sort the neighboring points 
 into different catagories, left, right and central neighbor. Note that the point $x_i$ is itself its neighbor in all sets of neighboring particles.  We restrict to  neighboring  points 
 within a radius $h$ in such a way that we have at least a minimum number of neighbors.  $h$ is usually chosen in  relation to $\Delta x$, compare \cite{TK07}. For a  first order approximation one can choose smaller values of $h$ than for a higher order approximation.  
 In order to guarantee a better accuracy we associate a weight function depending on  the 
 distance of the central point and its neighbors.   
Let $P(x) = x_{k},k=1, \ldots, m $ be the set of $m$ neighbor points of $x=x_i$ inside the radius $h$.  There are several choices of weight functions \cite{Sonar}. We choose a Gaussian weight function \cite{TK07, TKH09}
\begin{eqnarray} 
 w( {x}_k - x; h) =
\left\{ 
\begin{array}{l}  
\exp \left(- \alpha \frac{(x_k - x)^2}{h^2} \right), 
%\quad \mbox{if    }  \frac{ |x_k - x|} {h} \le 1 
\\ \nonumber
 0,  \qquad \qquad \mbox{else},
\end{array}
\right.
\label{weight}
\end{eqnarray}
with $ \alpha $ a  user defined positive constant. In our computation, we have chosen $\alpha = 6$. 

In order to approximate the derivatives we consider a second  order  Taylor expansion of $f(x_{k})$ around $x$ 
\begin{equation}
f(x_{k}) = f(x) + (x_{k}-x)\cdot \partial_x  f (x)+ (x_{k}-x)^T \partial_{xx } f (x) (x_{k}-x)   +e_{k}, 
\label{taylor}
\end{equation}
for $k=1, \ldots, m$, where $e_{k}$ is the error in the Taylor's expansion.   
The unknown $a = [\frac{\partial  f}{\partial x} (x), \frac{\partial^2  f}{\partial x^2}(x)]^T$ is now computed by minimizing the error $e_{k}$ for 
$k=1, \ldots, m$.  
The system of equations can be re-written in  vector form as 
\begin{equation}
\label{sys}
{ e } = { b} - D a, 
\end{equation}
where ${e} = [e_{1}, \ldots, e_{m}]^T$,  $b = [f_{1} - f(x), \ldots, f_{m} - f(x)]^T $ 
and
\begin{eqnarray}
D = \left( \begin{array}{cc}
dx_1 &  ~\frac{1}{2}dx^2_1    \\
\vdots  &\vdots    \\
dx_m &  ~\frac{1}{2}dx^2_m   
 \end{array} \right)
\label{matrixM}
\end{eqnarray}
with $dx_k = x_{k} - x$. %, \;  dy_j = y_{i_k}-y_i$. 

%\sout{In general, this system of equations is over-determined.} \giovanni{there is on system of equations so far} 
Imposing $e=0$ in (\ref{sys})   results in an overdetermined linear stems of algebraic equations, which in general has no solution. The unknown
$a$ is therefore obtained from the weighted least squares method by minimizing the quadratic form 
\begin{equation}
J = \sum_{k=1}^m w_{k} e_{k}^2 = (D a  - { b})^T W (D a - { b}),
\end{equation}
where $W=\mathrm diag(w_{1}, \ldots, w_m) $. % is the diagonal matrix. 
The minimization of $J$ formally yields 
\begin{equation}
a = (D^TWD)^{-1}(D^TW){ b}. 
\end{equation}

\subsubsection{First order upwind scheme}
We describe the procedure for simplicity only for one-dimensional physical space. We compute the partial derivatives of $g_1$ and $g_2$ in the following way. If $v-U > 0$, we compute the derivatives 
at $x_i$ from the set of left neighbors $P_L(x_i)$ lying within the radius $h$. Similarly, for $v-U < 0$ we use the set of right neighbors $P_R(x_i)$  lying within the radius $h$. 
Then we use the Taylor expansion (\ref{taylor})  to first order  and compute the derivatives in the corresponding set of neighboring  points.

\subsubsection{Second order WENO-type procedure}

When we apply a second order Taylor expansion, the scheme becomes unstable if the solution develops discontinuities. In this case we use the WENO idea 
in order to obtain higher order derivatives. 
We refer to   \cite{ADB,AVESANI2021113871,Z19} for  similar approaches for SPH-type particle methods.
For the sake of simplicity, we consider the one dimensional case to present our simplified WENO procedure. Let $P_L(x_i), P_R(x_i)$ and $P_C(x_i)$ be the sets 
of left, right and central neighbor points, see Fig. \ref{1d_grid}. Note that $P_C(x_i) = P_L(x_i) \cup P_R(x_i)$. 

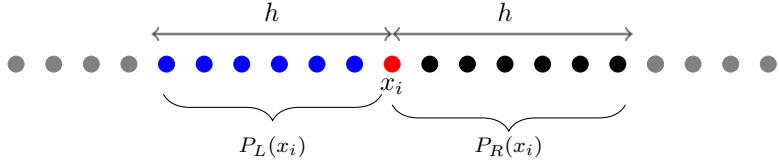
\begin{figure}[!t]
 \begin{center}
 \begin{tikzpicture}
               \filldraw[gray] (2,0.5) circle (3pt);
               \filldraw[gray] (2.5,0.5) circle (3pt);
               \filldraw[gray] (3,0.5) circle (3pt);
               \filldraw[gray] (3.5,0.5) circle (3pt);
                \filldraw[blue] (4,0.5) circle (3pt);
               \filldraw[blue] (4.5,0.5) circle (3pt);  
               \filldraw[blue] (5,0.5) circle (3pt);
                \filldraw[blue] (5.5,0.5) circle (3pt); 
               \filldraw[blue] (6,0.5) circle (3pt);  
                \filldraw[blue] (6.5,0.5) circle (3pt);
                 \filldraw[red] (7,0.5) circle (3pt);
                \filldraw[black] (7.5,0.5) circle (3pt);
                 \filldraw[black] (8.0,0.5) circle (3pt);
                  \filldraw[black] (8.5,0.5) circle (3pt);
                   \filldraw[black] (9.0,0.5) circle (3pt);
                    \filldraw[black] (9.5,0.5) circle (3pt);
                     \filldraw[black] (10.0,0.5) circle (3pt);
                      \filldraw[gray] (10.5,0.5) circle (3pt);
                       \filldraw[gray] (10.5,0.5) circle (3pt);
                        \filldraw[gray] (11.0,0.5) circle (3pt);
                         \filldraw[gray] (11.5,0.5) circle (3pt);
                          \filldraw[gray] (12.0,0.5) circle (3pt);
%                \draw (4,0.5) node[circle,fill=blue]{} circle(0.0001cm);
%               \draw (4.5,0.5) node[circle,fill=blue]{} circle(0.0001cm);
%               \draw (5,0.5) node[circle,fill=blue]{} circle(0.0001cm);
%                \draw (5.5,0.5) node[circle,fill=blue]{} circle(0.0001cm);
%                \draw (6,0.5) node[circle,fill=blue]{} circle(0.0001cm);
%                \draw (6.5,0.5) node[circle,fill=blue]{} circle(0.0001cm);
 %               \draw (7,0.5) node[circle,fill=red]{} circle(0.0001cm);
%                \draw (7.5,0.5) node[circle,fill=black]{} circle(0.0001cm);
%                \draw (8,0.5) node[circle,fill=black]{} circle(0.0001cm);
%                \draw (8.5,0.5) node[circle,fill=black]{} circle(0.0001cm);
%                \draw (9,0.5) node[circle,fill=black]{} circle(0.0001cm);
%                \draw (9.5,0.5) node[circle,fill=black]{} circle(0.0001cm);
%                \draw (10,0.5) node[circle,fill=black]{} circle(0.0001cm);
%                 \draw (10.5,0.5) node[circle,fill=gray]{} circle(0.0001cm);
%                  \draw (11,0.5) node[circle,fill=gray]{} circle(0.0001cm);
%                   \draw (11.5,0.5) node[circle,fill=gray]{} circle(0.0001cm);
%                    \draw (12,0.5) node[circle,fill=gray]{} circle(0.0001cm);
                  %    \draw (7,0.5) node[circle,fill=red]{} circle(3.2cm);
                     \draw [<->] [line width=1pt,color=black,opacity=0.5] (7,0.9) -- (10.2,0.9) ;
                      \draw[<->] [line width=1pt,color=black,opacity=0.5] (3.8,0.9) -- (7.0,0.9) ;
                      \node[color=black] at (7.0,0.2) {$x_i$ };
                        \node[color=black] at (8.5,1.2) {$h$ };
                          \node[color=black] at (5.4,1.2) {$h$ };
                        %  \node[color=black] at (5.4,0.1) {$P_L(x_i)$ };
                          % \node[color=black] at (8.5,1) {$P_R(x_i)$ };
                           \draw [decorate,decoration={brace,amplitude=10pt},xshift=-4pt,yshift=0pt]
(7.0,0.1) -- (4.1,0.1) node [black,midway,xshift=-0.0cm,yshift=-0.7cm] 
{\footnotesize $P_L(x_i)$};
\draw [decorate,decoration={brace,amplitude=10pt,mirror,raise=4pt},yshift=0pt]
(7.0,0.1) -- (10.1,0.1) node [black,midway,xshift=0.0cm,yshift=-0.7cm] {\footnotesize
$P_R(x_i)$};  
 \end{tikzpicture}
 \caption{Central, left and right neighbor points.}
 \label{1d_grid}
 \end{center}
 \end{figure}

Considering the Taylor expansions (\ref{taylor}) and applying the least squares method, we obtain the 
derivatives 
\[
{f_x}_L, {f_{xx}}_L, {f_x}_R, {f_{xx}}_R, {f_x}_C, {f_{xx}}_C
\]
using left, central and right neighbors, respectively. 
The desired first order derivative is obtained by the weighted sum 
\begin{eqnarray}
f_x &=& \omega_L {f_x}_L + \omega_C {f_x}_C + \omega_R {f_x}_R, 
\end{eqnarray}
where the weights are defined by 
\begin{equation}
\omega_k = \frac{\beta_k}{\beta_L + \beta_C + \beta_R}, \quad k = L, C, R
\end{equation}
with
\begin{equation}
\beta_k = \frac{C_k}{({f_x}_k^2 \Delta x^2 + {f_{xx}}_k^2 \Delta x^4  + \epsilon)^2}, \quad k = L, C, R
\end{equation}	
where $\epsilon = 10^{-6}$ and $\Delta x$ is the initial spacing of particles.  This is combined with the following choice of the coefficients $C_k$
depending on the sign of  $v - U$. 
If $v - U>0$ the values  are 
\begin{center}
\[
C_L = 0.5, C_C = 0.5, C_R = 0 
\]
\end{center}
and otherwise
\begin{center}
\[
C_L = 0, C_C = 0.5, C_R = 0.5. 
\]
\end{center}

In 2D we proceed in an analogous way. Here the derivatives $f_x$ and $f_y$ are required. They are obtained by determining the sets of points in the  left (L) and  right (R) half plane for the determination of $f_x$ and the sets in the top (T) and bottom (B) half plane for the determination of $f_y$, see Fig. \ref{2d_grids}. 
%In Fig. \ref{2d_grids} we have distinguish the 4 sets of neighbor lists of central points, which is the red one, with in the radius $h$. The central red point is common for all sets of neighbor. 

 \begin{figure}[ht]
 	\begin{center}
 		\begin{tikzpicture}[scale=0.7]   
		 %  \draw [-] [line width=2pt,color=black,opacity=0.5] (-5,0) -- (5,0) ;
		    \draw [-] [line width=2pt,color=black,opacity=0.5] (-6,-5) -- (-6,5);
		   \draw [-] [line width=0.8pt,color=black,opacity=0.5] (-6,0) -- (-3.6,2.1);
		   \node[color=black] at (-5,0.6) { h };
		     \node[color=blue] at (-9.5,0) {\bf L };
		      \node[color=blue] at (-2.5,0) {\bf R };
		     \draw (-6,0) node[]{} circle(3.1cm);
	          \filldraw[red] (-6,0) circle (3pt);
	          \filldraw[black] (-5,0.2) circle (3pt);	 
	          \filldraw[black] (-4,-0.2) circle (3pt);         
	          \filldraw[black] (-3,0.3) circle (3pt);	 
	         \filldraw[blue] (-8,-0.2) circle (3pt);	 
	          \filldraw[blue] (-9,0.2) circle (3pt);	 
	          \filldraw[blue] (-7,-0.3) circle (3pt);      
	         % \filldraw[blue] (-6,0) circle (4pt);	     
	          \filldraw[blue] (-6.3,1) circle (3pt);	 
	          \filldraw[black] (-5.5,2) circle (3pt);       
	          \filldraw[blue] (-6.4,1.7) circle (3pt);	 
	          \filldraw[blue] (-6.3,-1) circle (3pt);	 
	          \filldraw[black] (-5.8,-1.5) circle (3pt);      	
	          \filldraw[blue] (-6.9,-2) circle (3pt);	 
                  \filldraw[blue] (-7,.5) circle (3pt);	     
                  \filldraw[blue] (-7.5,1) circle (3pt);
                   \filldraw[blue] (-8,2) circle (3pt);
                   \filldraw[blue] (-8,.8) circle (3pt);	
                   \filldraw[blue] (-7,2.1) circle (3pt);	 
                    \filldraw[blue] (-6.9,1.2) circle (3pt);	   
                     \filldraw[blue] (-8.5,1.3) circle (3pt);
                     \filldraw[blue] (-6.5,2.7) circle (3pt);   
                   \filldraw[blue] (-7,-0.8) circle (3pt);    
                   \filldraw[blue] (-7.6,-1.6) circle (3pt);     
                   \filldraw[blue] (-8,-0.6) circle (3pt);    
                   \filldraw[blue] (-8.5,-1.4) circle (3pt);  
                   \filldraw[blue] (-6.5,-1.6) circle (3pt);       	 
                   \filldraw[blue] (-7.5,-2.3) circle (3pt);      
                    \filldraw[blue] (-6.5,-2.8) circle (3pt);    
                  \filldraw[black] (-5.5,1) circle (3pt);        
                   \filldraw[black] (-5,1.4) circle (3pt);       
                    \filldraw[black] (-4.5,1) circle (3pt);      
                     \filldraw[black] (-5.6,2.8) circle (3pt);     
                     \filldraw[black] (-5,2.3) circle (3pt);      
                      \filldraw[black] (-4.5,2) circle (3pt);     
                    \filldraw[black] (-5.4,-0.6) circle (3pt);     
                     \filldraw[black] (-5.5,-2.75) circle (3pt);       
                      \filldraw[black] (-5,-1.3) circle (3pt);   
                      \filldraw[black] (-5.6,-2) circle (3pt);    
                      \filldraw[black] (-4.8,-2.1) circle (3pt);        
                      \filldraw[black] (-4.7,-0.7) circle (3pt);     
                      \filldraw[black] (-3.8,-1) circle (3pt);     
                      \filldraw[black] (-3.8,-1.9) circle (3pt);   
                     \filldraw[black] (-3.6,1.1) circle (3pt);   
                       \node[color=black] at (-5.7,-0.3) {$x_i$ };
%%%%%%%%%%%%%%%%%%%%%%%%%%%%%%%%%%%
		   \draw [-] [line width=2pt,color=black,opacity=0.5] (-1.5,0) -- (8,0) ;
		\draw [-] [line width=0.8pt,color=black,opacity=0.5] (3,0) -- (4.9,2.4);
		   \node[color=black] at (4,0.6) { h };
		     \node[color=blue] at (3,3.5) {\bf T };
		      \node[color=blue] at (3,-3.5) {\bf B };
		     \draw (3,0) node[]{} circle(3.1cm);
	          \filldraw[red] (3,0.0) circle (3pt);
	          \filldraw[black] (4,-0.2) circle (3pt);	 
	          \filldraw[blue] (5,0.3) circle (3pt);         
	          \filldraw[blue] (6,0.3) circle (3pt);	 
	         \filldraw[black] (1,-0.2) circle (3pt);	 
	          \filldraw[blue] (0,0.2) circle (3pt);	 
	          \filldraw[black] (2,-0.3) circle (3pt);      
	         % \filldraw[blue] (-6,0) circle (4pt);	     
	          \filldraw[blue] (2.7,1) circle (3pt);	 
	          \filldraw[blue] (3.4,1.5) circle (3pt);       
	          \filldraw[blue] (2.7,1.8) circle (3pt);	 
	          \filldraw[black] (3,-1) circle (3pt);	 
	          \filldraw[black] (3,-2) circle (3pt);      	
	          \filldraw[black] (3,-3) circle (3pt);	 
                  \filldraw[blue] (2,.5) circle (3pt);	     
                  \filldraw[blue] (1.5,1) circle (3pt);
                   \filldraw[blue] (1,2) circle (3pt);
                   \filldraw[blue] (1,.8) circle (3pt);	
                   \filldraw[blue] (2,2.1) circle (3pt);	 
                    \filldraw[blue] (2.1,1.2) circle (3pt);	   
                     \filldraw[blue] (0.5,1.3) circle (3pt);
                     \filldraw[blue] (2.5,2.7) circle (3pt);   
                   \filldraw[black] (2,-0.8) circle (3pt);    
                   \filldraw[black] (1.4,-1.6) circle (3pt);     
                   \filldraw[black] (1,-0.6) circle (3pt);    
                   \filldraw[black] (0.5,-1.4) circle (3pt);  
                   \filldraw[black] (2.5,-1.6) circle (3pt);       	 
                   \filldraw[black] (1.5,-2.3) circle (3pt);      
                    \filldraw[black] (2.5,-2.8) circle (3pt);    
                  \filldraw[blue] (4.5,1) circle (3pt);        
                   \filldraw[blue] (4,1.4) circle (3pt);       
                    \filldraw[blue] (4.5,1) circle (3pt);      
                     \filldraw[blue] (3.4,2.4) circle (3pt);     
                     \filldraw[blue] (4,2.3) circle (3pt);      
                      \filldraw[blue] (4.5,2) circle (3pt);     
                    \filldraw[black] (4.6,-0.6) circle (3pt);     
                     \filldraw[black] (3.5,-2.75) circle (3pt);       
                      \filldraw[black] (4,-1.3) circle (3pt);   
                      \filldraw[black] (4.4,-2) circle (3pt);    
                      \filldraw[black] (4,-2.4) circle (3pt);        
                      \filldraw[black] (4,-0.7) circle (3pt);     
                      \filldraw[black] (5.2,-1) circle (3pt);     
                      \filldraw[black] (5.2,-1.9) circle (3pt);   
                     \filldraw[blue] (5.4,1.1) circle (3pt); 
                        \node[color=black] at (3.3,-0.3) {$x_i$ };
	       	\end{tikzpicture}
 	\caption{ Subdivision of the neighbors of a given point into  subsets, used for the computation of the polynomials adopted in the WENO reconstuction in $2D$.}
 \label{2d_grids}
  \end{center}
  \end{figure}
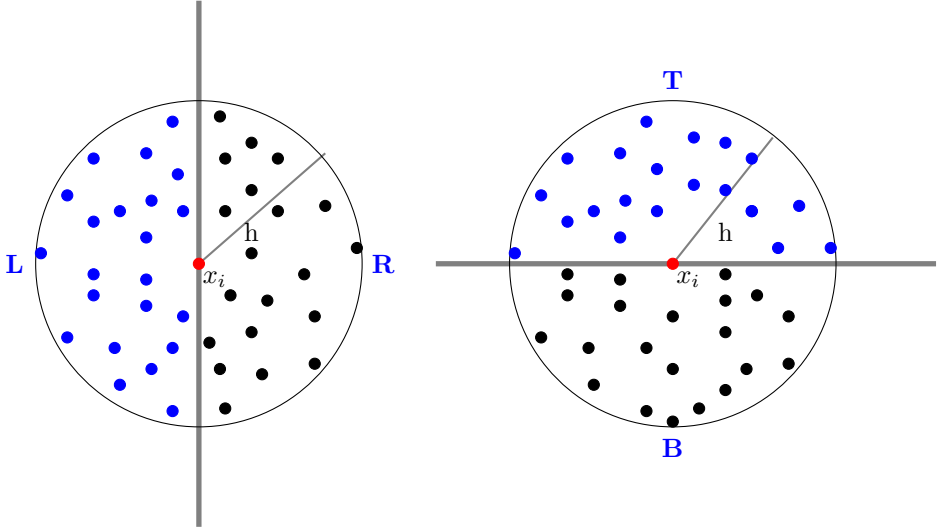

To compute the corresponding weights $w_k,k=L,C,R $ and $w_k, k = B,C,T$ respectively, we use 
the coefficients
\begin{equation}
\beta_k =
\frac{C_k}{({f_x}_k^2 \Delta x^2 + {f_y}_k^2 \Delta x^2 + 
              {f_{xx}}_k^2 \Delta x^4 + {f_{xy}}_k^2 \Delta x^4 + 
             {f_{yy}}_k^2 \Delta x^4 + \epsilon)^2}.
\end{equation}	
%\giovanni{
%In 2D  we distinguish  between 4 subsets of the central neighbour points see Figure \ref{fig:quadrants} and proceed in a similar way.
%Here we have to be more specific. How do we divide the plane? Is it something similar to the picture below, or are the four quadrants the usual one (NE,NW,SW,SE)? We should also explain how we actually compute the reconstruction in 2D
%I think few pictures accompanying the description of the method would make the paper more readable
%}

%\begin{figure}\begin{center}
%\includegraphics[width=0.4\textwidth]{figures/quadrants.pdf}
%\end{center}
%\caption{Subdivision of the neighbors of a given point into four subsets, used for the computation of the polynomials adopted in the WENO reconstuction}
%\label{fig:quadrants}
%\end{figure}
%

 \subsection{Management of grid points}
 
A very important aspect of the proposed ALE meshfree method is the grid management. It consists of three parts, which are presented in the following subsection, see  
 \cite{Kuhnert, DTKB08} for more details.
%In this paper, we consider the one and two dimensional physical spaces. For three dimensional case and other details we refer \cite{DTKB}. 
 
\subsubsection{Initialization of grid points}

The main parameter is the  
average distance between the particles $\Delta x$ which is approximately $\beta h$, where $\beta < 1$. 
First of all we initialize the boundary points %In $1D$ these are just points. In $2D$ boundaries are line segemnts or arcs.
 by establishing grid points on the boundaries at a  distance $\Delta x$. 
To initialize the interior grid points the algorithm starts with  the boundary particles. 
Then, a  first  layer inside 
the domain is constructed. Starting from this layer one proceeds as before  until
the domain is filled with   points having a minimal distance $\beta h$ and a maximal distance $h$. 
%After completion of initial filling, we apply once more the algirithms for hole filling and removing closed grid points, which are described in the following subsections.  
The  initial grid points are not  distributed  
on a regular lattice. Moreover, since the grid points  move, %, the regularity will be destroyed within few time iterations. Then
 they 
 may cluster or scatter in time. %Therefore, we need to monitor in every 
%or every two time steps to maintain the good 
In these cases, a proper quality of the  distribution of the grid points has to be guaranteed with the help of mechanisms to
add and remove points, see below.

\subsubsection{Neighbor search}
Searching neighboring grid points  at an arbitrary position is the most important and time consuming 
part of the  meshfree method. After the initialization, grid points are 
numbered from $1$ to $N$ with 
positions ${x}_i$. The fundamental operation to be done on the point cloud 
is to find for all points  at $x_i$ the neighbors inside a ball $B({x}_i, h)$ with given radius $h$.
To this purpose a voxel data structure containing the computational domain is constructed.  
The voxels form a regular grid of squares with side length $h$. %For each voxel the indices of points contained in it. 
Three types of lists are established. The first one contains the voxels of  all 
points. This is of complexity $O(N)$. 
The second list is  obtained from the first list by sorting with respect 
to the voxel indices. This is of complexity $O(N\log N)$. 
%Of course, in time 
%dependent problems, one generates the list just at the beginning  from scratch and, in 
%the following, updates them from the previous time step. 
Finally,  for each of the points ${x}_i$, all points inside the ball $B({x}_i, h)$ 
have to be determined. 
%The second
% list gives direct access to the voxel, where 
%this point is contained in. 
This is done by testing  all points in the voxel and its 8  neighboring voxels  for being inside the ball
using lists 1 and 2.
%The indices of the points in these voxels are given by the first list. 
Since each voxel contains $O(1)$ points, this operation is of constant 
effort. Hence the total complexity of a neighborhood search for every 
point is $O(N\log N)$. Finally, the neighborhood information is saved in the 
third list.

\subsubsection{Adding and removing points}
%\giovanni{I think the thresholding process is not explained in detail. 
%One possibility is to fix it now, the other is to leave as it is and to fix 
%later, possibly after the referee requires it.}
Determining whether the point-cloud is sufficiently uniform or not and correcting it is more complicated.
%a more tricky task, 
%since only the information are provided where grid points are and not 
%provider where they are not. 
%We note that the center of a void 
%circle of maximal radius must be a Voronoi point.
To determine whether points have to be added, one considers  the  Voronoi cells
\cite{Voronoi}  
of each point ${x}_i$, i.e.  the set of all points closer to ${x}_i$ 
than to any other point.
% In general
 %$3$ Voronoi cells meet in one point called the Voronoi point.
%Finding all Voronoi points  can be done in $O(N \log N)$ time, see, for example \cite{Fortune}
%and  \cite{She99} for details.
%As for hole filling the 
%Voronoi points are
%only to be detected locally, 
%Subdivision algorithms can reduce the 
%effort to $O(N\log N)$, see \cite{She99} for details.
%Obviously constructing the full Voronoi 
%diagram can be too expensive 
%for a large number of points.  If the Voronoi diagram was easy to compute, 
%so would be the Delaunay triangulation.
We note that the existing voxel (or octree) 
structure can be used to construct local, partially overlapping, 
Voronoi diagrams. If the point cloud is not too deformed, such an  approach 
successfully identifies regions with an insufficient number of grid particles  in $O(N)$ time, since the number of 
points considered locally is of order $O(1)$. Once  these regions  are 
identified, new points are inserted. After the insertion of new grid points we use the 
moving least squares interpolation for the approximation of the 
particle distribution function. 
Particles which are too clustered are removed by merging pairs of close by 
points into a single one, see  \cite{Kuhnert, DTKB08} for more details.
% \giovanni{In any case, here a few words of explanations have to be written, or we have to specify where it is explained how this is done.} \tiwari{It is refered in the beginning of particle management}. 
 By an iterative application, also large clusters 
can be thinned out. The two closest points can be found in $O(N)$ time by 
looping over all points and for each point finding its closest neighbor by 
checking all points in its circular neighborhood. With the same procedure, 
one can find all points closer than
a given distance. 
%Introducing a 
%new point is of complexity $O(1)$, since it is one interpolation problem 
%using all neighboring points.
%One prescribes a minimum distance between two points. 
If  two particles, that are closer than this distance, are detected,
both are removed and replaced by a new particle inserted at the center of mass of the two particles under consideration.
The distribution function is interpolated from the neighboring grid points with the help of the moving least squares method.  %%%%%%%%%%%%%%%%%%%%%%%%%%%%%%%%%%%%%

\section{Numerical results}
\label{sec:numerics}
We consider a variety of numerical test cases ranging from smooth  and non-smooth 1D  and 2D solutions of the BGK equation to 1D and 2D moving boundaries with one-way and two-way coupling of moving objects and   gas flow.

\subsection{Example 1: The 1D-BGK model with smooth solution}
\label{1d_conv_study}

For  the convergence study we consider the BGK model  (\ref{ALE1}-\ref{ALE3})  with 1D space and 3D velocity space 
 for short time, compare \cite{PP}. The computational domain 
is  $\Omega = [-1, 1]$. The initial distribution is given by 
\[
f(0,x,v) = \frac{\rho_0}{(2 \pi R T_0)^{3/2}}  \exp {\left(-\frac{|{ v}- { U_0}|^2}{2RT_0}\right)}
\]
with non-dimensional variables and with $R = 1$. Then we choose  $\rho_0 = 1, T_0 = 1$ and $U_0 = (U_0^{(x)},0,0)$, where 
\[
U_0^{(x)} = \frac{1}{\sigma} \left(\exp(\left(-(\sigma x-1)^2\right ) - 2 \exp\left(-(\sigma x + 3)^2\right) \right ), \quad \sigma = 10.
\]
%We note that the initial distribution is locally perturbed in velocity, but still smooth.  
The convergence study is performed up to 
time $t = 0.04$, where the solution is still smooth. We consider a fixed relaxation time $\tau = 10^{-5}$.

In Table \ref{1d_convergence_table1}  the $L^1$ and $L^2$ errors of the temperature determined from the numerical solutions of the first order scheme  are shown. Table \ref{1d_convergence_table2_small_tau}  shows the convergence rate for the ARS(2,2,2) scheme
from section \ref{ARS}. The ARS(2,2,1) scheme from section \ref{RK2} gives very similar results.
Moreover, for larger relaxation time $\tau = 0.1$ and $1$ the convergence rates for  the ARS(2,2,1)  scheme are shown in Tables \ref{1d_convergence_table2_tau0dot1} and \ref{1d_convergence_table2_tau1}. The scheme still produces second order convergence, 
as well as  the ARS(2,2,2) scheme.
The reference 
solution is the solution obtained from a grid with $N_x = 2/{\Delta x}= 801$, where $799$ points are interior points and $2$ are grid points. 
For the convergence study  we used the grid size  $ \Delta x = 0.35 \cdot h$.   
In order to compute the errors, we have generated a mesh with $100$  points 
 and approximated the fluid quantities on this mesh  with the help of 
MLS interpolation from the surrounding grid points. 
For the velocity discretization we use a  uniform grid  with $N_v = 20$   and  
the finite velocity interval $[-v_{\rm max}, v_{\rm max}]$ with $v_{\rm max} = 10$.  
 The time step is always chosen such that the CFL condition 
 $$\Delta t = C \Delta x / v_{\rm max}$$
  with $C=0.5$ is fulfilled for all grid sizes. Noting that the CFL condition for the ALE scheme with a fixed velocity grid is
  \[
  \max_{x,v} \frac{|v-U|\Delta t}{h}<CFL,
  \]
  the above simplified condition essentially means that the difference 
 $ |v-U|$ does not exceed  $2 v_{max}$, which is fulfilled for all examples.
Note that in principle one could have a much better stability condition and use larger time steps, if, as suggested in subsection
\ref{sec:velocity_discretization},  the velocity grid is centered in $U$.

 We observe that all 
  schemes have the expected order of convergence.
  
 \begin{table}
 	\begin{center}
 		\begin{tabular}{|r|r|r|r|r|r|r|}
 			\hline
 			$\Delta t$ &	$h $ & $N_x$ & $L^1$-error & Order& $L^2$ error & Order  \\ \hline 
 			$ 4\cdot10^{-3}$ & $0.28$ &26 &$ 2.50\cdot10^{-2}$   &    $--$ &           $3.02\cdot10^{-2}$ &$--$\\  
 			$ 2\cdot 10^{-3}$ &	$0.14$ & 51&$1.58\cdot 10^{-2}$  &   $0.66$ &       $1.91\cdot 10^{-2}$ &$0.66$\\ 
 			$ 1\cdot 10^{-3}$ &	$0.07$ &101 &$9.02\cdot 10^{-3}$  &   $0.80$ &       $1.10\cdot 10^{-3}$ &$0.78$\\ 
 			$ 5\cdot 10^{-4}$ &	$0.035$ &201 &$4.39\cdot10^{-3}$ &   $1.04$ &       $5.40\cdot 10^{-3}$ &$1.03$\\  
 			$ 2.5\cdot 10^{-4}$ &	$0.0175$ &401 &$1.56\cdot10^{-3}$ &   $1.49$ &       $1.93\cdot 10^{-3}$ &$1.48$\\  
 			\hline
 		\end{tabular}
 		\caption{Example 1: 1-D smooth solution. Convergence  of temperature for $\tau = 10^{-5}$ at time  $t =0.04$ from the first order scheme.}
 		\label{1d_convergence_table1}
 	\end{center}
 \end{table}

   \begin{table}
   	\begin{center}
   		\begin{tabular}{|r|r|r|r|r|r|r|}
   			\hline
   			$\Delta t$ &	$h $ & $N_x$& $L^1$-error & Order& $L^2$ error & Order  \\ \hline 
   			$ 4\cdot10^{-3}$ & $0.28$ &26 &$ 2.24\cdot10^{-2}$   &    $--$ &           $2.90\cdot10^{-2}$ &$--$\\  
   			$ 2\cdot 10^{-3}$ &	$0.14$ & 51&$1.05\cdot 10^{-2}$  &   $1.09$ &       $1.32\cdot 10^{-2}$ &$1.14$\\ 
   			$ 1\cdot 10^{-3}$ &	$0.07$ &101 &$3.29\cdot 10^{-3}$  &   $1.68$ &       $3.89\cdot 10^{-3}$ &$1.76$\\ 
   			$ 5\cdot 10^{-4}$ &	$0.035$ &201 &$6.38\cdot10^{-4}$ &   $2.37$ &       $7.60\cdot 10^{-4}$ &$2.35$\\  
   			$ 2.5\cdot 10^{-4}$ &	$0.0175$ &401 &$1.57\cdot10^{-4}$ &   $2.03$ &       $1.51\cdot 10^{-4}$ &$2.33$\\  
   			\hline
   		\end{tabular}
   		\caption{Example 1: 1-D smooth solution. Convergence  of temperature for $\tau = 10^{-5}$ at time  $t =0.04$ for the ARS(2,2,2)  scheme.}
   		\label{1d_convergence_table2_small_tau}
   	\end{center}
   \end{table}
   
%%%%%%%%%%%%%%%%%%%%%%%%  
   \begin{table}
   	\begin{center}
   		\begin{tabular}{|r|r|r|r|r|r|r|}
   			\hline
   			$\Delta t$ &	$h $ & $N_x$ & $L^1$-error & Order& $L^2$ error & Order  \\ \hline 
   			$ 4\cdot10^{-3}$ & $0.28$ &26 &$ 1.68\cdot10^{-2}$   &    $--$ &           $2.13\cdot10^{-2}$ &$--$\\  
   			$ 2\cdot 10^{-3}$ &	$0.14$ & 51&$7.36\cdot 10^{-3}$  &   $1.19$ &       $8.96\cdot 10^{-3}$ &$1.25$\\ 
   			$ 1\cdot 10^{-3}$ &	$0.07$ &101 &$2.41\cdot 10^{-3}$  &   $1.61$ &       $2.78\cdot 10^{-3}$ &$1.69$\\ 
   			$ 5\cdot 10^{-4}$ &	$0.035$ &201 &$5.10\cdot10^{-4}$ &   $2.23$ &       $5.82\cdot 10^{-4}$ &$2.26$\\  
   			$ 2.5\cdot 10^{-4}$ &	$0.0175$ &401 &$9.93\cdot10^{-5}$ &   $2.36$ &       $1.21\cdot 10^{-4}$ &$2.27$\\  
   			\hline
   		\end{tabular}
   		\caption{Example 1: 1-D smooth solution. Convergence of temperature for $\tau = 0.1$ at time  $t =0.04$ for the  ARS(2,2,1) scheme.}
   		\label{1d_convergence_table2_tau0dot1}
   	\end{center}
   \end{table}

%%%%%%%%%%%%%%%%%%%%%%%%%%%%  
   \begin{table}
   	\begin{center}
   		\begin{tabular}{|r|r|r|r|r|r|r|}
   			\hline
   			$\Delta t$ &	$h $ &  $N_x$  & $L^1$-error & Order& $L^2$ error & Order  \\ \hline 
   			$ 4\cdot10^{-3}$ & $0.28$ &26 &$ 1.65\cdot10^{-2}$   &    $--$ &           $2.08\cdot10^{-2}$ &$--$\\  
   			$ 2\cdot 10^{-3}$ &	$0.14$ & 51&$7.24\cdot 10^{-3}$  &   $1.19$ &       $8.79\cdot 10^{-3}$ &$1.24$\\ 
   			$ 1\cdot 10^{-3}$ &	$0.07$ &101 &$2.41\cdot 10^{-3}$  &   $1.59$ &       $2.78\cdot 10^{-3}$ &$1.66$\\ 
   			$ 5\cdot 10^{-4}$ &	$0.035$ &201 &$5.27\cdot10^{-4}$ &   $2.20$ &       $5.98\cdot 10^{-4}$ &$2.22$\\  
   			$ 2.5\cdot 10^{-4}$ &	$0.0175$ &401 &$1.06\cdot10^{-4}$ &   $2.32$ &       $1.25\cdot 10^{-4}$ &$2.26$\\  
   			\hline
   		\end{tabular}
   		\caption{Example 1: 1-D smooth solution. Convergence of temperature for $\tau = 1$ at time  $t =0.04$ for the  ARS(2,2,1) scheme.}
   		\label{1d_convergence_table2_tau1}
   	\end{center}
   \end{table}

%%%%%%%%%%%%%%%%%%%%%%%%%%%%  
In Figure \ref{bgk_sol} we have plotted density, mean velocity and  temperature for   $N_x = 160$ grid points obtained from both schemes at time $t=0.04$ together with the reference solution.
In all figures the improved approximation quality of the  ARS  schemes can be clearly observed.

\begin{figure}[!t]
 \centering
 \includegraphics[keepaspectratio=true, width=.329\textwidth]{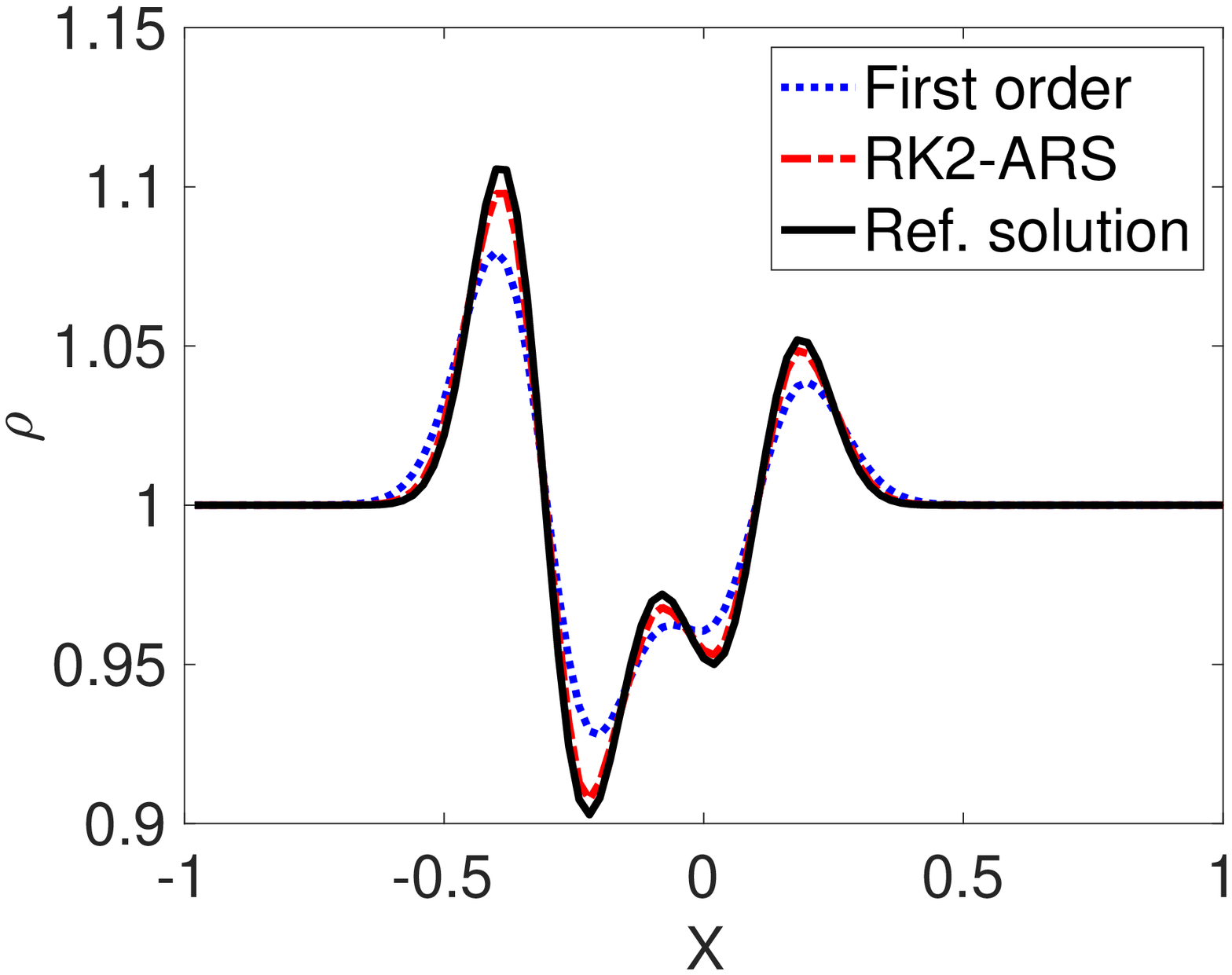}
 \includegraphics[keepaspectratio=true, width=.329\textwidth]{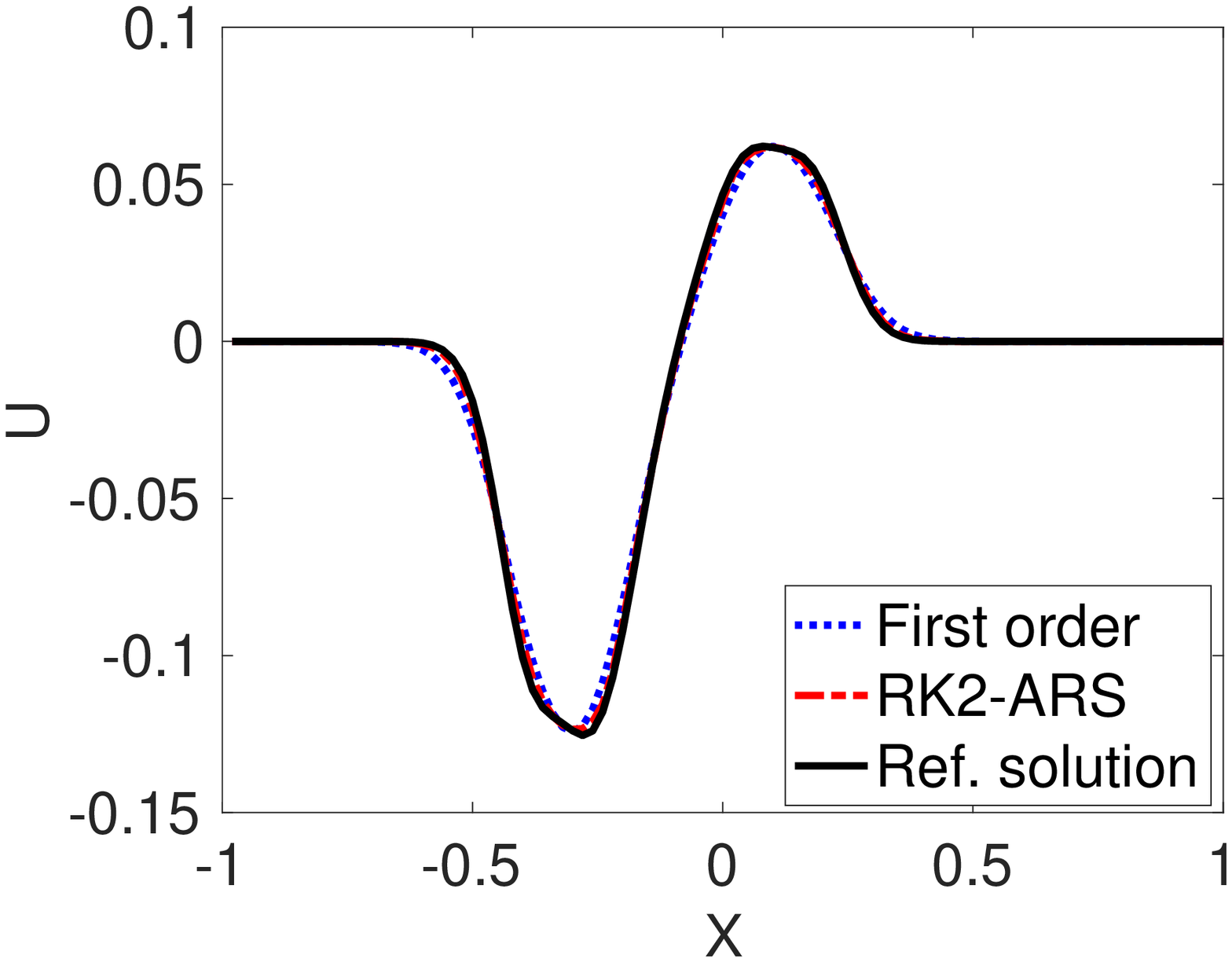}
 \includegraphics[keepaspectratio=true, width=.329\textwidth]{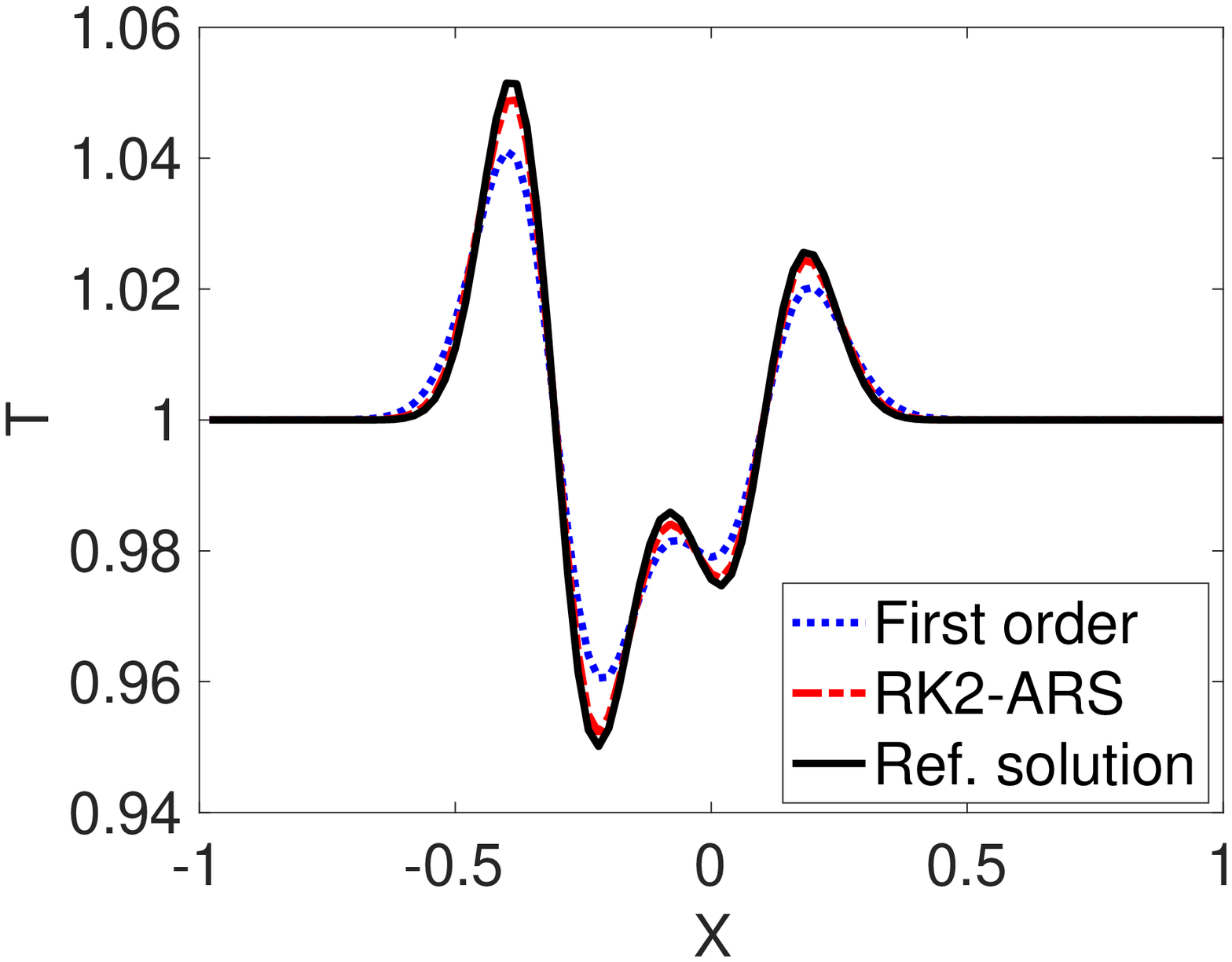}
%\includegraphics[width=5.2cm,angle=0]{figures/bgk_rho_t_0dot04} 
%\includegraphics[width=5.2cm,angle=0]{figures/bgk_U_t_0dot04} 
%\includegraphics[width=5.2cm,angle=0]{figures/bgk_T_t_0dot04}
% \vspace{-0.5cm}
  \caption{Example 1: 1-D smooth solution. Comparison of density, mean velocity and temperature computed from the reference solution and from the solutions obtained from the first order, ARS  schemes for  $N_x = 100$ initial grid points and $N_v = 30$ at time $T_{final} = 0.04$.}
  \label{bgk_sol}
\end{figure}
%%%%%%%%%%%%%%%%%

 \subsection{Example 2: The 1D-BGK model for a Riemann  problem}
 \label{sod}

We consider  a Riemann problem similar to Sod's shock tube problem \cite{SOD}  to validate the  numerical schemes for discontinuous solutions. On the one hand, we compare the first and second order numerical solutions of equations (\ref{ALE1}-\ref{ALE3}) with
a very small value of the relaxation time  $\tau$ to 
 the hydrodynamic limit solution, i.e. the solution of the Euler equations. 
 On the other hand, the numerical solutions of the BGK equation for larger values of $\tau$ are considered and compared to 
 other numerical results and to DSMC solutions.
 
We consider  the computational domain $ [0, 1]$.  
The initial condition is  a Maxwellian distribution with the initial parameters 
\[
\rho_l = 10^{-3}, \;\;U_l^{(x) }= 0, \;\;T_l = 273  \; \; 
\mbox{for} \;\; 0\le x < 0.5
\]
\[
\rho_r = 0.125\times10^{-3},\;\; U_r^{(x)}= 0,\;\; T_r = 273  \;\; 
\mbox{for}\;\; 0.5\le x \le 1.
\]
Diffuse reflection  boundary conditions are   applied and SI units with  the gas constant $R=208$ are  chosen.
The initial values of $\lambda$ and $\tau$ on the left half of the domain are computed according to equations (\ref{lambda}) and (\ref{tau}). We obtain 
$\lambda= 1.110\times 10^{-4}$ and $\tau = 3.69\times 10^{-7}$, respectively. The values on the right half of the domain are $8$ times larger.  
During the time evolution we consider variable relaxation times given by equations (\ref{tau}).

We use a uniform velocity grid with $N_v = 30$. 
We have chosen  the time step $\Delta t = 5 \cdot 10^{-4}$ which leads again to a CFL condition with constant $0.5$.
%, therefore, the 
%time step is much larger than the other schemes, like IMEX \cite{PP} for the BGK model. 
 The computation is performed up to $t = 0.0008$. Initially $N_x = 400$ grid points are 
generated uniformly with spacing $\Delta x = 1/N_x$. The  radius $h$ fulfills  again  $ \Delta x = 0.3 \cdot h$.   
In Figure \ref{fig:sod} we have plotted the numerical solutions obtained by first and second order  schemes together with the analytical solutions of 
the compressible Euler equations. The improved accuracy of the  second order scheme is clearly observed.

 \begin{figure}[!t]
 \centering
 \includegraphics[keepaspectratio=true, width=.329\textwidth]{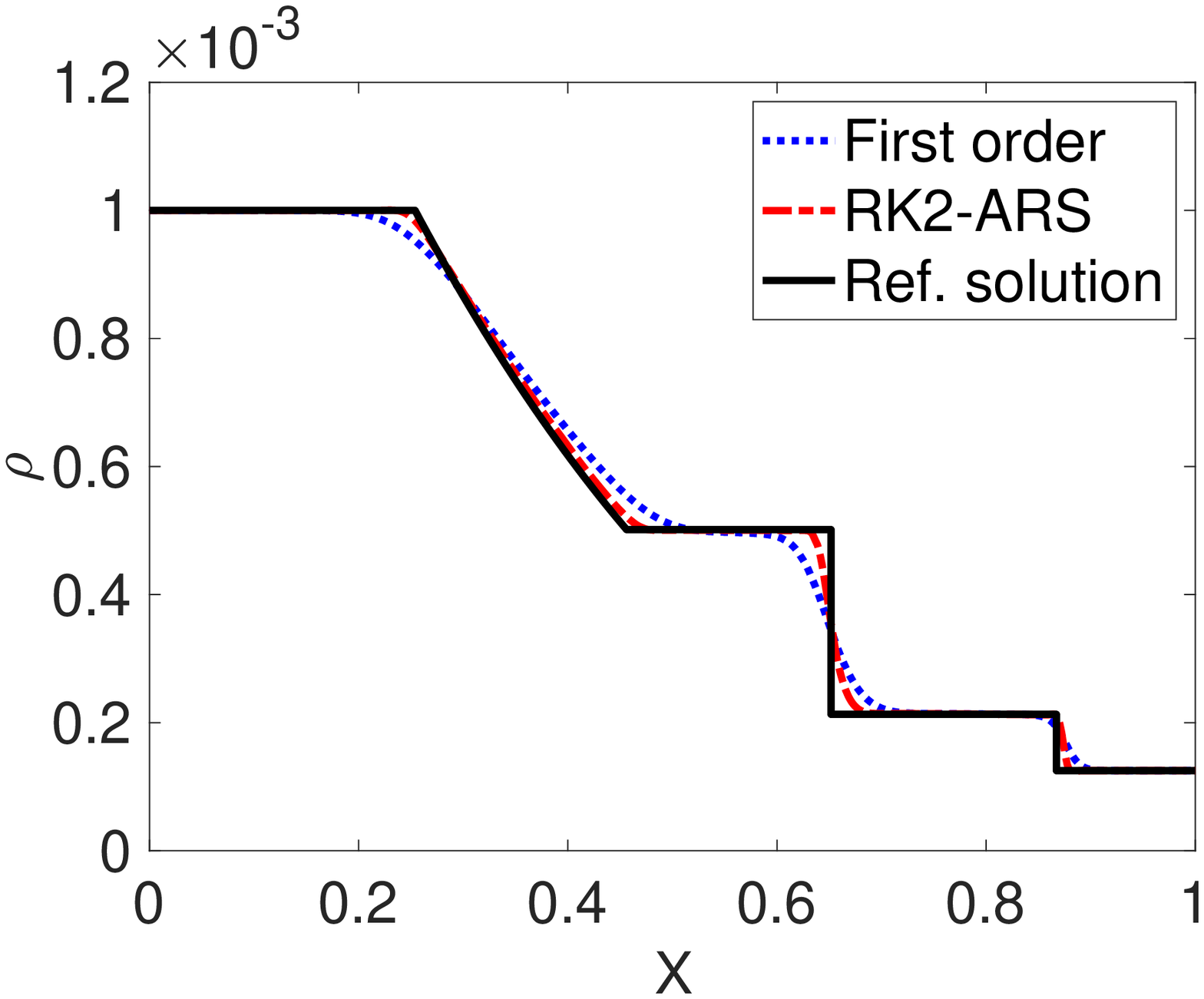}
 \includegraphics[keepaspectratio=true, width=.329\textwidth]{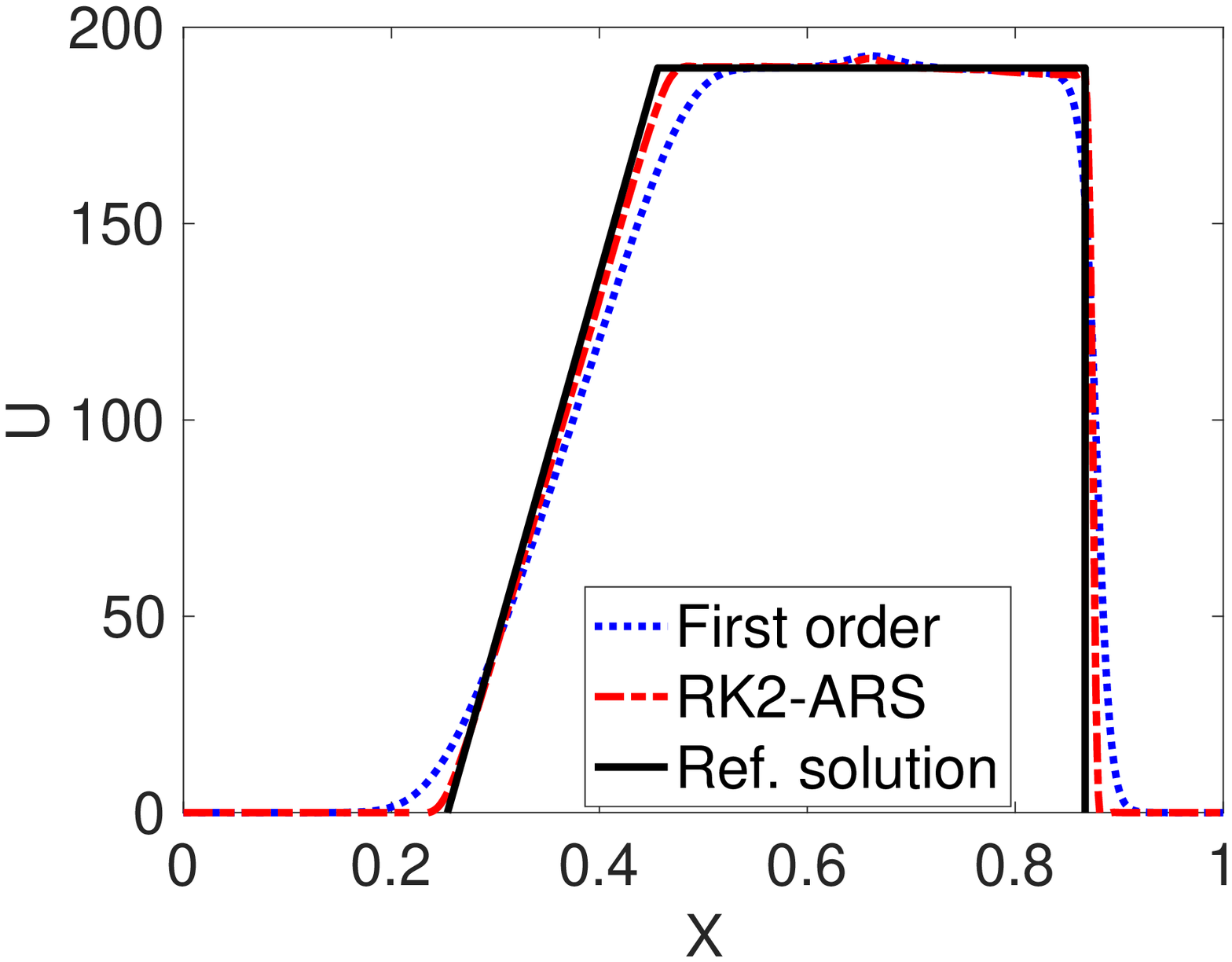}
 \includegraphics[keepaspectratio=true, width=.329\textwidth]{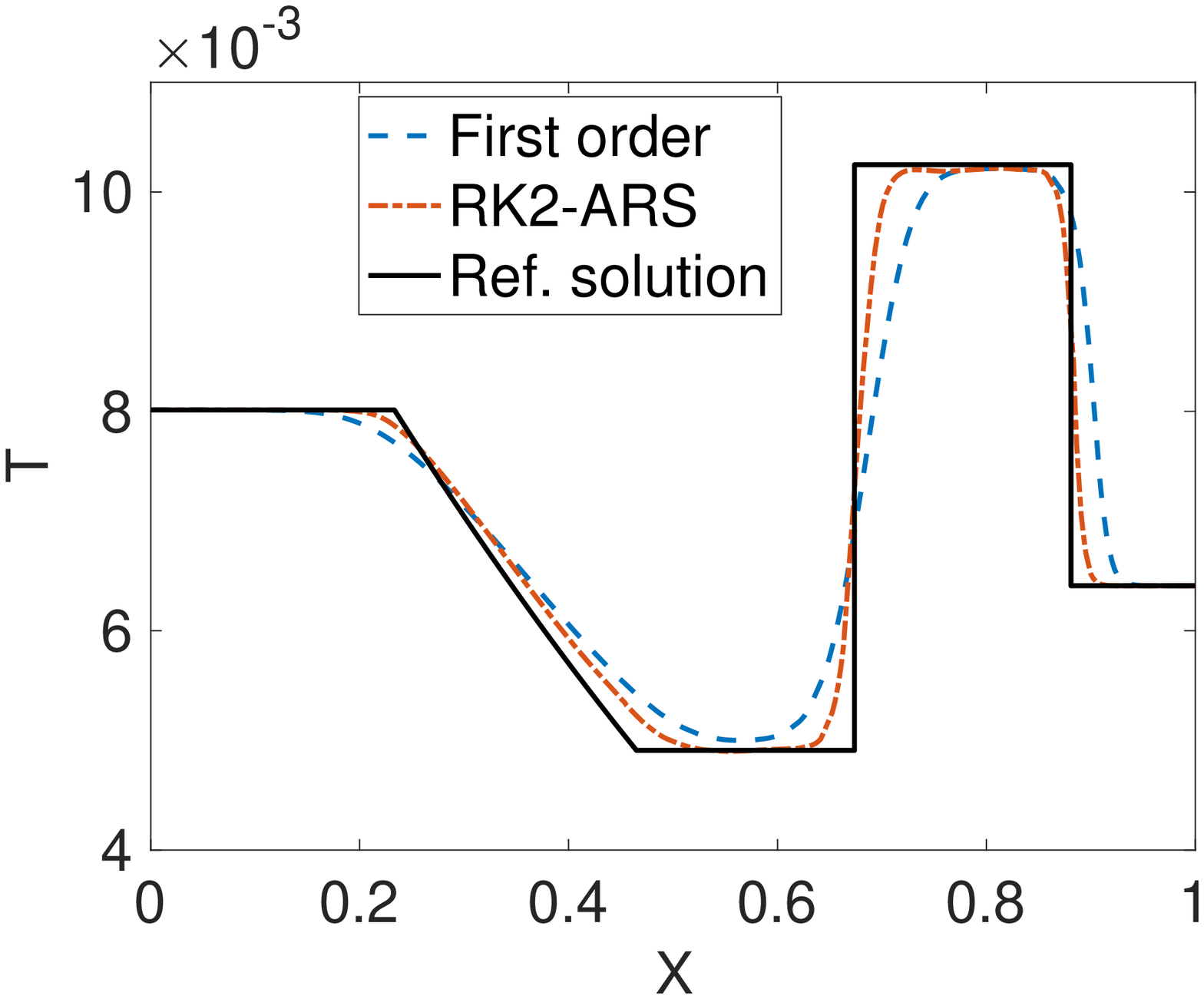}
%\includegraphics[width=5.2cm,angle=0]{figures/bgk_shock_rho} %\hspace{0.5cm}
%\includegraphics[width=5.2cm, angle = 0]{figures/bgk_shock_U} %\hspace{0.5cm}
%\includegraphics[width=5.2cm,angle=0]{figures/bgk_shock_P} %\hspace{0.5cm}
%\vspace{-1cm}
  \caption{Example 2: 1D shock tube. Comparison of the exact solutions of  the Euler equations and the numerical solutions 
  	of the BGK problem with $\tau = 3.69\times 10^{-7}$ initially for a shock tube problem with initial $N_x= 400$ grid points and $N_v = 30$.}
  \label{fig:sod}
\end{figure}
  
%%%%%%%%%%%%%%%%%%%%%%%%%  

As already stated, we use this example also  to consider the solutions of the BGK model for larger values of  $ \tau$ and compare them  with those of the full Boltzmann  equation. 
As before we use   relaxation times  $\tau$ according to equation
(\ref{tau}).
 The density ratio between  left and right part of the domain is again ${\rho_l}/{\rho_r} = 8$,
 and we consider two more rarefied cases with $\rho_l = 10^{-6},  \rho_l = 10^{-4}$, respectively, with corresponding 
 values of the initial  relaxation times  $\tau_l = 3.69 \times 10^{-4}, 3.69 \times 10^{-6}$
 determined from (\ref{tau}). 
In the following figures \ref{fig:sod_case1} to \ref{fig:sod_case2} we have plotted the density, velocity and pressure obtained from the Boltzmann equation and the BGK model at the final time $0.0008$ .  For the Boltzmann equation we consider a hard sphere monatomic gas. The solutions of the Boltzmann equation are obtained from a DSMC simulation averaging  20 independent runs. 
One observes in Fig.~\ref{fig:sod_case1}  and Fig.~\ref{fig:sod_case2} that the solutions of the BGK model coincide  with those of the Boltzmann equation for both  values of the relaxation time $\tau$. Note that for  the larger value of $\tau$, see Fig. \ref{fig:sod_case1}, we have used a number of velocity grid points equal to $N_v = 200$ to avoid  oscillating solutions of the BGK model.\footnote{This behaviour is typical for problems with large Knudsen number: the interaction among among gas particles is weaker and a greater resolution in velocity is needed to resolve the distribution in phase space and avoid spurious oscillations.} 

\begin{figure}[!t]
 \centering
 \includegraphics[keepaspectratio=true, width=.329\textwidth]{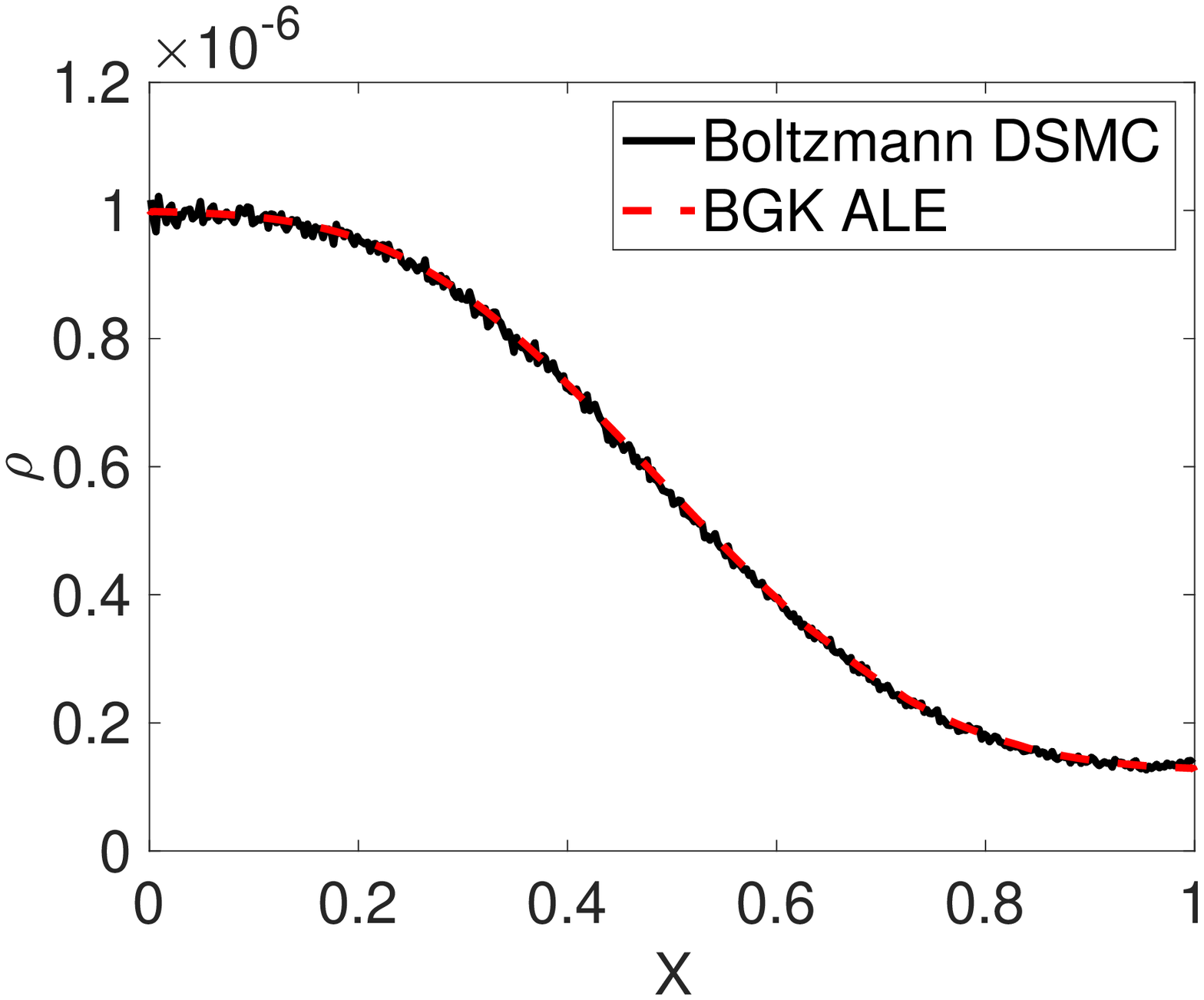}
 \includegraphics[keepaspectratio=true, width=.329\textwidth]{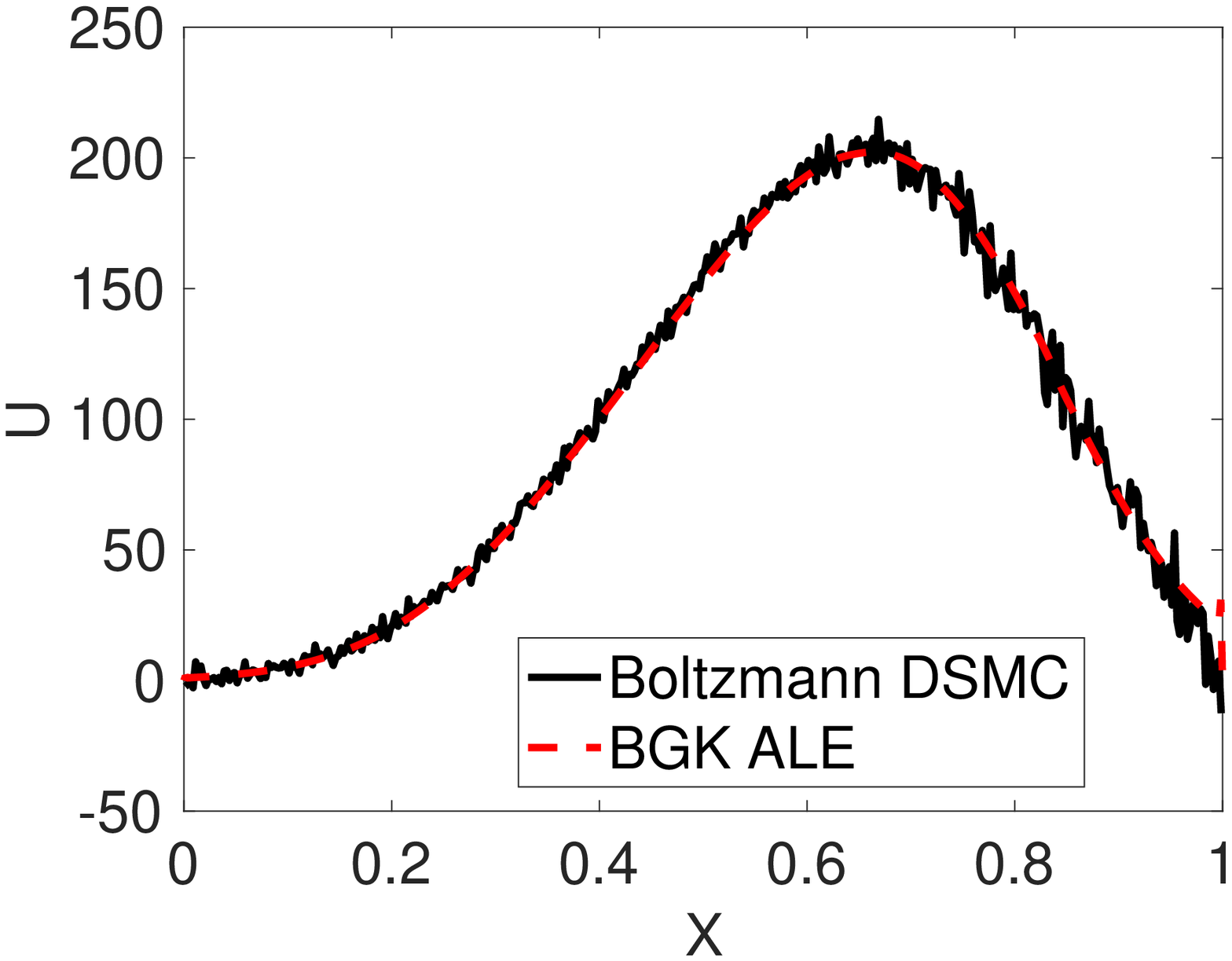}
 \includegraphics[keepaspectratio=true, width=.329\textwidth]{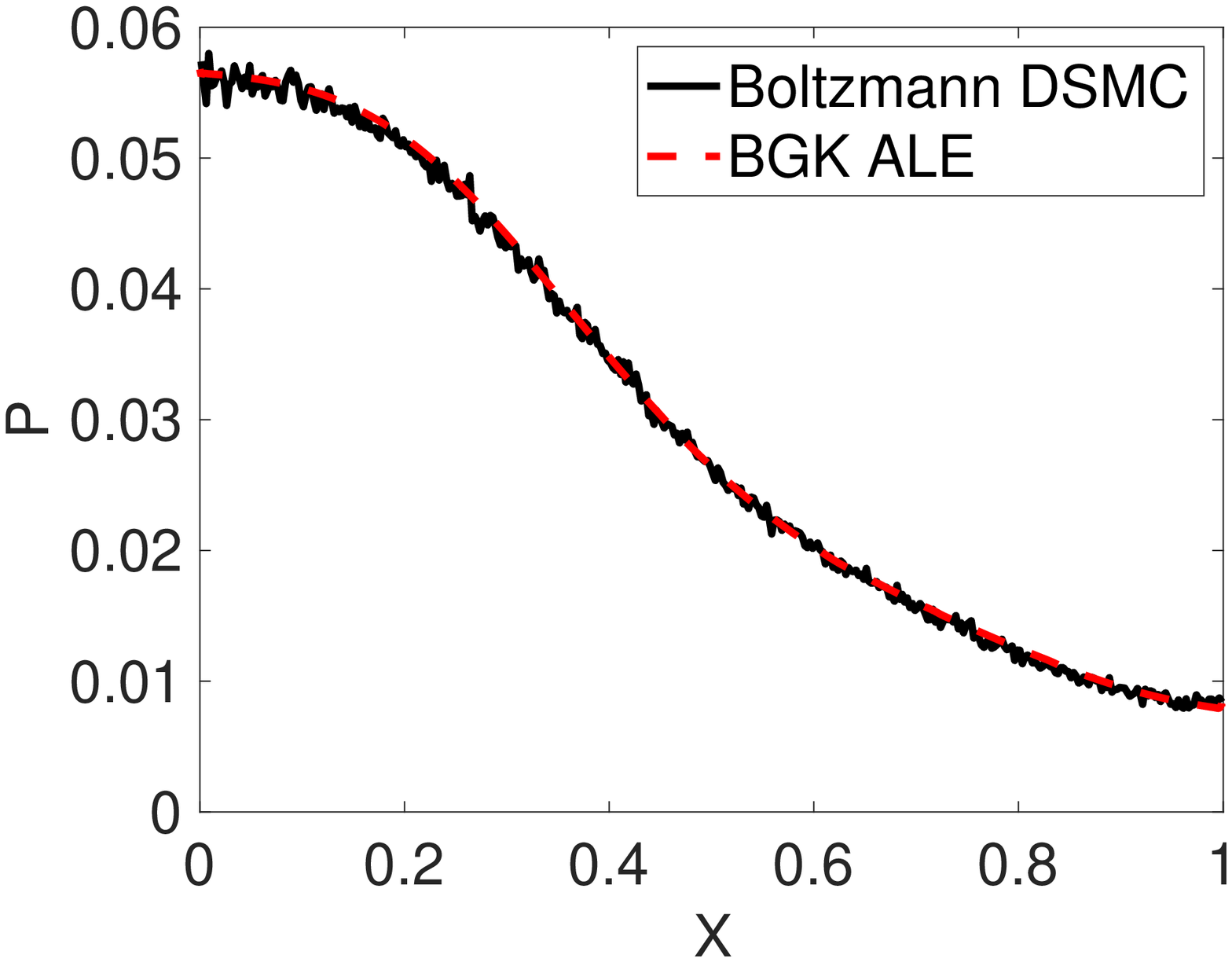}
 \caption{Example 2: 1D shock tube. Comparison of the solutions obtained from the Boltzmann equation with DSMC and the BGK model with $N_v = 200$ for 
 $ \rho_l = 10^{-6}$ and the corresponding 
  initial  relaxation times  $\tau_l = 3.69 \times 10^{-4}$.}
  \label{fig:sod_case1}
\end{figure}
  
\begin{figure}[!t]
 \centering
\includegraphics[keepaspectratio=true, width=.329\textwidth]{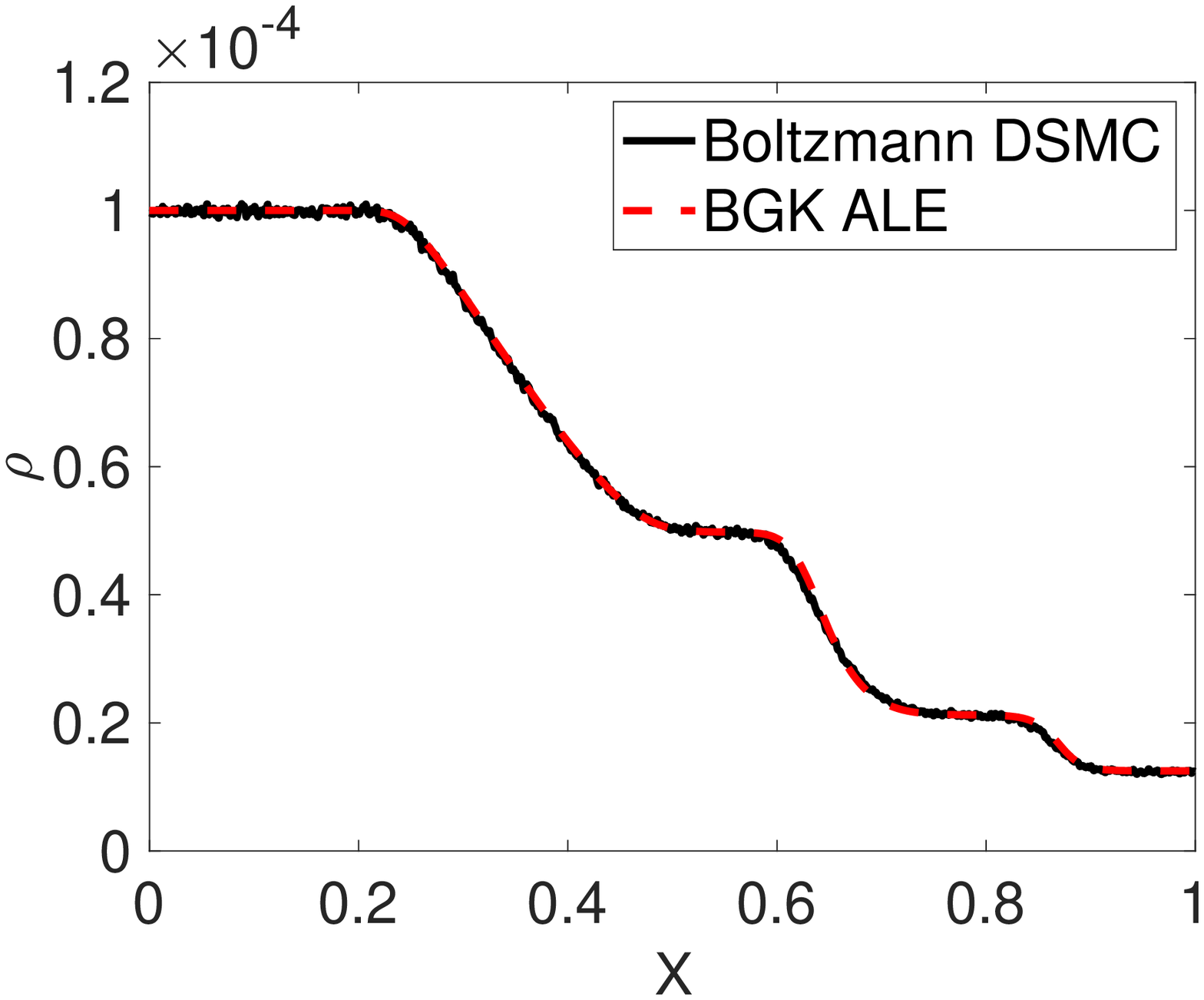} 
\includegraphics[keepaspectratio=true, width=.329\textwidth]{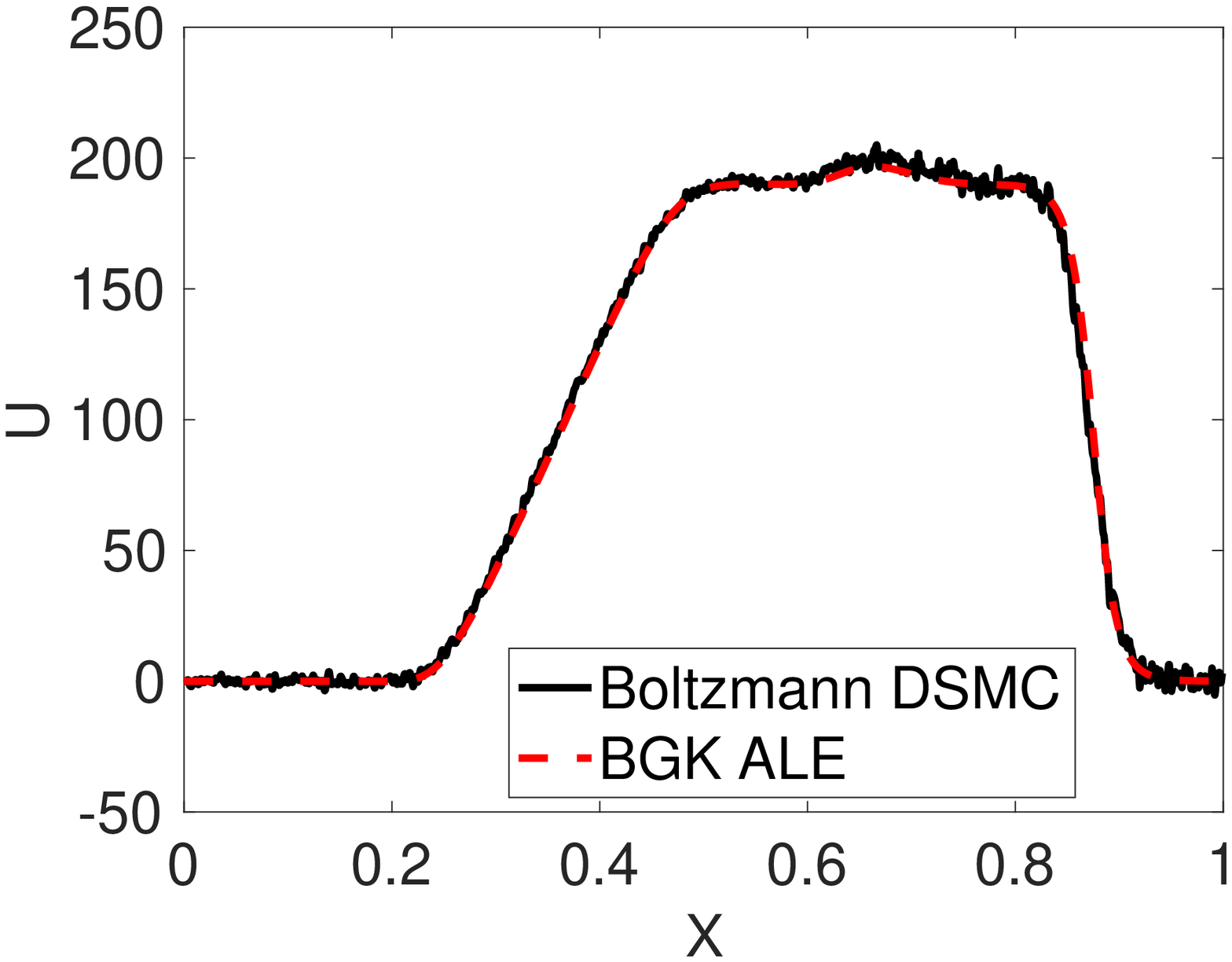} 
\includegraphics[keepaspectratio=true, width=.329\textwidth]{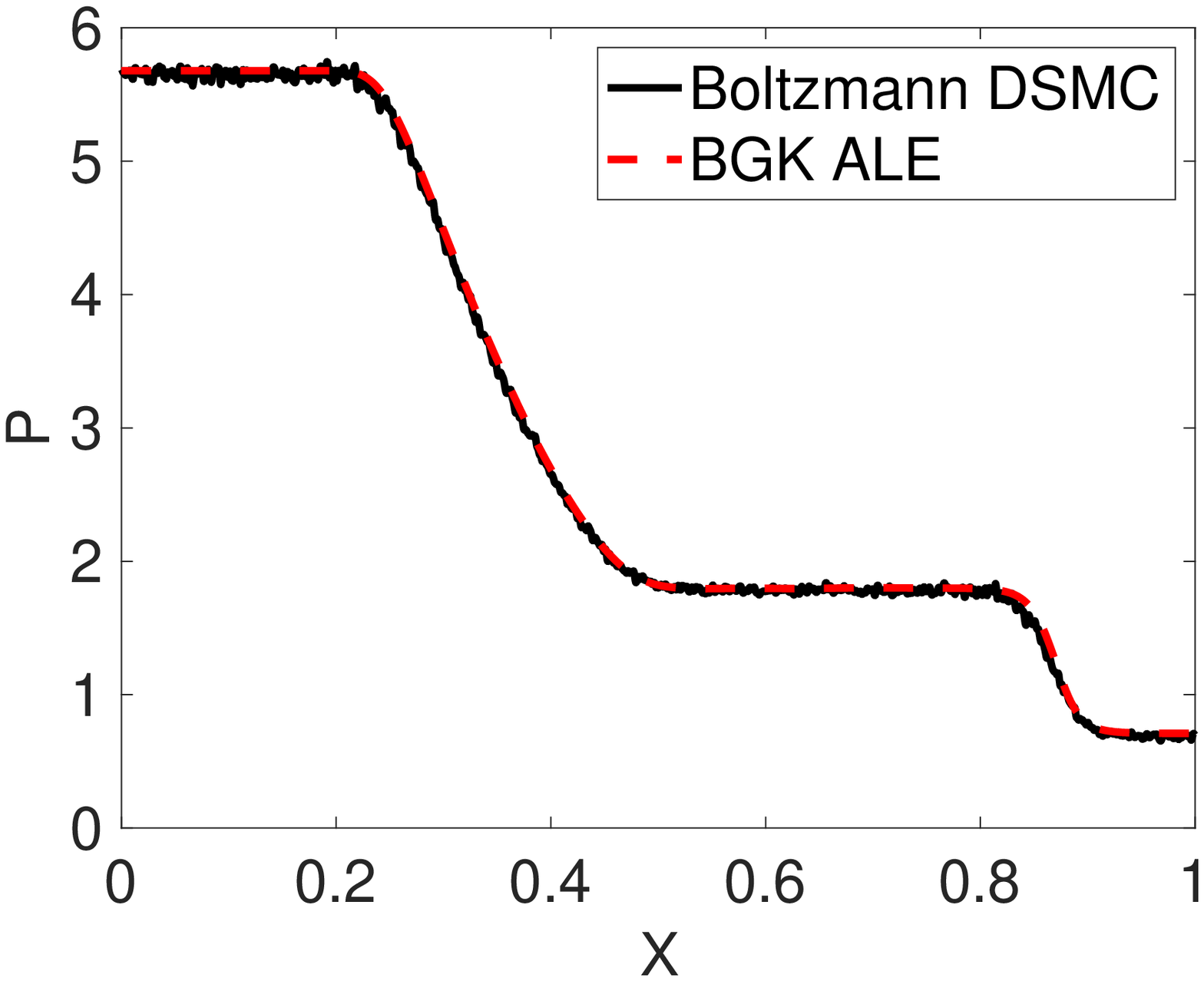} 
 \caption{Example 2: 1D shock tube. Comparison of the solutions obtained from the Boltzmann equation with DSMC and the BGK model with $N_v = 30$ for $\rho_l = 10^{-4}$ and the corresponding 
 	initial  relaxation times  $ 3.69 \times 10^{-6}$.}
  \label{fig:sod_case2}
\end{figure}  

%%%%%%%%%%%%%%%%%%%%%%%

Furthermore, we compare the solutions of the BGK model obtained from the ALE method presented here with a higher order semi-Lagrangian (SL) scheme, see \cite{CBRY1,CBRY2}.  
 We consider  the initial densities $\rho_l = 10^{-4}$.  In  Fig. \ref{fig:sod_SL_ALE2a} we have compared the densities obtained from ALE and SL scheme  for different spatial resolutions. The solutions match perfectly  well for  a larger number  of spatial grid points
like  $N_x = 400$, see   Fig. \ref{fig:sod_SL_ALE2a} on the right. We use this  solution as the reference solution and compare it to the 
ALE and SL solutions for coarser grids.
One observes that for $N_x = 50$ and $N_x=100$ the solutions obtained from the ALE method deviates slightly from the reference solution, whereas the higher order SL solutions are still very near to the reference solution.

 \begin{figure}[!t]
 \centering
\includegraphics[keepaspectratio=true, width=.329\textwidth]{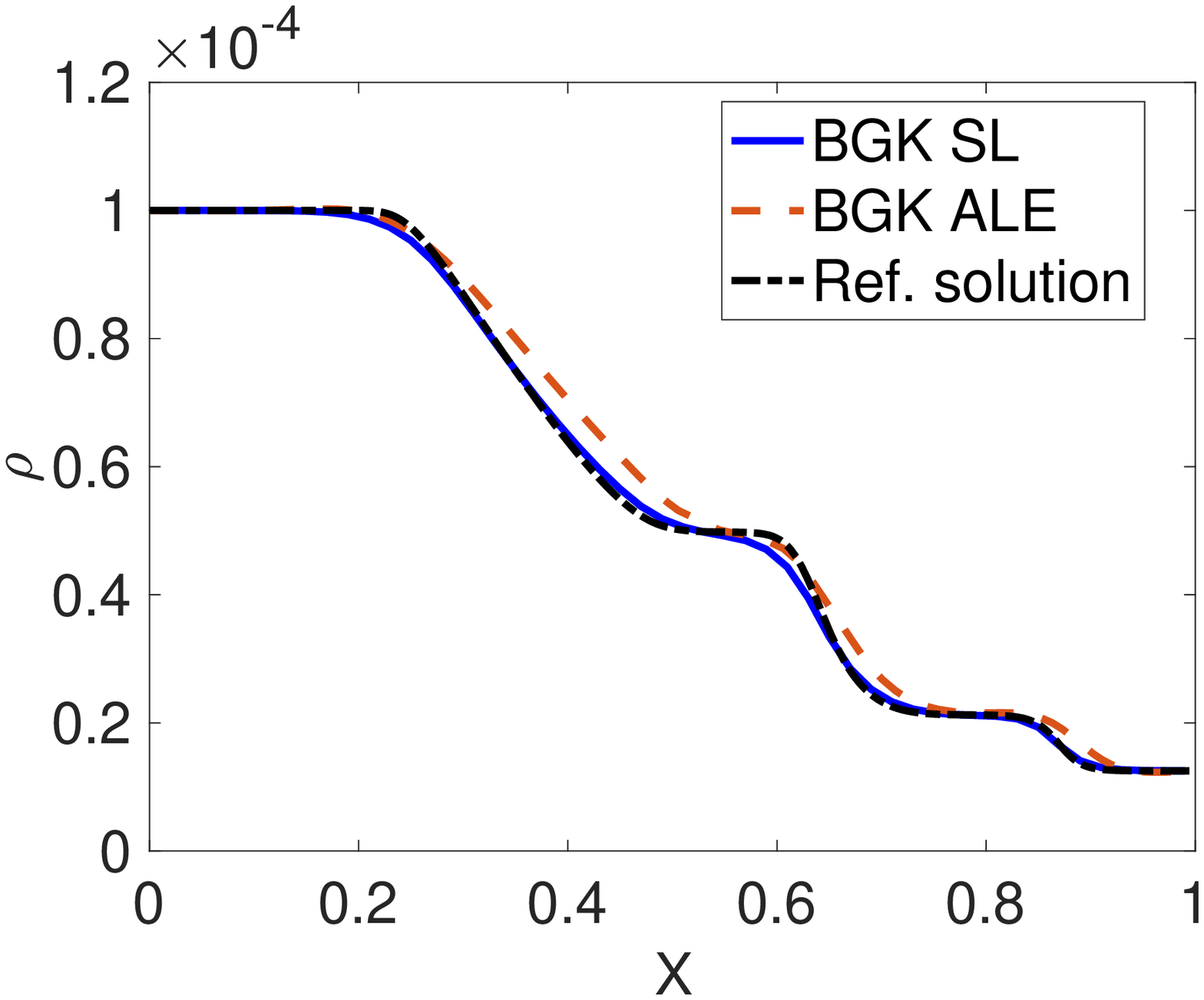} 
\includegraphics[keepaspectratio=true, width=.329\textwidth]{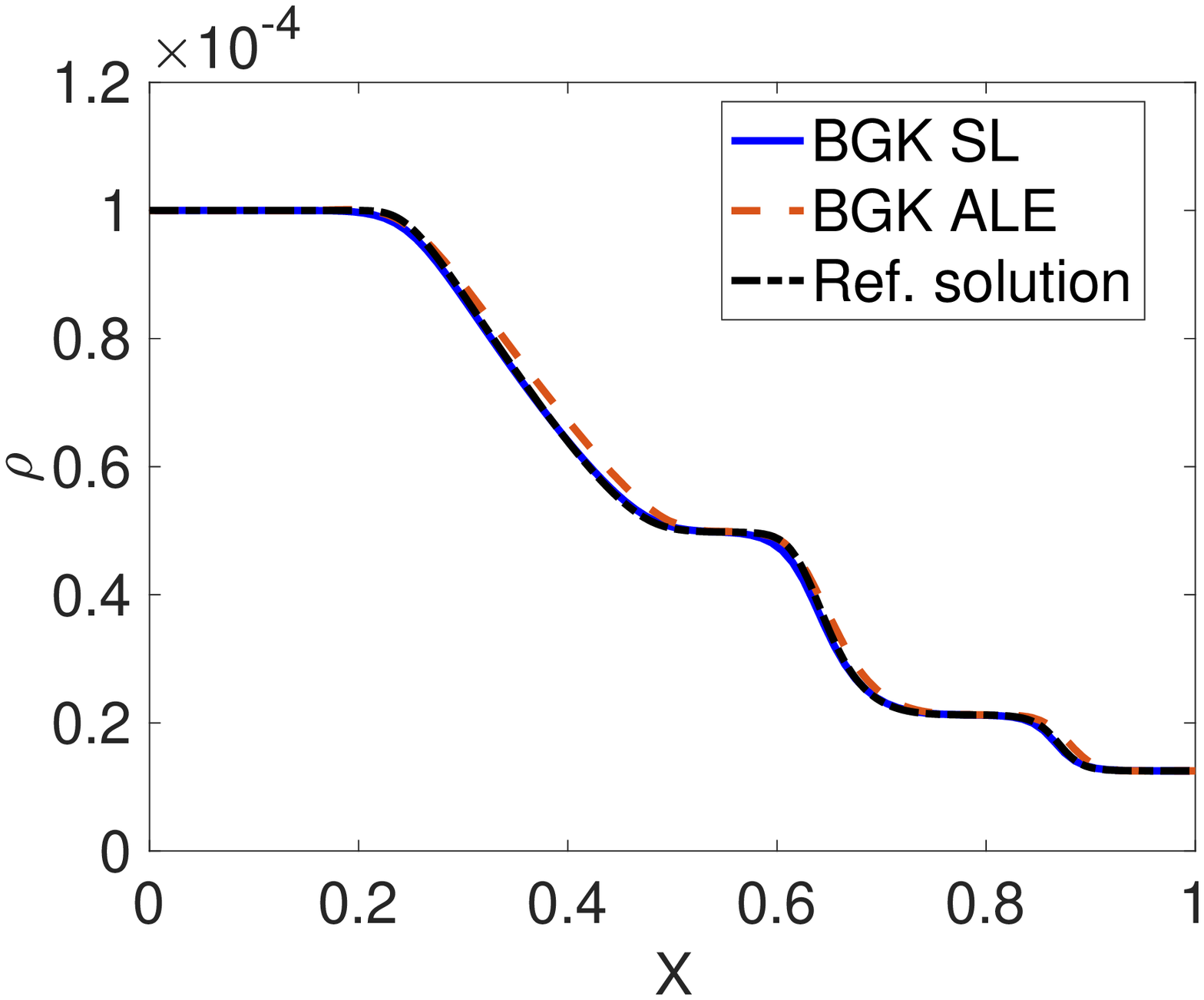} 
\includegraphics[keepaspectratio=true, width=.329\textwidth]{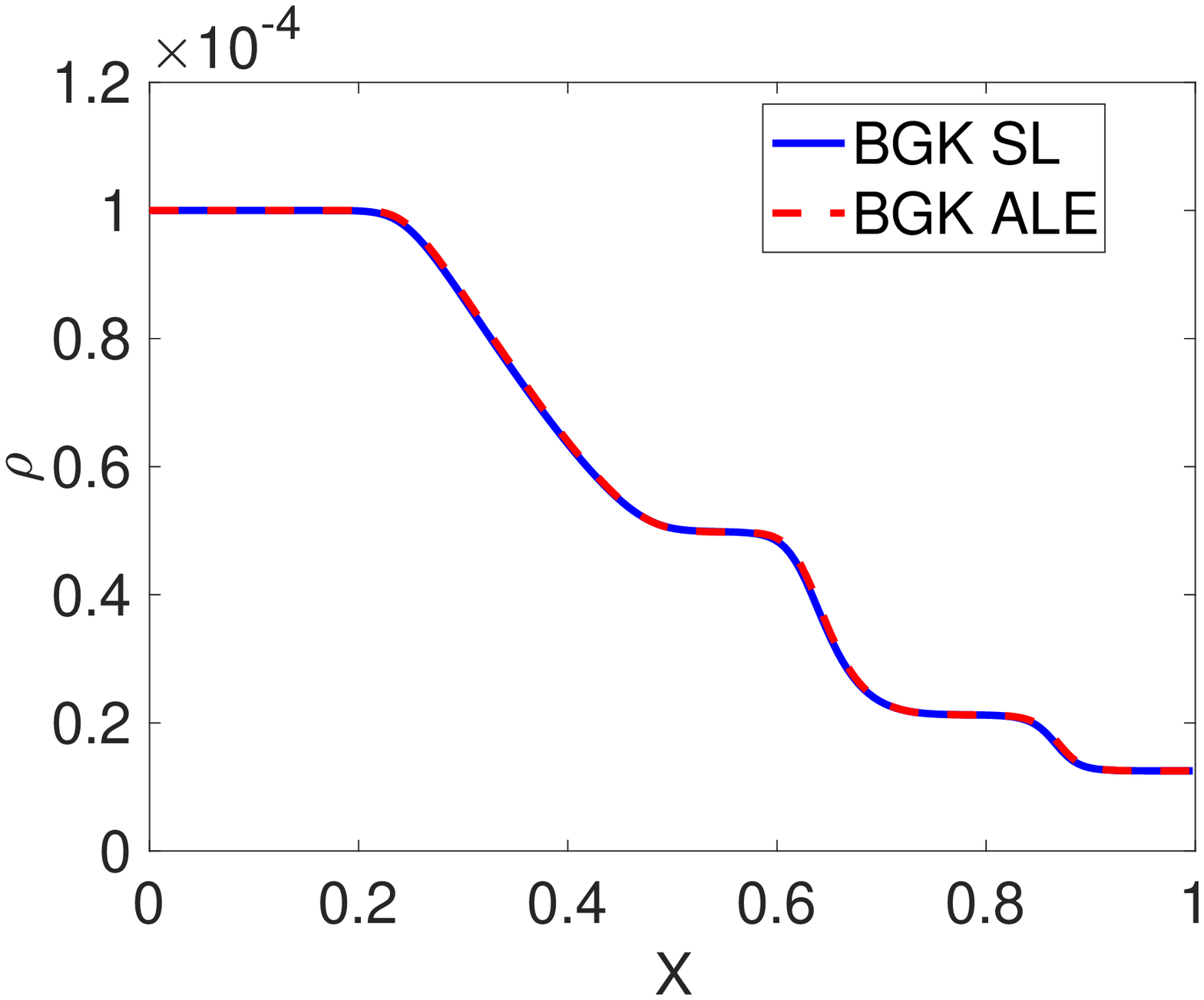} 
 \caption{Example 2: 1D shock tube. Comparison of the solutions obtained from the semi-Lagrange scheme and the ALE method with $N_x = 50$ (left), $N_x = 100$ (middle) and $N_x = 400$ (right) for the BGK model with $N_v = 30$ for 
 $ \rho_l = 10^{-4}$ and the corresponding  initial  relaxation times  $\tau_l = 3.69 \times 10^{-6}$.}
  \label{fig:sod_SL_ALE2a}
\end{figure}

%%%%%%%%%%%%%%%%%%%%%%%%%%%
\subsection{Example 3: Moving piston with prescribed velocity}
 
This problem has been considered in \cite{DM, RF} in a larger domain. We consider the one-dimensional domain $\Omega=[0, 20]$.
  Initially the piston is positioned at $x = 2$. We consider a 
total number $N_x=300$  grid points in physical space and $N_v = 20$ grid points in velocity space. 
The left boundary moves with velocity 
\[
u_p = 0.25*\sin(t).
\]
Again we use non-dimensional variables with $R=1$. The initial velocity is $U_0 = 0$, the density $\rho_0 = 0.001$ and the temperature $T_0 = 1$. The minimum and maximum  of the velocity are 
$v_{\rm min} = -10 $ and $v_{\rm max}=10$. The initial distribution is 
the Maxwellian with the above  initial macroscopic quantities. 
%The position of the piston is the boundary on the left side. 
Initially particles are generated  in the interval $[2, 20]$. 
We have considered a fixed value of $\tau=    1.83 \cdot 10^{-2}  $, 
%  an  initial mean free path is $2.29\cdot10^{-2}$
 a final time $t_{final} = 4 $ and a  time step  $\Delta t = 0.001$. 
%, see Figure \ref{piston_geom} for the physical setup of the problem. 
As in the previous section we compare the solutions obtained by the numerical method for the BGK equations 
to 
the solution obtained from a DSMC simulations of the full Boltzmann equation with a moving geometry, see  \cite{STKH}.
For the DSMC method we use  $\Delta x=20/900 = 2.22\cdot 10^{-2}$. In order to obtain a smooth 
solution for the DSMC simulations we  have performed 50 independent runs. 

Figure \ref{moving_piston1} to Figure \ref{moving_piston4_zoom} show the results for different times.
When the piston starts to move in time,  two situations occur:  when the velocity is positive, the grid points are approaching  each other. In this case 
one has to  remove the grids points which are too close. We replace two grid points by a new one and locate it in the center between the two. When the velocity is negative new grid points have to be added. In both cases the  distribution functions have to be updated in the  additional grid  points. This is done with the help of a least squares interpolation.

%To validate the numerical results of the time splitting ALE scheme for the BGK model, we compare it with the results of a DSMC simulations of the BGK model \cite{Bird, Babovsky, NS, DP}.  
%  
 
 \begin{figure}[!t]
\includegraphics[keepaspectratio=true, width=.329\textwidth]{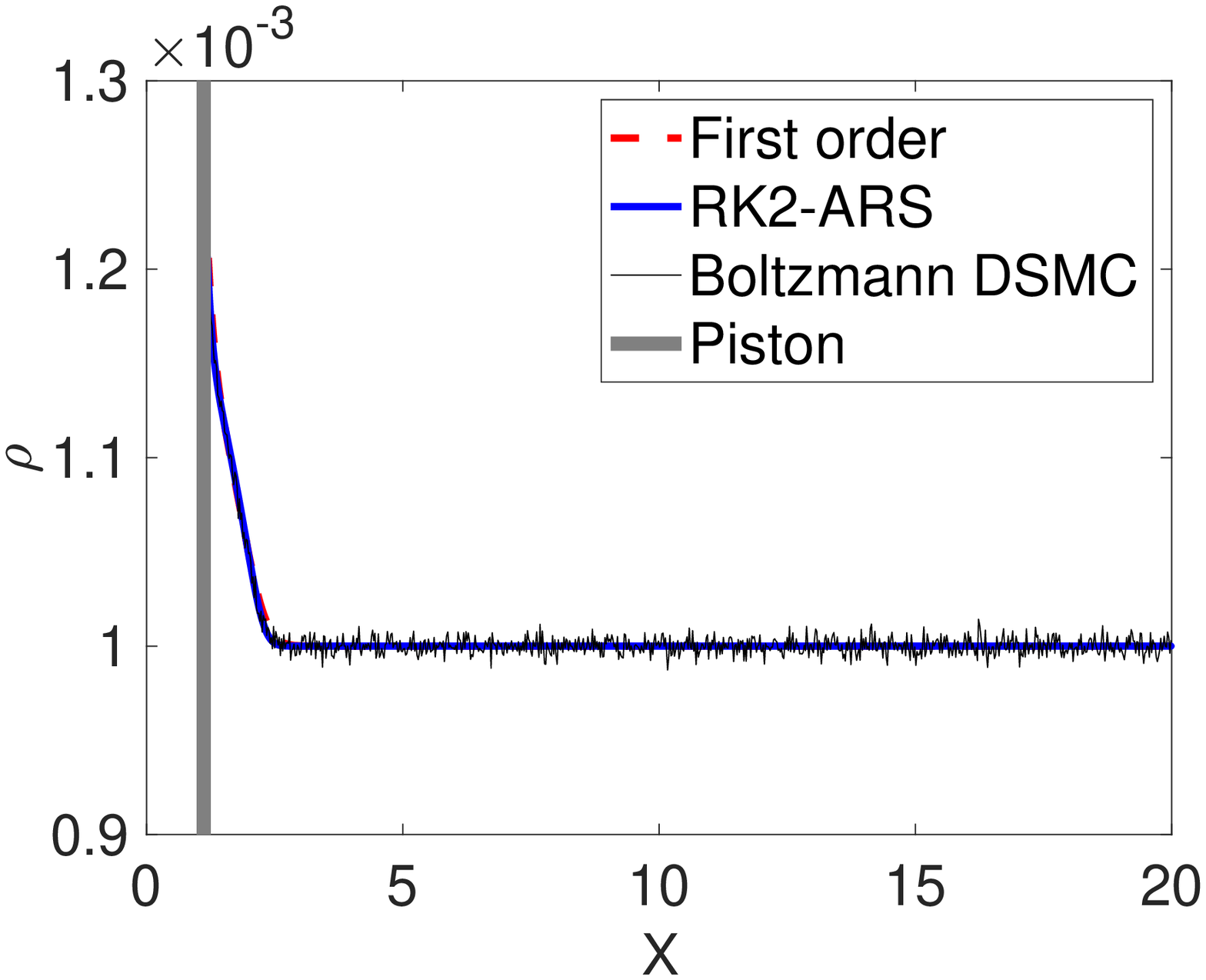}
\includegraphics[keepaspectratio=true, width=.329\textwidth]{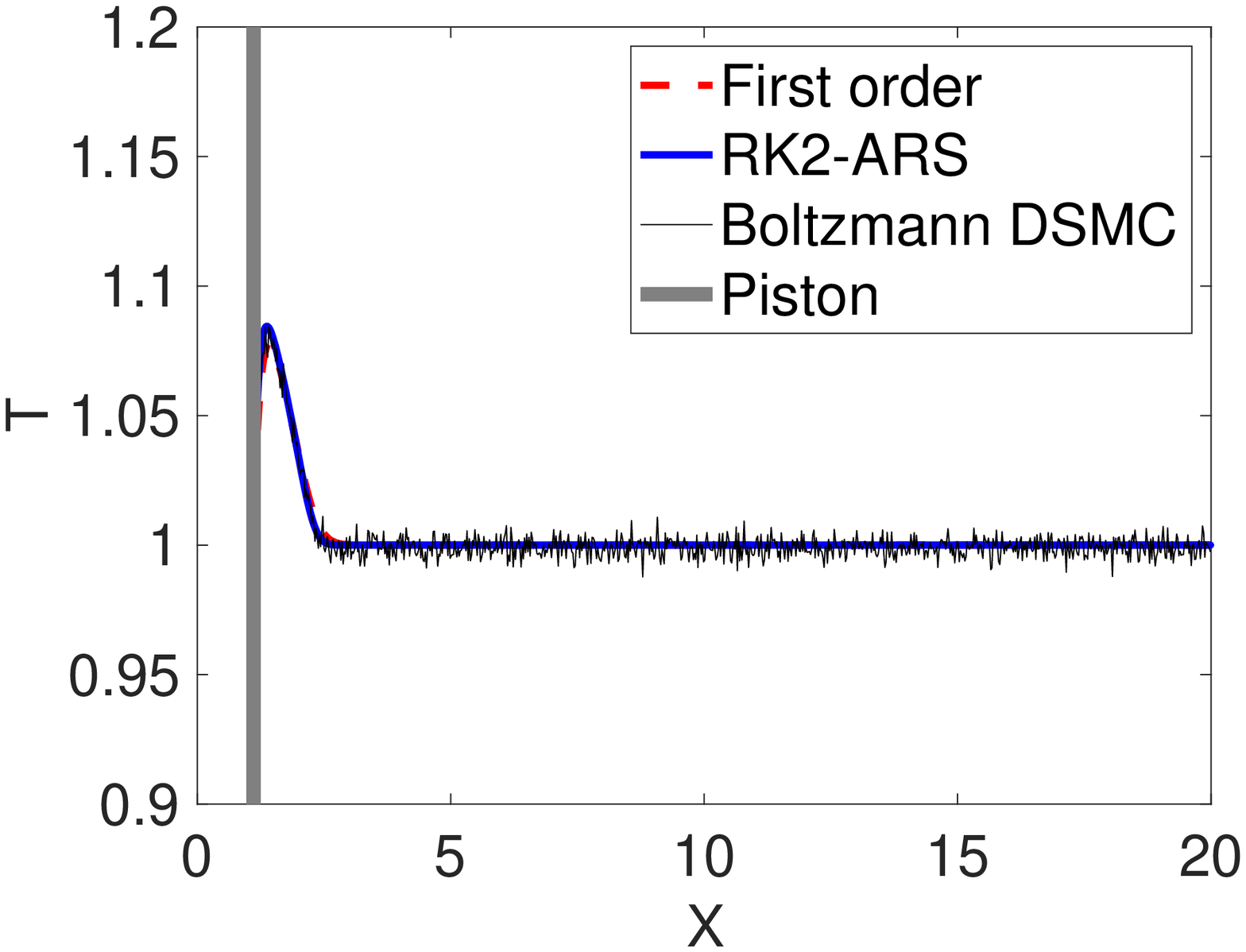}
\includegraphics[keepaspectratio=true, width=.329\textwidth]{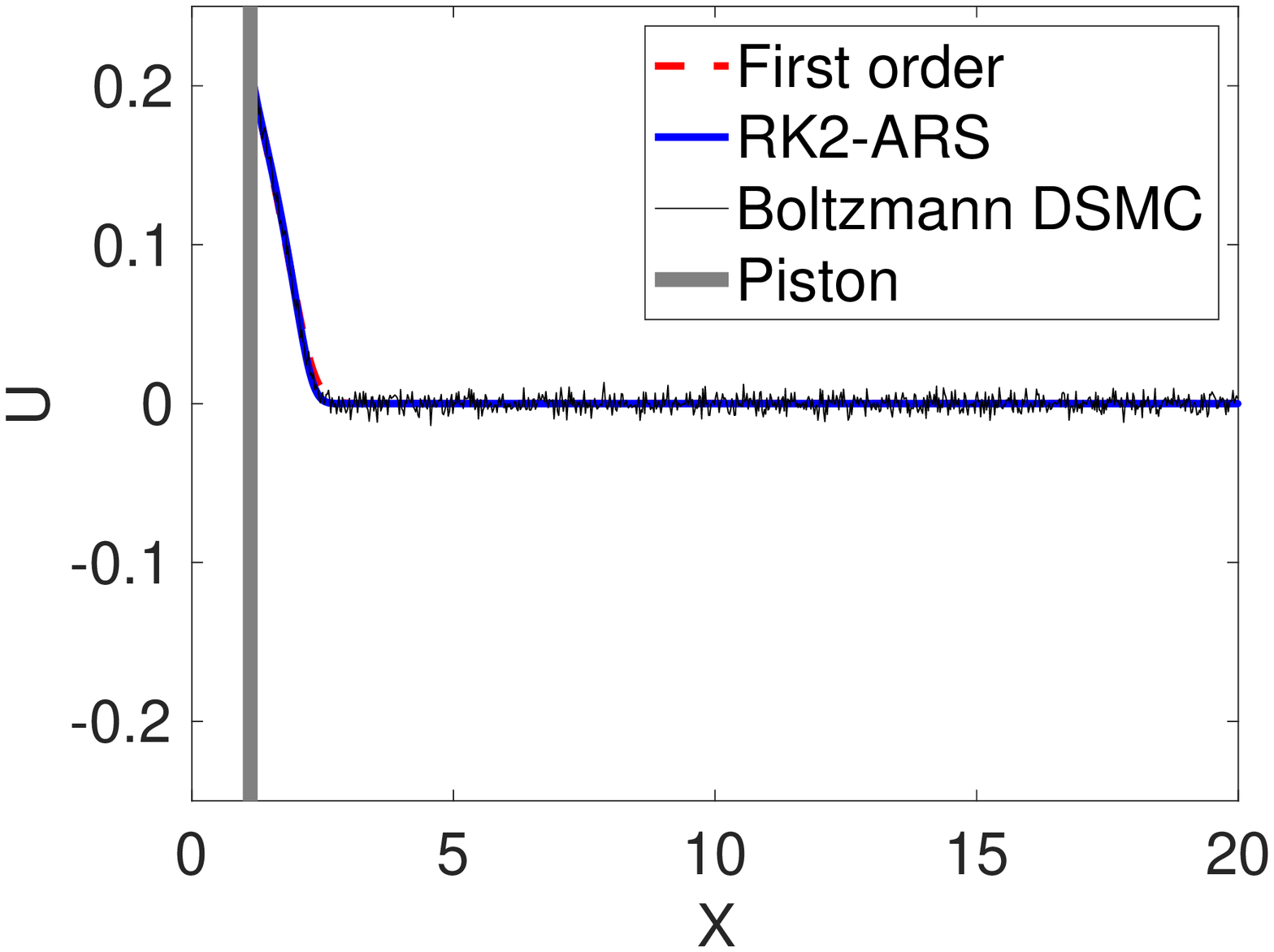}
\caption{Example 3: Moving piston. Comparison of ALE and DSMC methods at time $t=1$.}
\label{moving_piston1}
\end{figure}     

 \begin{figure}[!t]
\includegraphics[keepaspectratio=true, width=.329\textwidth]{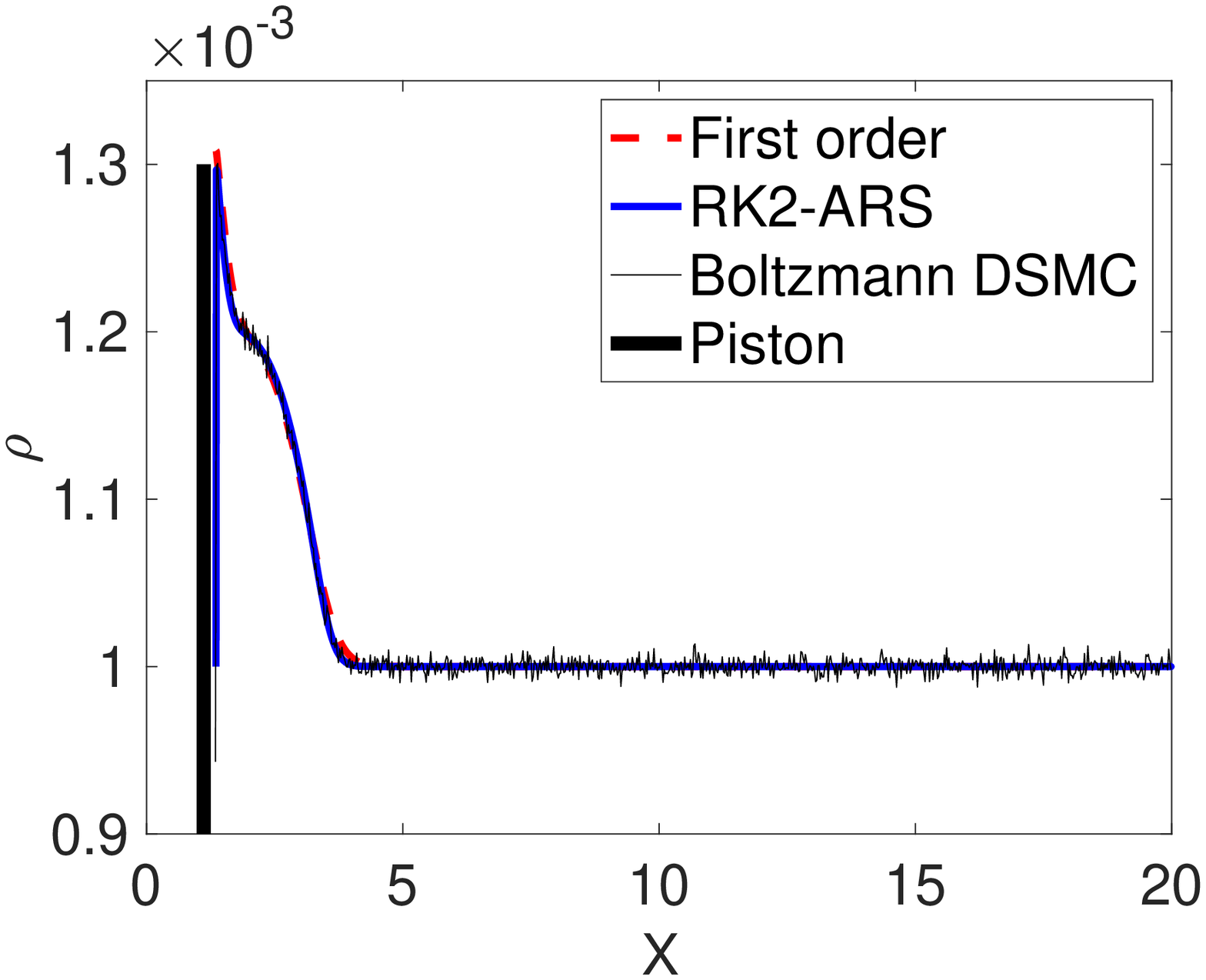}
\includegraphics[keepaspectratio=true, width=.329\textwidth]{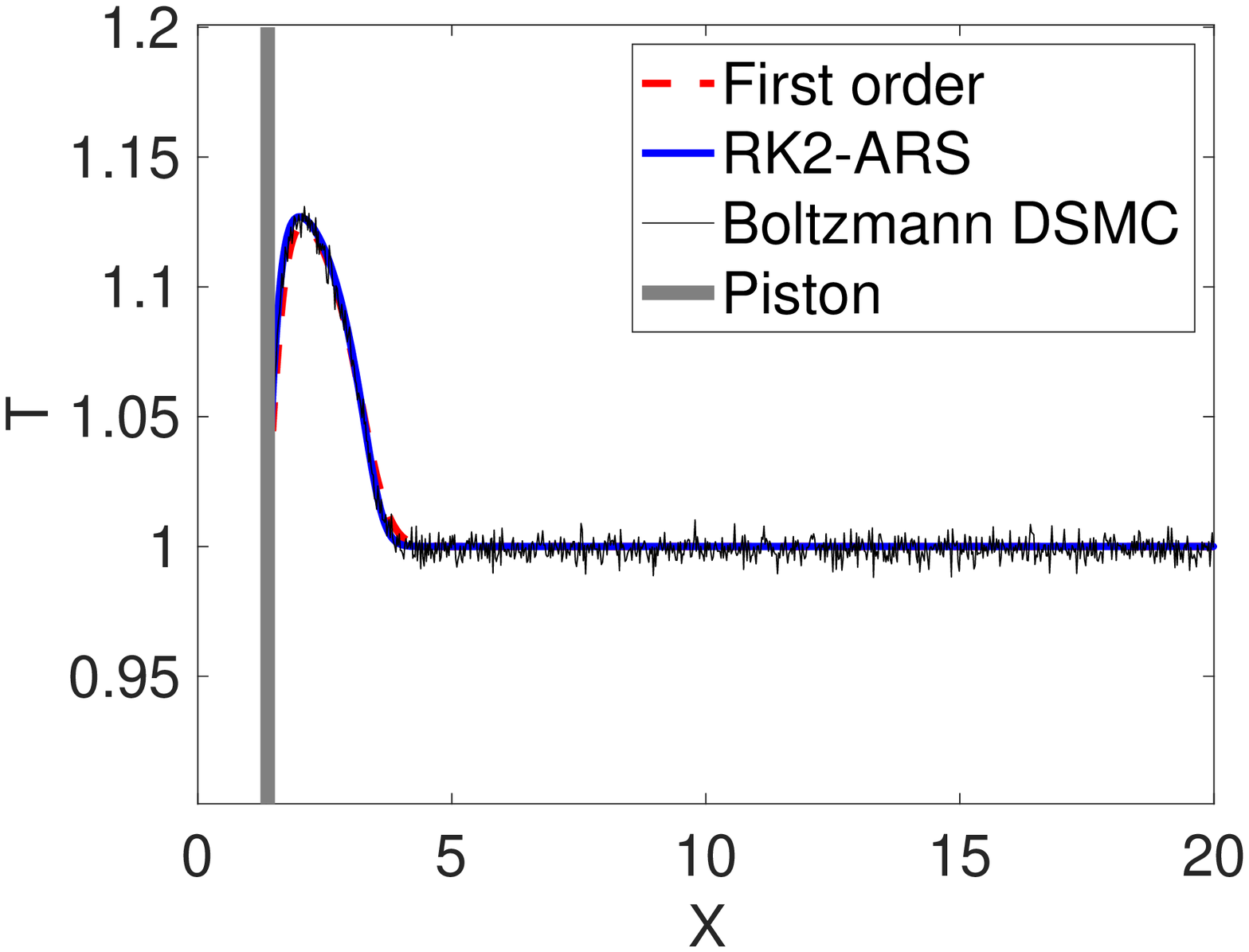}
\includegraphics[keepaspectratio=true, width=.329\textwidth]{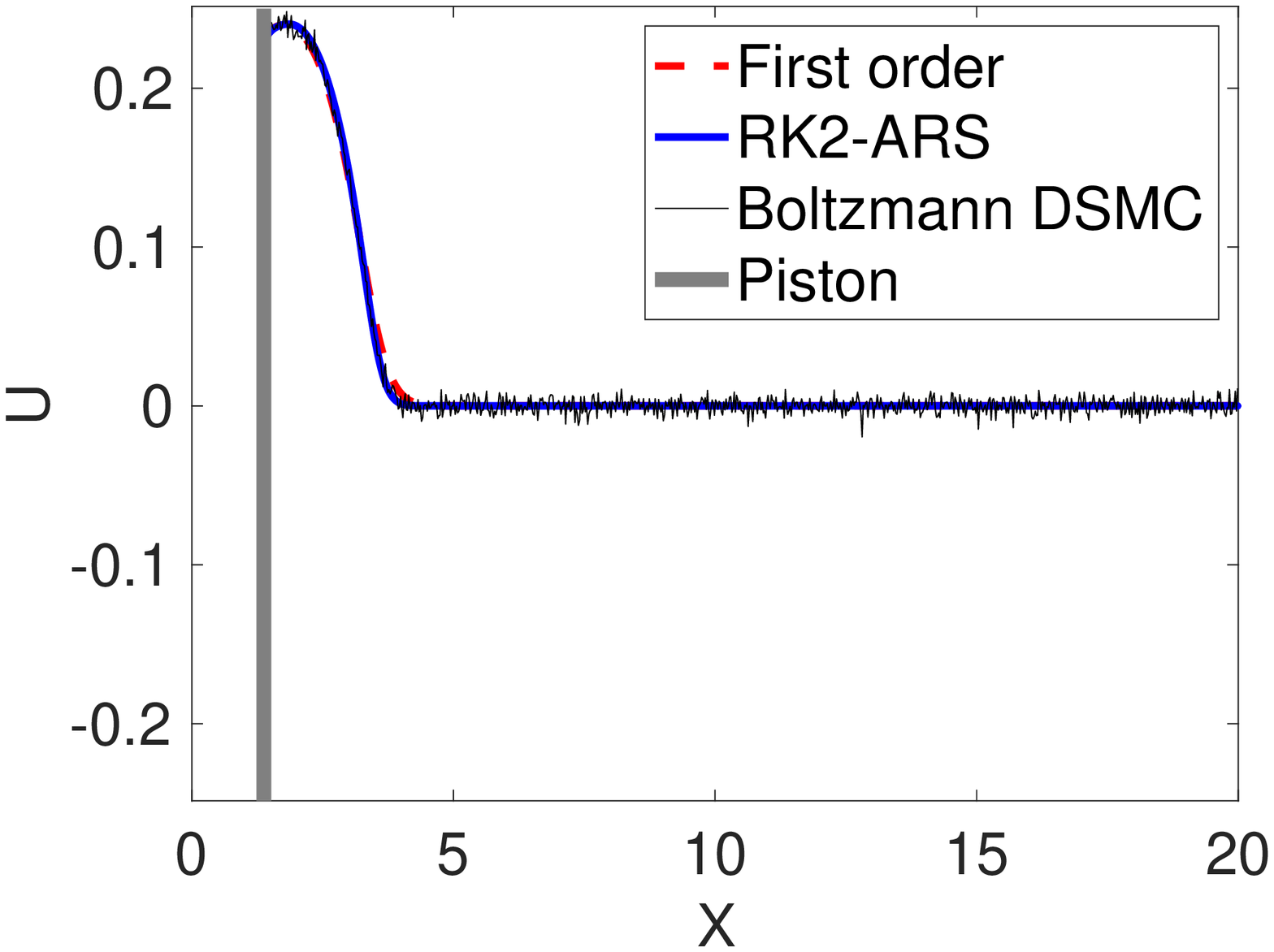}
\caption{Example 3: Moving piston. Comparison of ALE and DSMC methods at time $t=2$.}
\label{moving_piston2}
\end{figure}     

 \begin{figure}[!t]
\includegraphics[keepaspectratio=true, width=.329\textwidth]{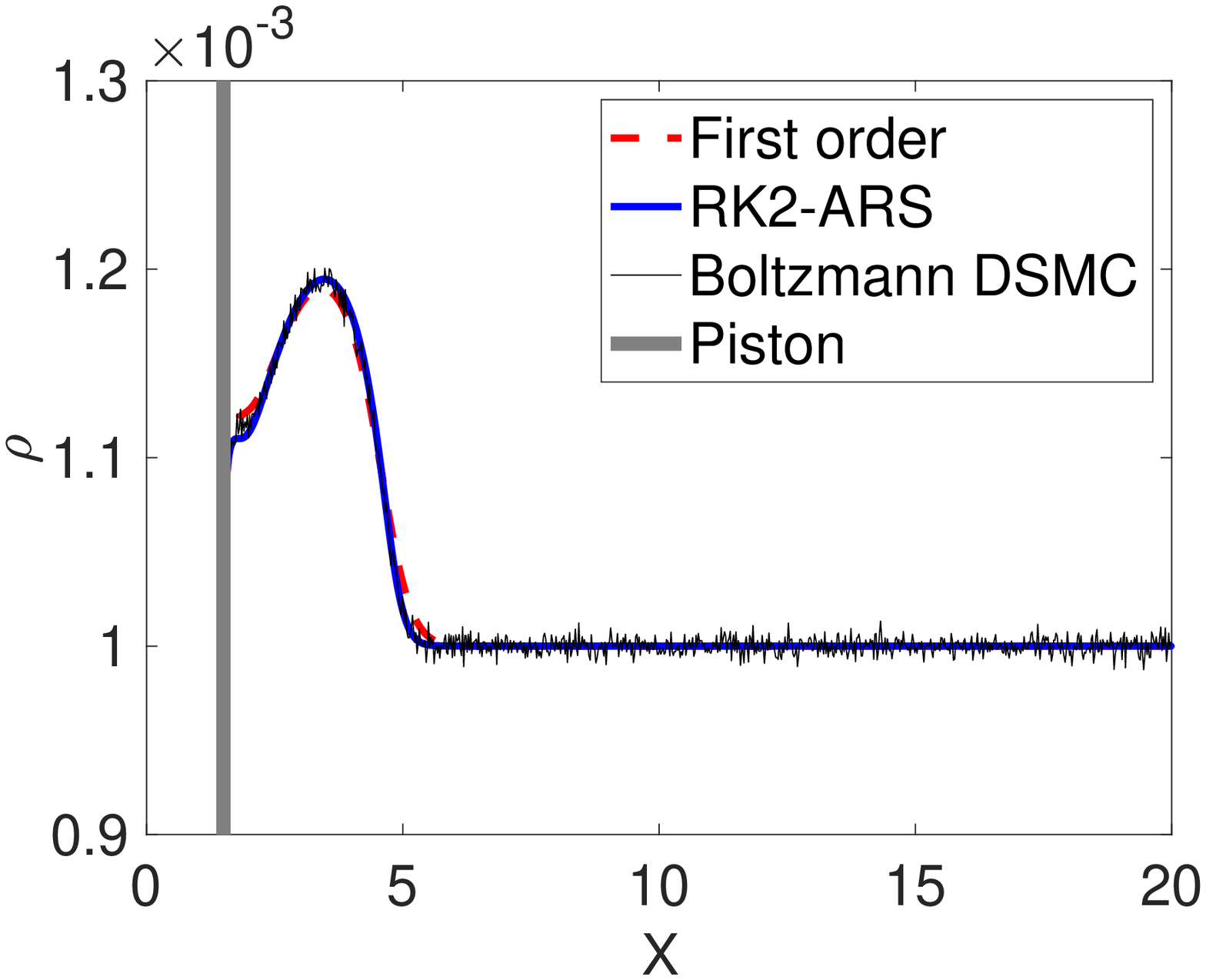}
\includegraphics[keepaspectratio=true, width=.329\textwidth]{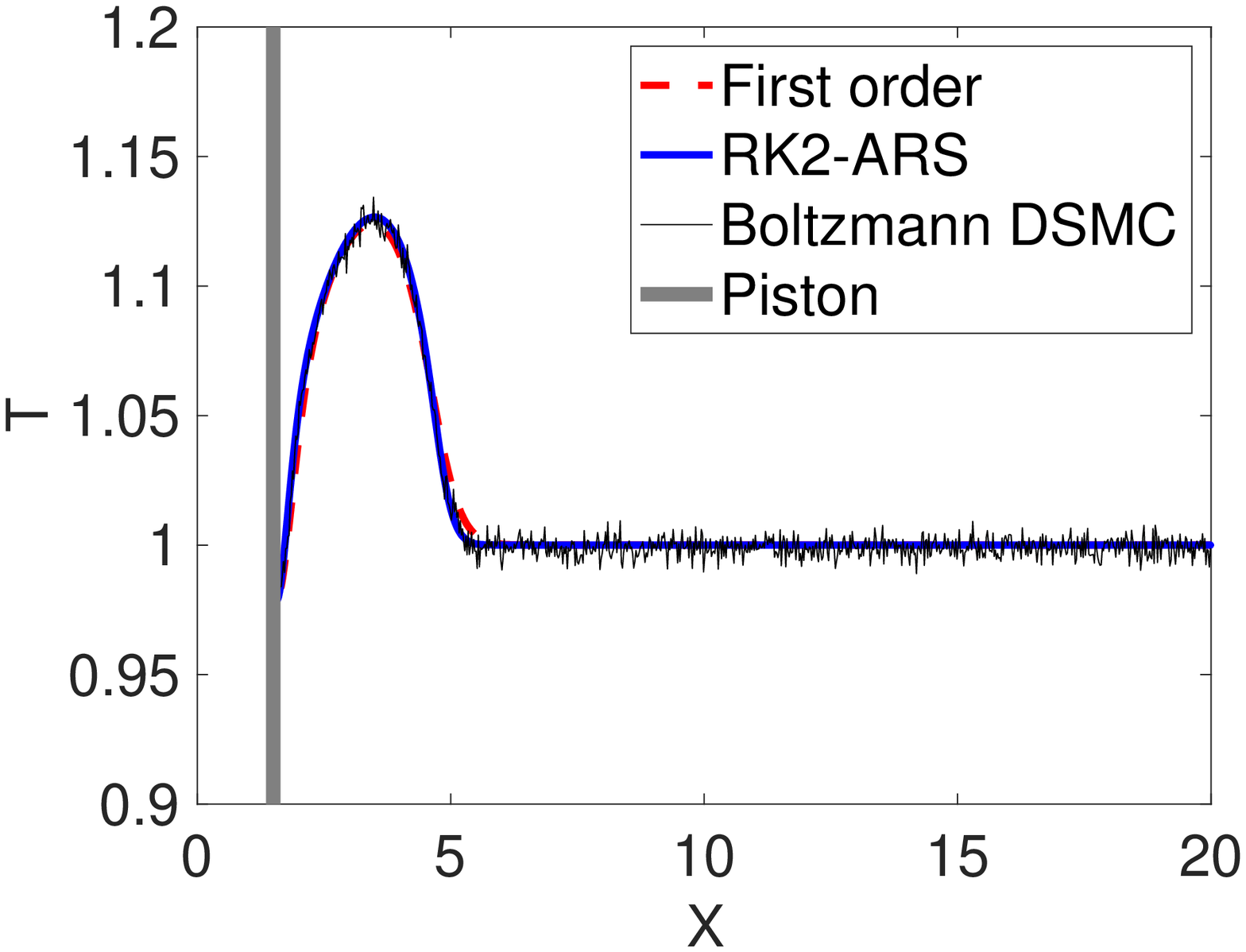}
\includegraphics[keepaspectratio=true, width=.329\textwidth]{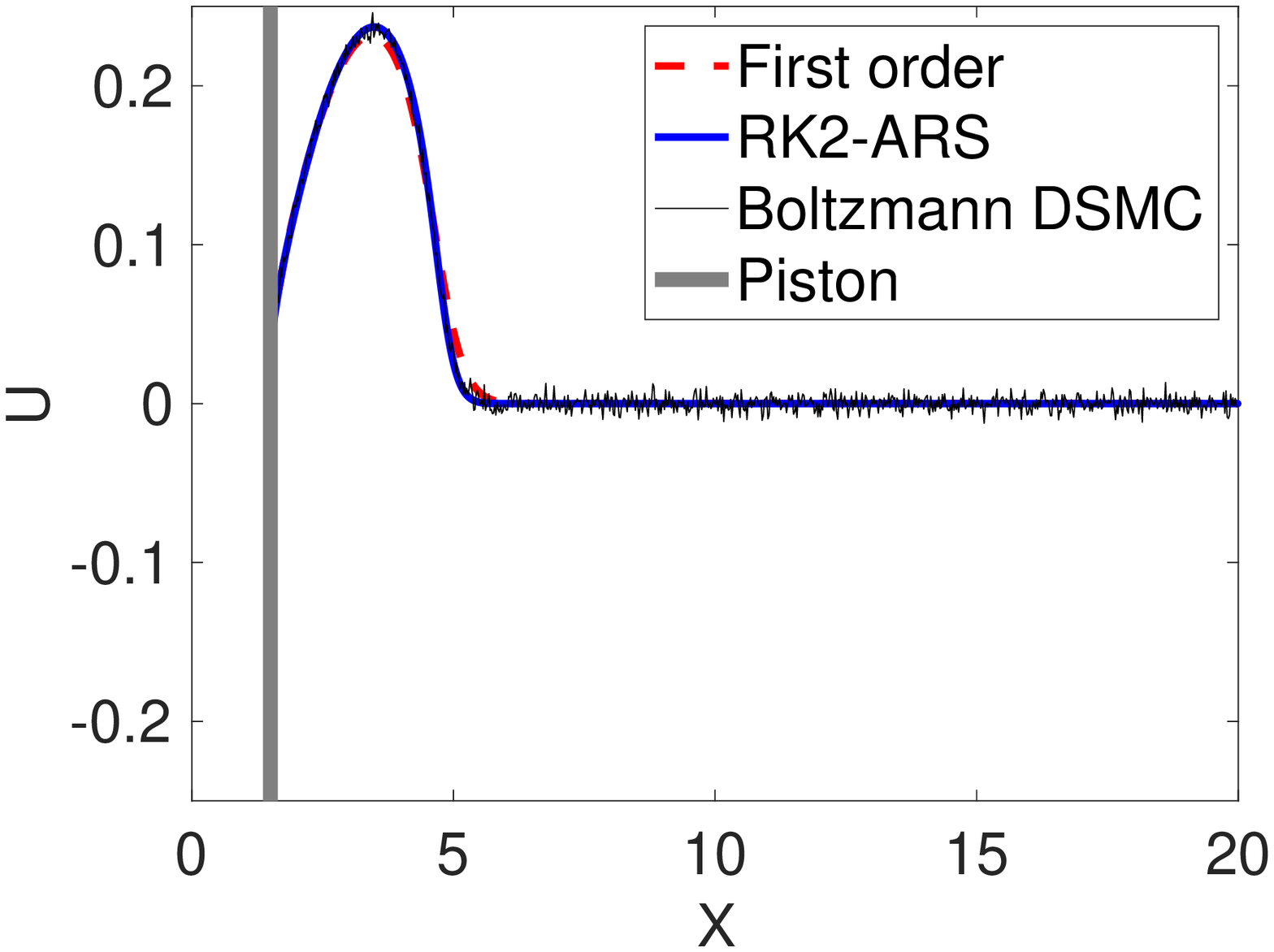}
\caption{Example 3: Moving piston. Comparison of ALE and DSMC methods at time $t=3$.}
\label{moving_piston3}
\end{figure}     

 \begin{figure}[!t]
\includegraphics[keepaspectratio=true, width=.329\textwidth]{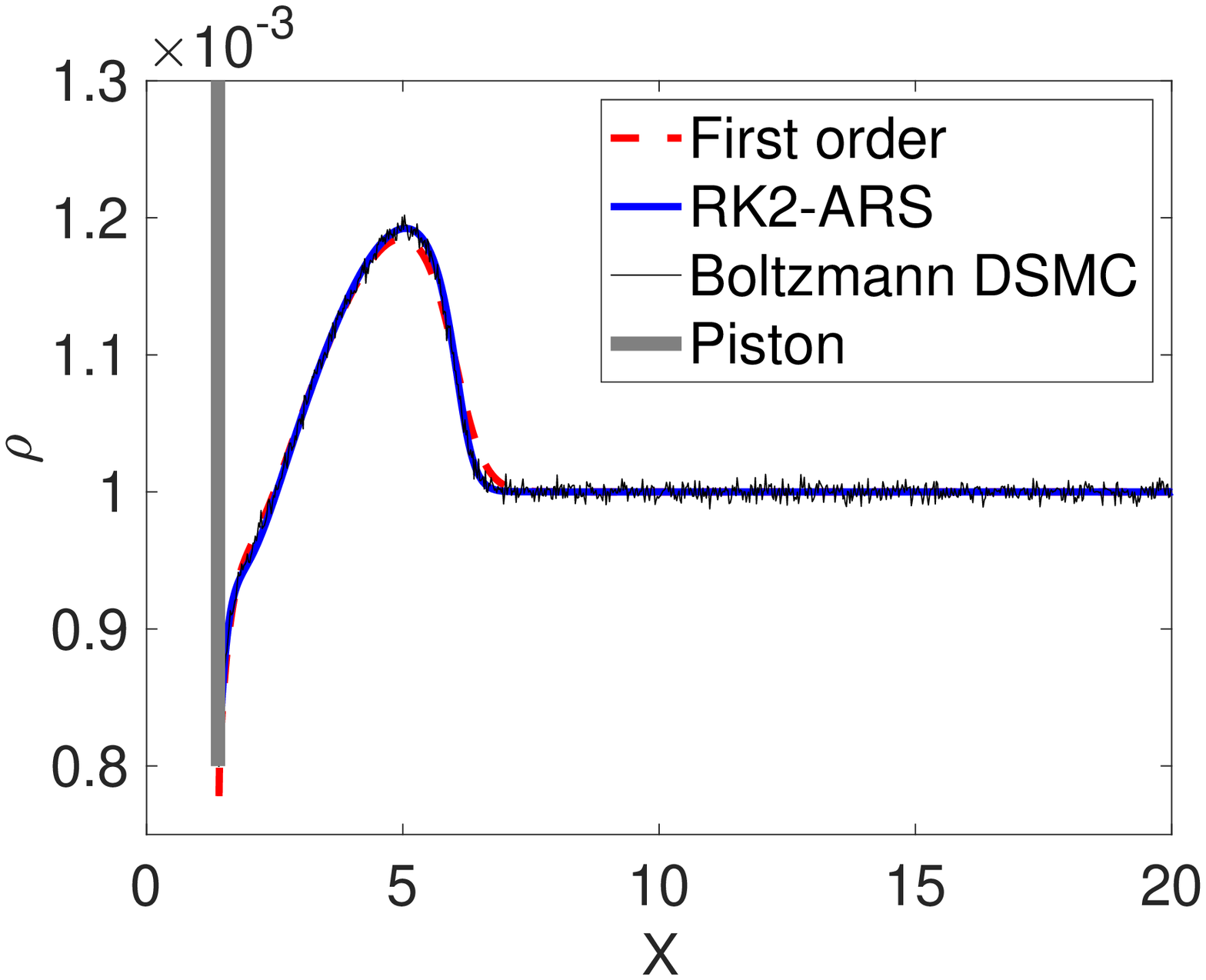}
\includegraphics[keepaspectratio=true, width=.329\textwidth]{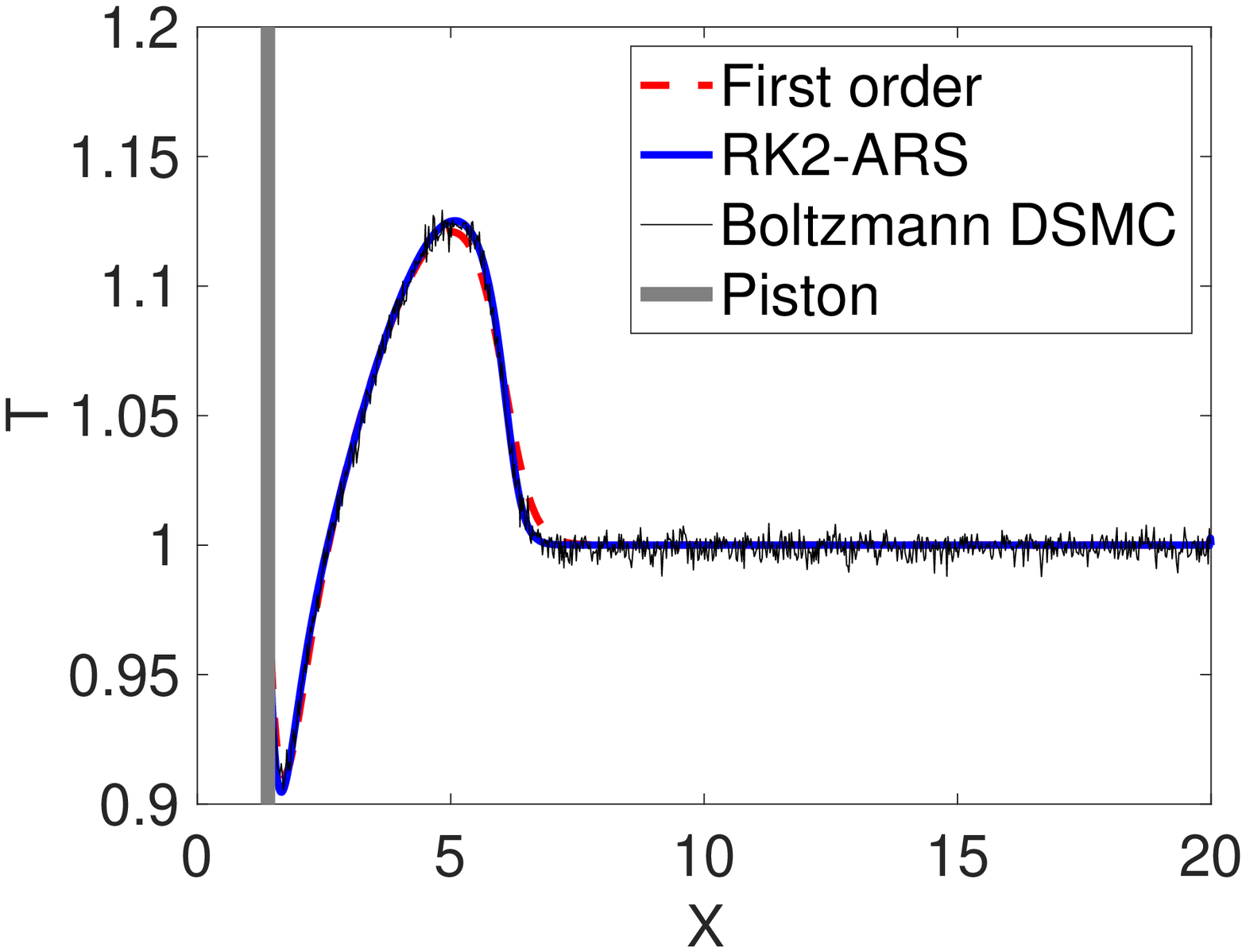}
\includegraphics[keepaspectratio=true, width=.329\textwidth]{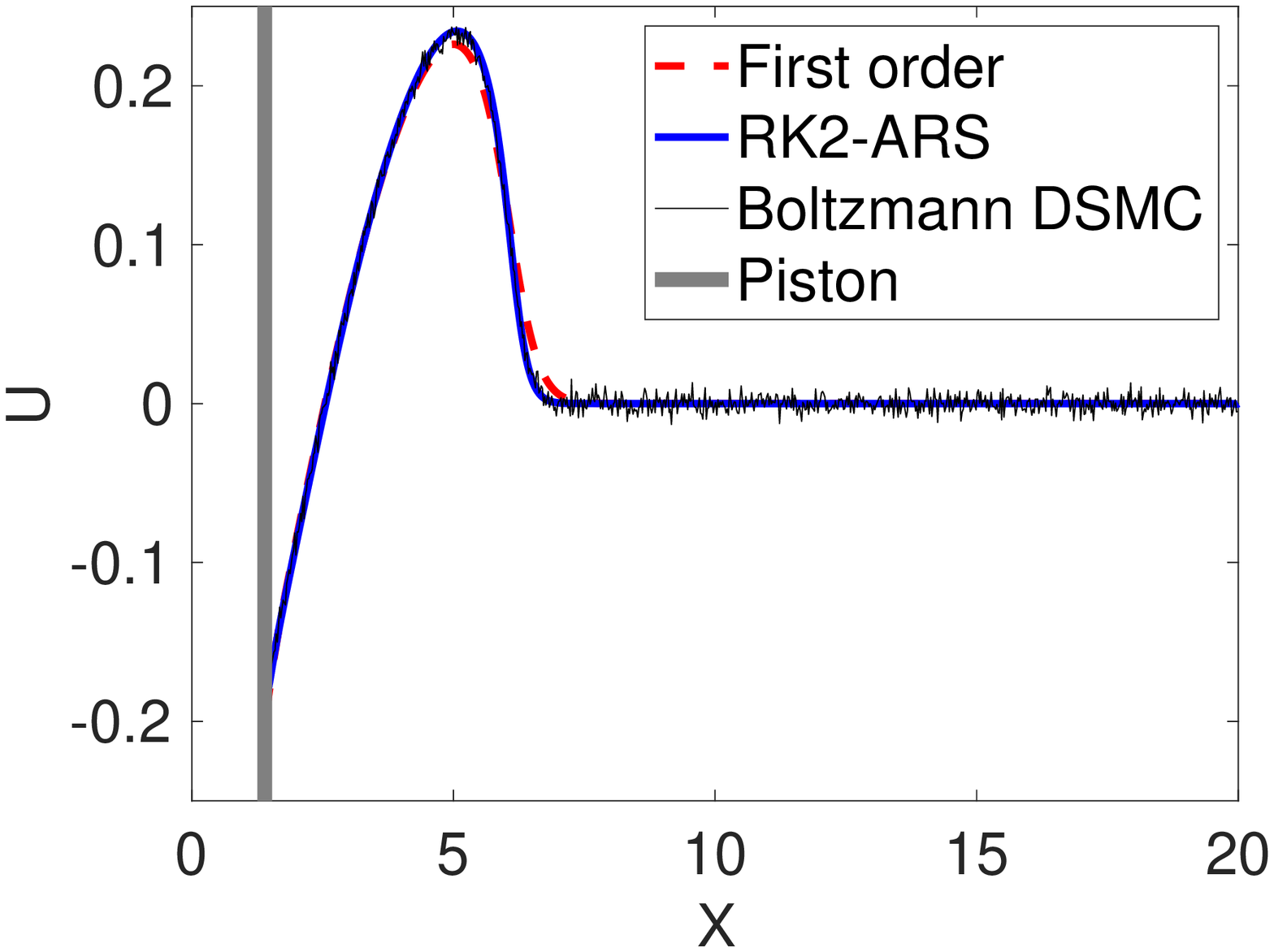}
\caption{Example 3: Moving piston. Comparison of ALE and DSMC methods at time $t=4$.}
\label{moving_piston4}
\end{figure}     

\begin{figure}[!t]
\includegraphics[keepaspectratio=true, width=.329\textwidth]{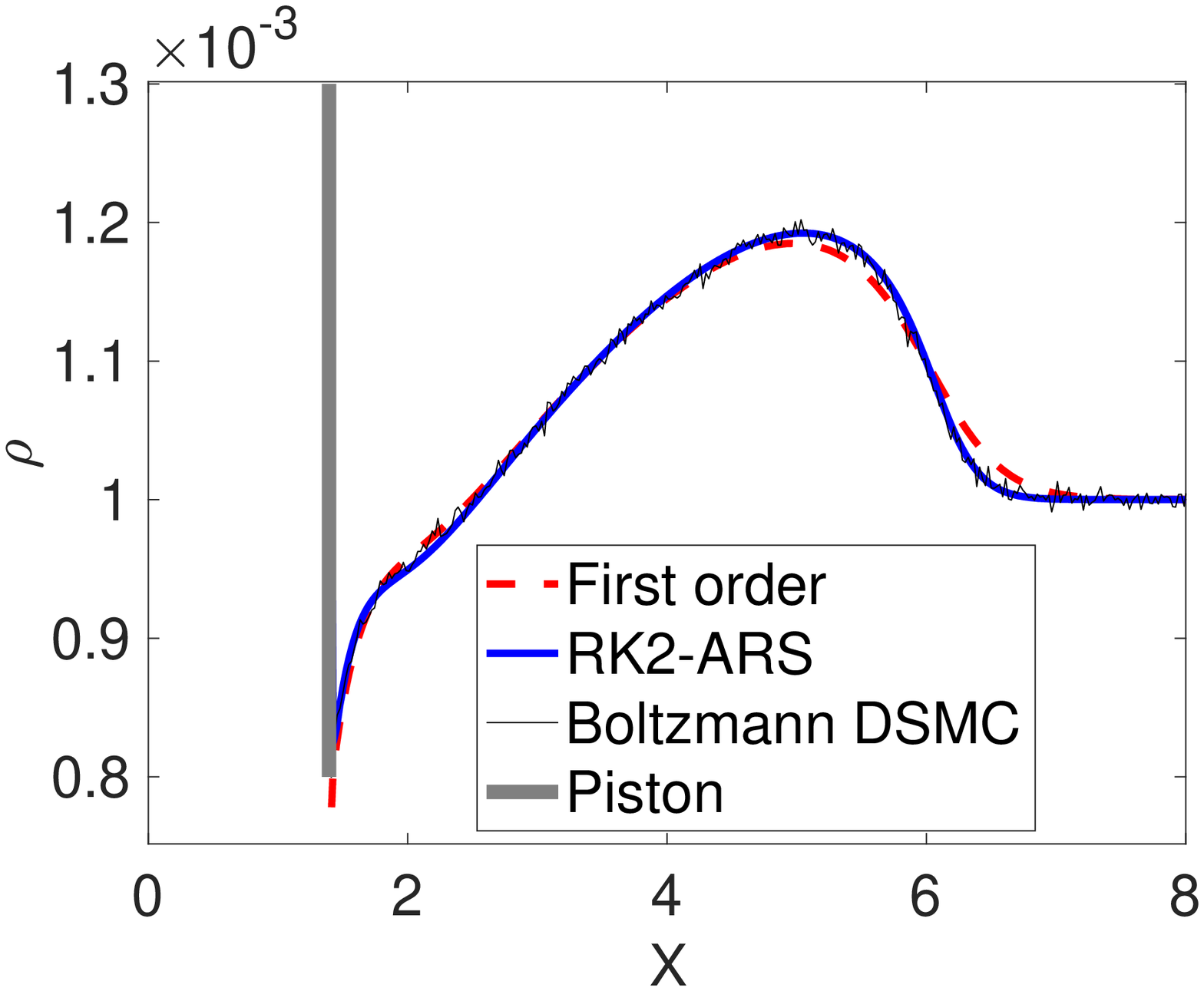}
\includegraphics[keepaspectratio=true, width=.329\textwidth]{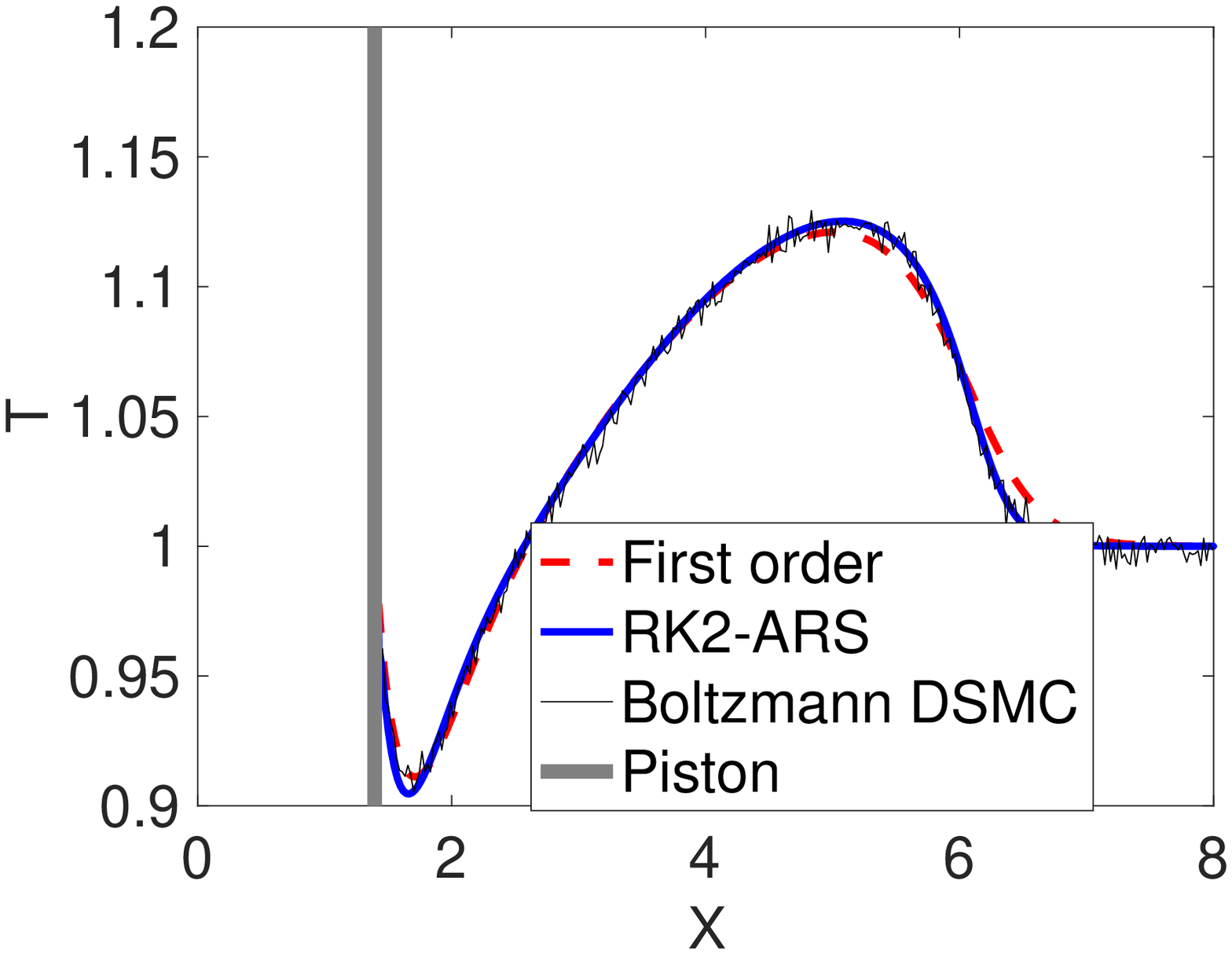}
\includegraphics[keepaspectratio=true, width=.329\textwidth]{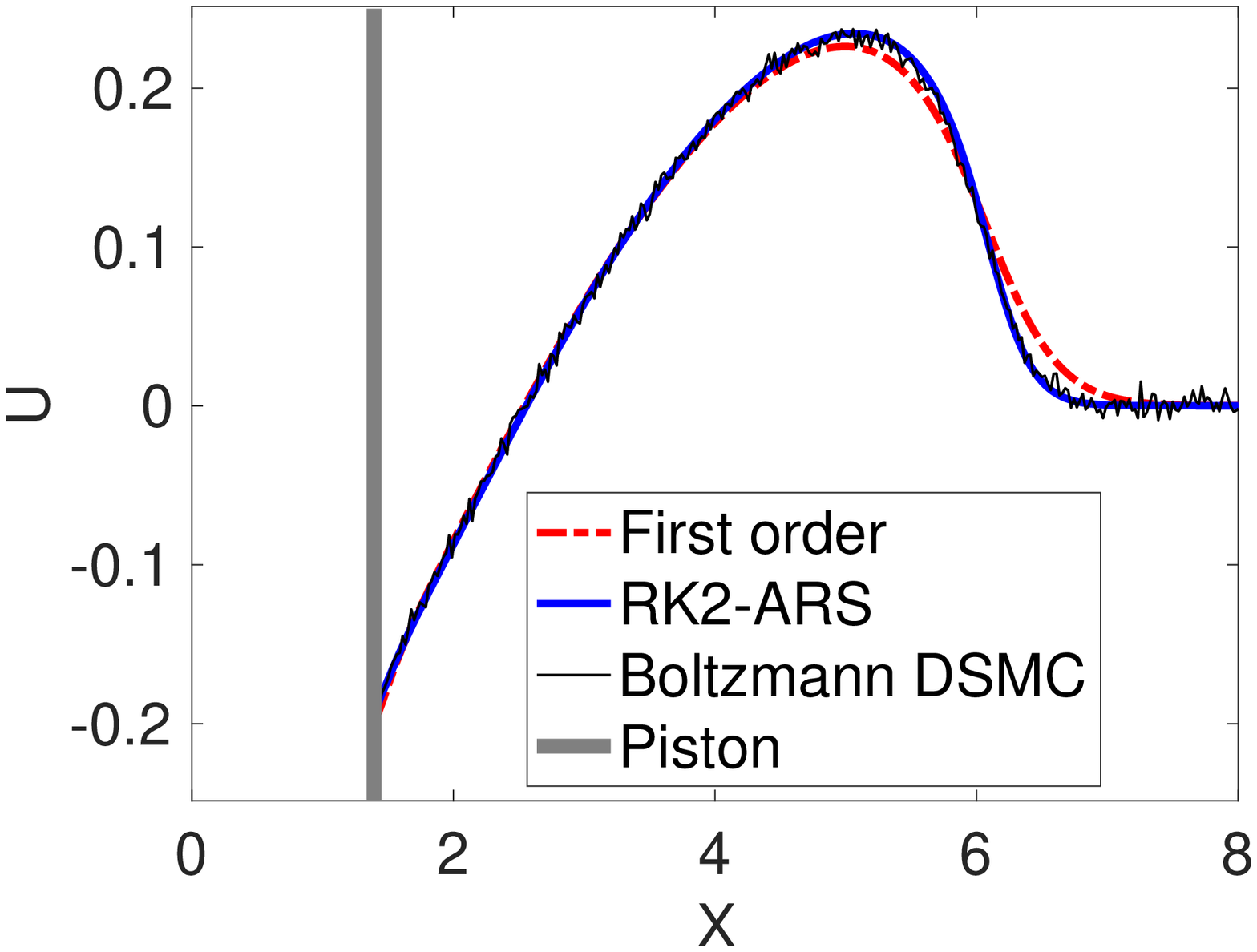}
\caption{Example 3: Moving piston. Zoom of the data obtained from ALE and DSMC methods at time $t=4$.}
\label{moving_piston4_zoom}
\end{figure}

 \subsection{Example 4: Movement of a plate with pressure differences}

We consider a computational domain as described in Figure \ref{Fig-Actuator} with $L=1$ and $l=0.1$. 
\begin{figure}[ht]
	\begin{center}
		\begin{tikzpicture} [scale=0.8]
		\draw [fill=black,opacity=0.4] (5,0) rectangle (7,2);
		\draw [fill=gray,opacity=0.4] (0,0) rectangle (7,2);	
		\draw [fill=gray,opacity=0.4] (7,0) rectangle (12,2);
		  \node[color=black] at (3,	1) {\bf Gas };
		   \node[color=black] at (9,	1) {\bf Gas };
%		\draw[line width=3pt,color=black,opacity=1.0] (0,0) -- ( 10,0);
%		\draw[line width=3pt,color=black,opacity=1.0] (0,2) -- ( 10,2);
%		\draw[line width=3pt,color=black,opacity=1.0] (0.08,0) -- ( 0.08,2);
%		\draw[line width=3pt,color=black,opacity=1.0] (9.94,0) -- ( 9.94,2);
		
		% Scale
		\draw[line width=1pt,color=black,opacity=0.5] (0,-0.7) -- ( 12,-0.7);
		\draw[line width=1pt,color=black,opacity=0.5] (0.0,-0.5) -- ( 0.0,-0.9);
		\draw[line width=1pt,color=black,opacity=0.5] (12,-0.5) -- ( 12,-0.9);
		\draw[line width=1pt,color=black,opacity=0.5] (6,-0.5) -- (6,-0.9);
		
		% Text for scale
		\draw[line width=0.4pt,<->] (5,-0.2) -- (7,-0.2);
		\node[color=black] at (6,-0.4) {$2l$ };
		\node[color=black] at (0.08,-1.1) {$-(L+l)$ };
		\node[color=black] at (12,-1.1) {$(L+l)$ };
		\node[color=black] at (6,-1.1) {$0$ };
		
		% Text for Plate
		\draw[line width=0.4pt,<->] (5,2.2) -- (7,2.2);
		\node[color=black] at (6,2.5) {Moving plate };
		% Text for walls
		\node[color=black] at (-0.35,1) {$T_0$ };
		\node[color=black] at (4.7,1) {$T_0$ };
		\node[color=black] at (7.3,1) {$T_w$ };
		\node[color=black] at (12.43,1) {$T_w$ };		
		\end{tikzpicture}
	\caption{Example 4: Movement of a plate. Schematic view of a piston separating two subdomains with different temperature.
	% Like in piston problem, the black circles are active grids and grey circles are inactive grids.
	}
 \label{Fig-Actuator}
  \end{center}
  \end{figure}
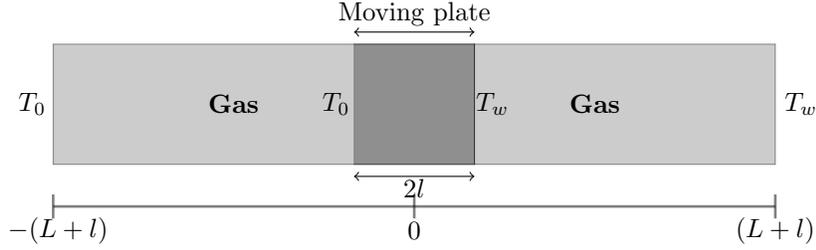
Initially, the center of mass of the plate is located at $X_c = 0$. The gas and the plate are at rest. This problem has been studied in \cite{DM, TKR}. We reconsider it  as  benchmark problem since an  analytical expression 
%\giovanni{of what? Perhaps this ``analytical solution'' should be better explained} 
is available for  the equilibrium state. Using 
SI units, the initial temperature is $T_0=270$, gas constant $R=208$ and the initial pressures $P_0$ are the same on both sides of the plate  and are  equal to $0.0386$.  The initial density $\rho_0$ is obtained from the equation of state. The initial Knudsen number is $0.08$ based on the characteristic length $2L$ and the relaxation time is fixed as $\tau = 5.398\cdot 10^{-4}$. There are four boundary points, two are at the boundary of the domain and 
two are at the left and right end of the plate. The interior grid points are  initialized with the spacing $\Delta x = 2.2/200$ on the left and right of the plate. No 
grid points are initialized on the plate.  The neighbor radius is given by $ \Delta x=0.35 h$ and the constant time step $\Delta t = 2\cdot 10^{-6}$ is considered. 
We prescribe a higher temperature $T_w= 330$ on the right side of plate and on the right boundary of the computational domain. On the left boundary of the plate and on the left boundary of the computational domain the temperature is fixed  to $T_0$. Due to the high temperature on the right wall, the pressure on the right hand side starts to increase and  the plate starts to move to the left hand side. The density of the plate is $10$ times larger than the density of the gas. This means, the mass of the plate is equal 
to $M = 3.4366\cdot 10^{-5}$. The motion of the plate is computed from the Newton-Euler equations, where only a translational force is computed for  the one dimensional case. Since the plate has two opposite normals $\pm 1$, from equations (\ref{force_torque_bgk}) and (\ref{pressure_tensor}) the total force is given as the difference of pressure 
\begin{equation}
{F} = (\varphi_{\rm left} - \varphi_{\rm right}) A,
\end{equation}
where $A$ is the area of the plate and $\varphi = \int_{\mathbb R}(v-U)^2 g_1 dv$. The plate starts oscillating and finally reaches the equilibrium position \cite{DM}
\begin{equation}
\label{0.1}
x_{\rm equi} = L\frac{(T_0-T_w)}{(T_0+T_w)} = -0.1.
\end{equation}

We have compared the  dynamics of the plate obtained from the ALE method with first and second order ARS schemes
  with a Boltzmann solution using the DSMC method. 
We observe that the oscillation of the plate obtained from both methods match. The simulations are performed up to the final time $t=0.6$ and the piston 
already reached the equilibrium at this time, see Figure \ref{moving_plate1}. At  the final time the simulated equilibrium position obtained from the first order method is $-9.639 \cdot 10^{-3}$ and one given by the second order method is $-9.963 \cdot 10^{-3}$ 
compared to the analytical solution which gives a value of $-0.1$, see (\ref{0.1}).
This yields  an error of  $3.7\%$ and $0.37\%$, respectively. 
\begin{figure}[!t]
	\centering
 	\includegraphics[keepaspectratio=true, width=0.45\textwidth]{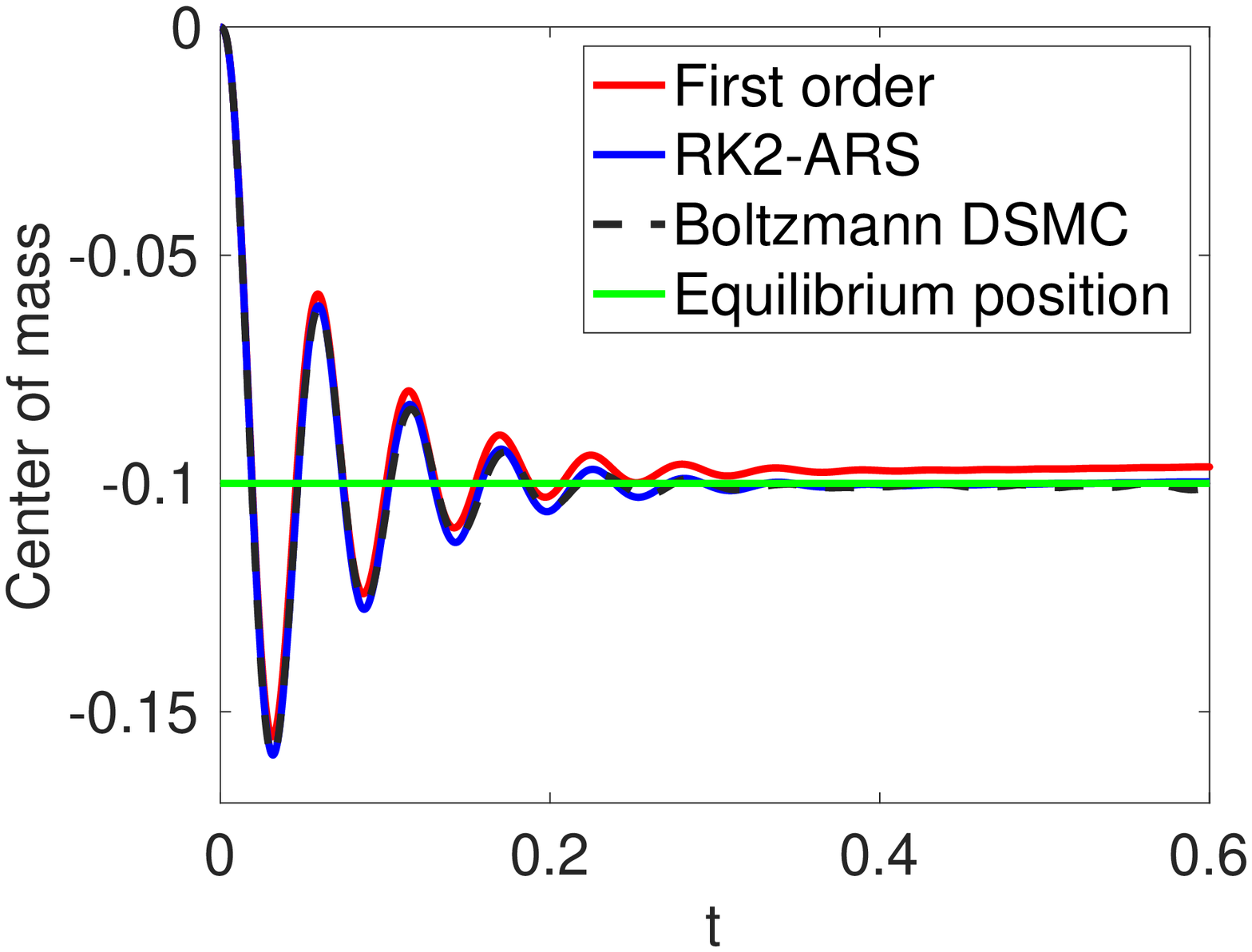}
        \includegraphics[keepaspectratio=true, width=0.45\textwidth]{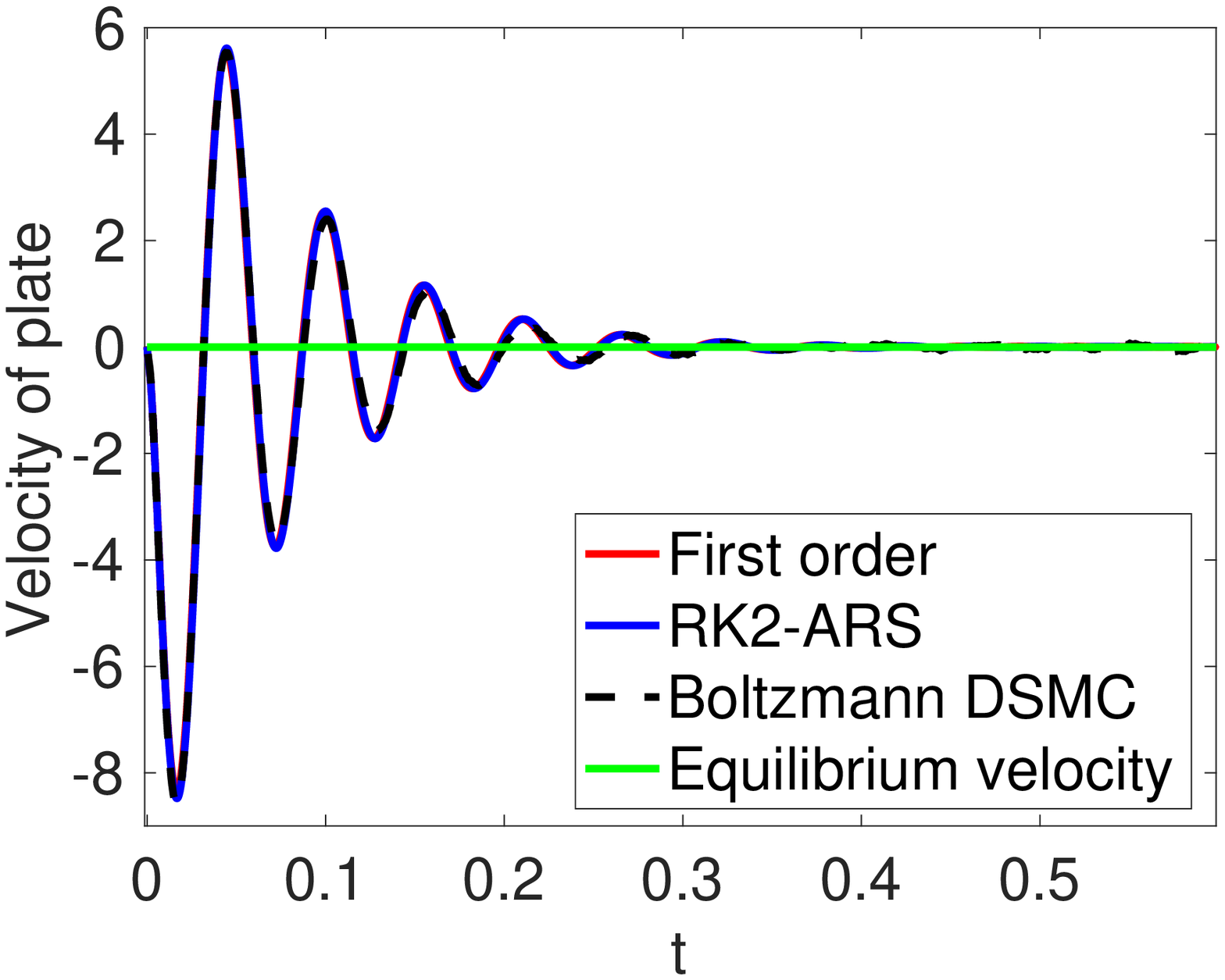}
   	\caption{Example 4: Movement of a plate.  Comparison of position and velocity vs time of piston obtained from ALE and DSMC method}
	\label{moving_plate1}
\end{figure}
                  
%%%%%%%%%%%%%%%%%%%%%   2D cases     %%%%%%%%%%%%%%%%
\subsection{Example 5: The 2D-BGK model with smooth solution}

For the convergence study we consider  the BGK model with two-dimensional  space and velocity domain for short time for a situation extending the one in   section  \ref{1d_conv_study} to 2-D. 
The computational domain 
is  $\Omega = [-1, 1]\times [-1, 1]$. The initial distribution is again the Maxwellian distribution and is given by 
\[
f(0,x,v) = \frac{\rho_0}{(2 \pi R T_0)^{3/2}}  \exp {\left(-\frac{({v}- {U_0})^2}{2RT_0}\right)}
\]
with $\rho_0 = 1, T_0 = 1, R = 1$ and $U_0=(U_0^{(x)},U_0^{(y)},0)$ with
\begin{eqnarray*}
U_0^{(x)} &=& \frac{1}{\sigma} \left(\exp(\left(-(\sigma \sqrt{(x-0.2)^2+y^2} - 1)^2\right ) - 2 \exp\left(-(\sigma \sqrt{(x+0.2)^2+y^2} - 1)^2\right) \right )
\\
U_0^{(y)} &=& \frac{1}{\sigma} \left(\exp(\left(-(\sigma \sqrt{x^2+(y-0.2)^2} - 1)^2\right ) - 2 \exp\left(-(\sigma \sqrt{x^2+(y+0.2)^2} - 1)^2\right) \right ),  
\end{eqnarray*}
where $\sigma = 10$.
We have chosen again $\tau = 10^{-5}$.
Far field boundary conditions are applied on the boundaries with initial 
density, temperature and zero mean velocities. 
In order to perform the convergence study the time integration is carried out up to time $t=0.0208$, where the solution is still smooth. Different numbers of grid points are considered depending on the size of $h$. The initial grid spacing is $\Delta x = 0.4~h$. The reference solution is the solution obtained from a grid  with $h = 0.013$, which corresponds to an initial number of grid points equal to $148996$. For the reference solution we use a time step equal to $\Delta t = 2.6\cdot 10^{-5}$, which corresponds to a CFL condition with constant $0.5$. We refer to subsection \ref{1d_conv_study} for a discussion of the CFL condition used here.
%\giovanni{We should specify better what we mean by CFL and we should explain that with ALE in principle we should have a better stability condition and use larger time steps. However, I suspect this would be true if, as we suggest somewhere in the paper, we centered the velocity grid in $U$. If we use a fixed velocity from $-v_{\rm max}$ to $+v_{\rm max}$ actually the CFL condition may become MORE restrictive, since we have to compute $\Delta t$ imposing that 
%\[
%	\max_{x,v} \frac{|v-U|\Delta t}{h}<CFL
%\]
%In any case, we have to better specify what we do.}
%
%\axel{What is understood by CFL  is shortly exlained  at the end of 4.1. Essentially a velocity 2 vmax is used for the CFL definition, which means that $U$ should not be larger than vmax. This is in all cases fulfilled. However, the estimate is very 
%coarse and could be certainly improved as you are stating it. I used your sentences at the end of 4.1 to explain this  a bit better}
This CFL number is also used for all other grid-sizes.

The convergence rate is determined by interpolating the temperature on $100$ grid points along $y=0$  for all grid sizes. In Figure \ref{2d_temp_along_y_0} we have plotted the temperature  obtained from the first order scheme and the ARS(2,2,1) scheme.
Again, the ARS(2,2,2) scheme gives equivalent results. We note that we gain some computation time by using the ARS(2,2,1) scheme  due to the additional function
evaluations in the ARS(2,2,2) scheme.
%\giovanni{This is not 100\% clear to me. Why ARS(2,2,1) requires less computation than ARS(2,2,2)?}
%\axel{Since it requires less function evaluations/computations of derivatives (31) compared to (37)}

 \begin{figure}[!t]
 \centering
 \includegraphics[keepaspectratio=true, width=0.45\textwidth]{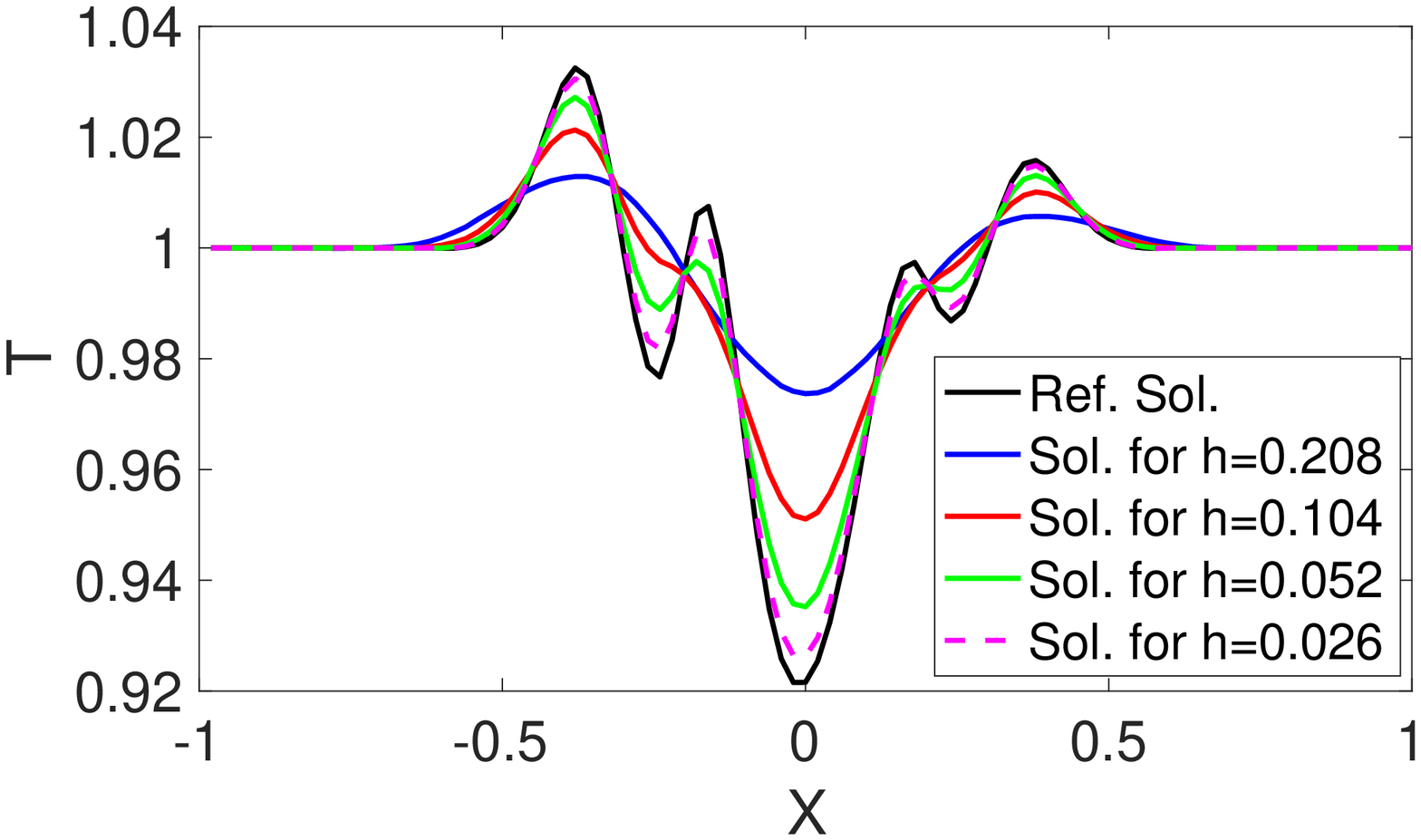}
 \includegraphics[keepaspectratio=true, width=0.45\textwidth]{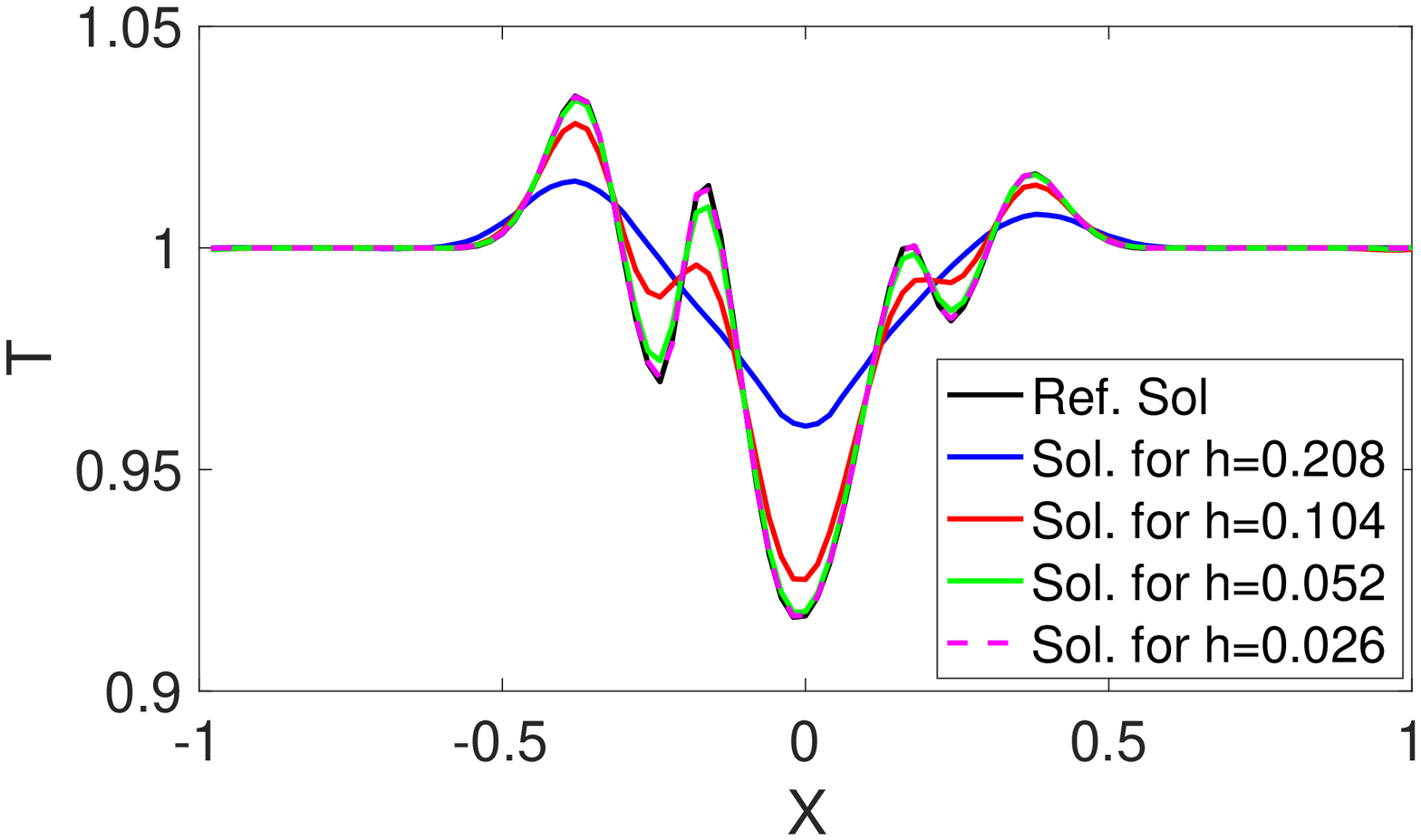} 
 \caption{Example 5: 2D smooth solution. Temperature at $t = 0.0208$ along $y=0$ for different $h$ obtained from the first order (left) and second order (right) schemes. }
 \label{2d_temp_along_y_0}  
 \end{figure}

In Tables \ref{2d_convergence_table1} and \ref{2d_convergence_table2} we have presented the corresponding errors and 
the rate of convergence. It can be observed that the rates of convergence are as expected  for the corresponding  schemes.
  
  \begin{table}
  \begin{center}
  \begin{tabular}{|r|r|r|r|r|r|r|}
  \hline
  $\Delta t$ &	$h $ &  $N_x$ & $L^1$-error & Order& $L^2$ error & Order  \\ \hline 
  $ 4.16\cdot10^{-4}$ & $0.208$ &676 &$ 5.57\cdot10^{-2}$   &    $--$ &           $2.10\cdot10^{-2}$ &$--$\\  
  $ 2.08\cdot 10^{-4}$ &	$0.104$ & 2500&$1.01\cdot 10^{-2}$  &   $0.64$ &       $1.30\cdot 10^{-2}$ &$0.69$\\ 
  $ 1.04\cdot 10^{-4}$ &	$0.052$ &9604 &$5.18\cdot 10^{-3}$  &   $0.97$ &       $6.62\cdot 10^{-3}$ &$0.98$\\ 
  $ 5.20\cdot 10^{-5}$ &	$0.026$ &37636 &$1.93\cdot10^{-3}$ &   $1.43$ &       $2.47\cdot 10^{-3}$ &$1.42$\\  
  \hline
  \end{tabular}
  \caption{Example 5: 2D smooth solution. Convergence of temperature at time  $t =0.0208$ from the first order scheme.}
  \label{2d_convergence_table1}
  \end{center}
  \end{table}
  
  %%%%
   \begin{table}
   	\begin{center}
   		\begin{tabular}{|r|r|r|r|r|r|r|}
   			\hline
   			$\Delta t$ &	$h $ &  $N_x$  & $L^1$-error & Order& $L^2$ error & Order  \\ \hline 
   			$ 4.16\cdot10^{-4}$ & $0.208$ &676 &$ 1.48\cdot10^{-2}$   &    $--$ &           $1.89\cdot10^{-2}$ &$--$\\  
   			$ 2.08\cdot 10^{-4}$ &	$0.104$ & 2500&$5.71\cdot 10^{-3}$  &   $1.38$ &       $7.56\cdot 10^{-3}$ &$1.33$\\ 
   			$ 1.04\cdot 10^{-4}$ &	$0.052$ &9604 &$1.15\cdot 10^{-3}$  &   $2.31$ &       $1.60\cdot 10^{-3}$ &$2.34$\\ 
   			$ 5.20\cdot 10^{-5}$ &	$0.026$ &37636 &$2.49\cdot10^{-4}$ &   $2.21$ &       $3.33\cdot 10^{-4}$ &$2.27$\\  
   			\hline
   		\end{tabular}
   		\caption{Example 5: 2D smooth solution. Convergence of temperature at time  $t =0.0208$ from the ARS(2,2,1) scheme.
	%	\giovanni{I guess the oreder is relative to $h$v and $\Delta t$, right?}
	}
   		\label{2d_convergence_table2}
   	\end{center}
   \end{table}
%%%%%%%%%%%%%%%%%%%%%%%%%%%%%%%%%%%%%%%%%%%%%%%%%

\subsection{Example 6: Moving 2D shuttle with prescribed velocity} 

This example is an  extension of Example 5 to two space dimensions. We use a 2D velocity space.
We have taken this problem from the paper by Frangi et al. \cite{FFL},  where the authors have studied 
the biaxial accelerometer produced by STMicroelectronics with a surface micro-machining process.  The authors have analysed the problem by considering a two-dimensional simplification. 
In Figure \ref{shuttle} we have sketched the computational domain in details. The shuttle lies initially in the middle of the domain. In the rest of the domain a gas flow is 
taking place. The shuttle oscillates with the velocity $v = v_0\cos(2\pi \nu t)$, where $\nu$ is the frequency. We use SI units in the following. We set $v_0 = 1$.
The parameters mentioned in  Figure \ref{shuttle} are $L_1=19.2 \cdot 10^{-6}, ~ d_1 = 4.2\cdot 10^{-6}, ~ d_2 = 2.6\cdot 10^{-6}, ~ d_3 = 5\cdot 10^{-6}, ~ d_4 = 3.9\cdot 10^{-6},~d_5 = 18.8\cdot 10^{-6}$.  We have changed the parameter $\nu$  in  \cite{FFL} and have chosen $\nu = 40\cdot 4400$ Hz such that the maximum amplitude of the oscillations of the shuttle is half of the distance $d_2$ and the shuttle is  not touching the boundaries of the domain. 
The initial pressure of the gas is equal to $0.125$ bar, which corresponds to  an initial density $\rho_0 = 0.2$.   These parameters give a  relaxation time $\tau = 1.73\times 10^{-9}$ which is fixed for all times.

The initial distribution $f_0$ of the gas is the Maxwellian with zero mean velocity, initial temperature $T_0 = 293$ and initial density $\rho_0$. 
A diffuse reflection boundary condition with wall temperature $T_0$ is applied on the solid lines and a far field boundary condition $f_0$ is applied on the dotted lines. 
We note that in the present investigation  the time dependent motion of  the shuttle is resolved, while in  \cite{FFL} the authors solve  stationary 
equations with assigned non zero velocity on the boundary. 

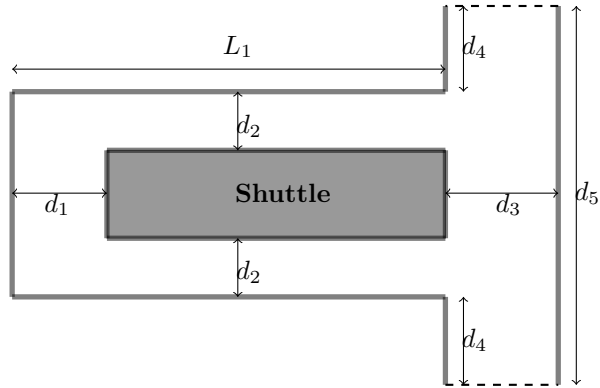
\begin{figure}[!t]
	\centering
		\begin{tikzpicture} [scale=0.3]
		\draw[line width=2pt,color=black,opacity=0.5] (4.2,6.5) -- ( 19.2,6.5);
		\draw[line width=2pt,color=black,opacity=0.5] (19.2,6.5) -- ( 19.2,10.4);
		\draw[line width=2pt,color=black,opacity=0.5] (19.2,10.4) -- ( 4.2,10.4);
		\draw[line width=2pt,color=black,opacity=0.5] (4.2,10.4) -- ( 4.2,6.5);
		\draw [fill=black,opacity=0.4] (4.2,6.5) rectangle (19.2,10.4);
		\node[color=black] at (12,8.5) {\bf Shuttle };
		\draw[line width=2pt,color=black,opacity=0.5] (0,3.9) -- ( 19.2,3.9);
		\draw[line width=2pt,color=black,opacity=0.5] (19.2,3.9) -- ( 19.2,0);
		\draw[dashed, thick] (19.2,0) -- ( 24.2,0);
		\draw[line width=2pt,color=black,opacity=0.5] (24.2,0) -- ( 24.2,16.8);
		\draw[dashed,thick]  (24.2,16.8) -- ( 19.2,16.8);
		\draw[line width=2pt,color=black,opacity=0.5] (19.2,16.8) -- ( 19.2,13);
		\draw[line width=2pt,color=black,opacity=0.5] (19.2,13) -- ( 0,13);
		\draw[line width=2pt,color=black,opacity=0.5] (0,13) -- ( 0,3.9);		     
		\draw[line width=0.4pt,<->] (0,8.5) -- (4.2,8.5);
		\node[color=black] at (2,8) {$d_1$ };
		\draw[line width=0.4pt,<->] (0,14) -- (19.2,14);
		\node[color=black] at (10,15) {$L_1$ };
		%  \draw[line width=0.4pt,<->] (4.2,6) -- (19.2,6);
		%\node[color=black] at (12.5,5.5) {$L_1$ };
		\draw[line width=0.4pt,<->] (10,13) -- (10,10.4);
		\node[color=black] at (10.5,11.5) {$d_2$ };
		\draw[line width=0.4pt,<->] (10,6.5) -- (10,3.9);
		\node[color=black] at (10.5,5) {$d_2$ };
		\draw[line width=0.4pt,<->] (19.2,8.5) -- (24.2,8.5);
		\node[color=black] at (22,8) {$d_3$ };
		\draw[line width=0.4pt,<->] (20,13) -- (20,16.8);
		\node[color=black] at (20.5,15) {$d_4$ };
		\draw[line width=0.4pt,<->] (20,3.9) -- (20,0);
		\node[color=black] at (20.5,2) {$d_4$ };
		\draw[line width=0.4pt,<->] (25,0) -- (25,16.8);
		\node[color=black] at (25.5,8.5) {$d_5$ };
		%   \node[color=black] at (22,-0.5) {Far field  };
		%    \node[color=black] at (22,17.3) {Far field  };
		\end{tikzpicture}
		\vspace{0.3cm}
		\caption{Example 6: 2D shuttle. Geometry setup for moving $2D$ shuttle. }
		\label{shuttle}
	\centering
\end{figure}

%\giovanni{Congratulations for the figure! It must have taken a long time to draw!}

In Figure \ref{shuttle_velocity} we have plotted the velocity vector fields as well as $x$- and $y$- components of the velocity at times $t = 1.2\cdot 10^{-6}$. 
Notice that the period of oscillations here is $T\approx 5.61\cdot 10^{-6}$. The total number of grid points is approximately $7000$ which gives   $h = 2.5 \Delta x= 4\cdot10^{-7}$. The first order Euler scheme is used for 
the time integration with the time step $\Delta t = 1.5\cdot10^{-10}$. 

\begin{figure}[!t]
\centering
\includegraphics[keepaspectratio=true, width=0.45\textwidth]{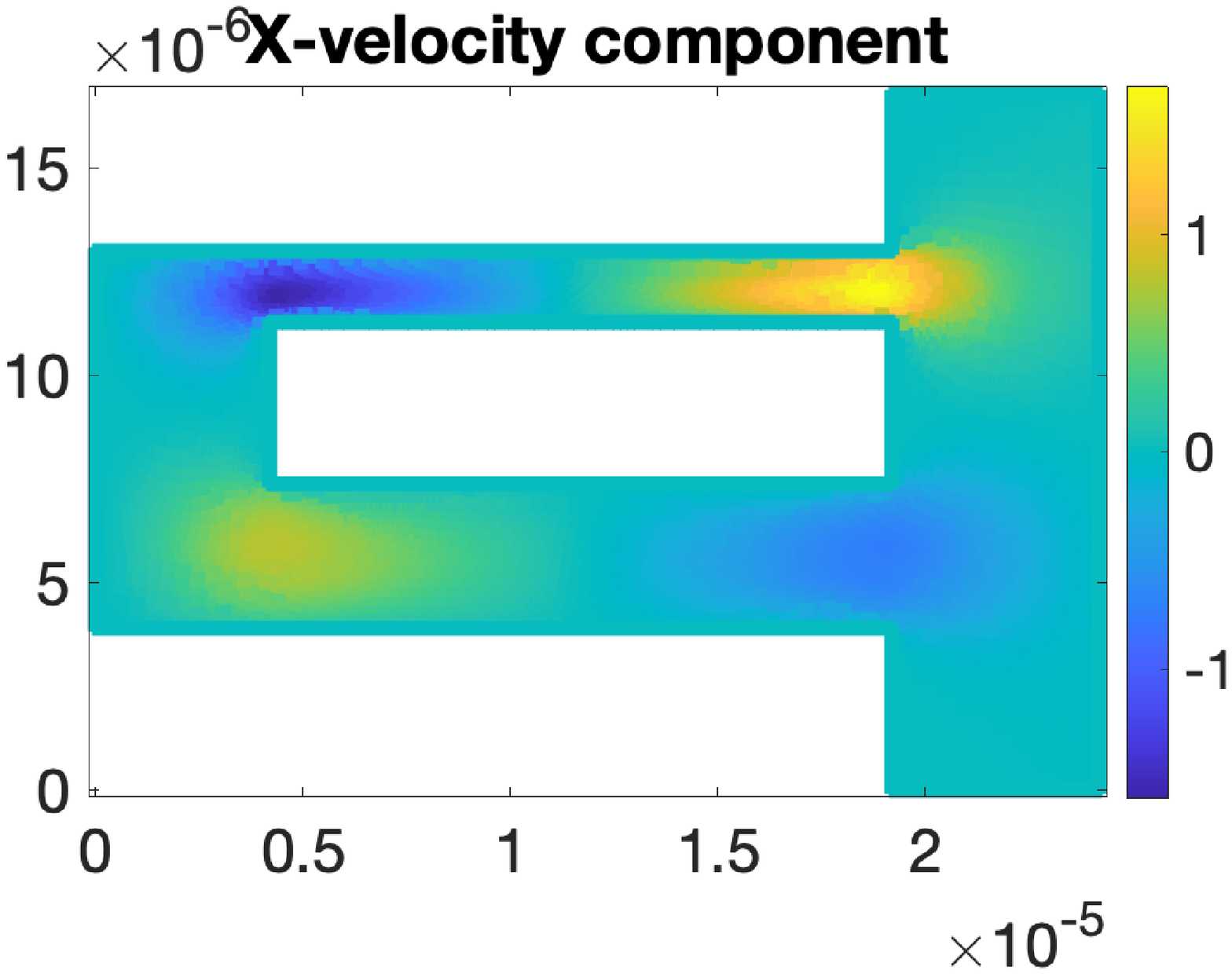}
\includegraphics[keepaspectratio=true, width=0.45\textwidth]{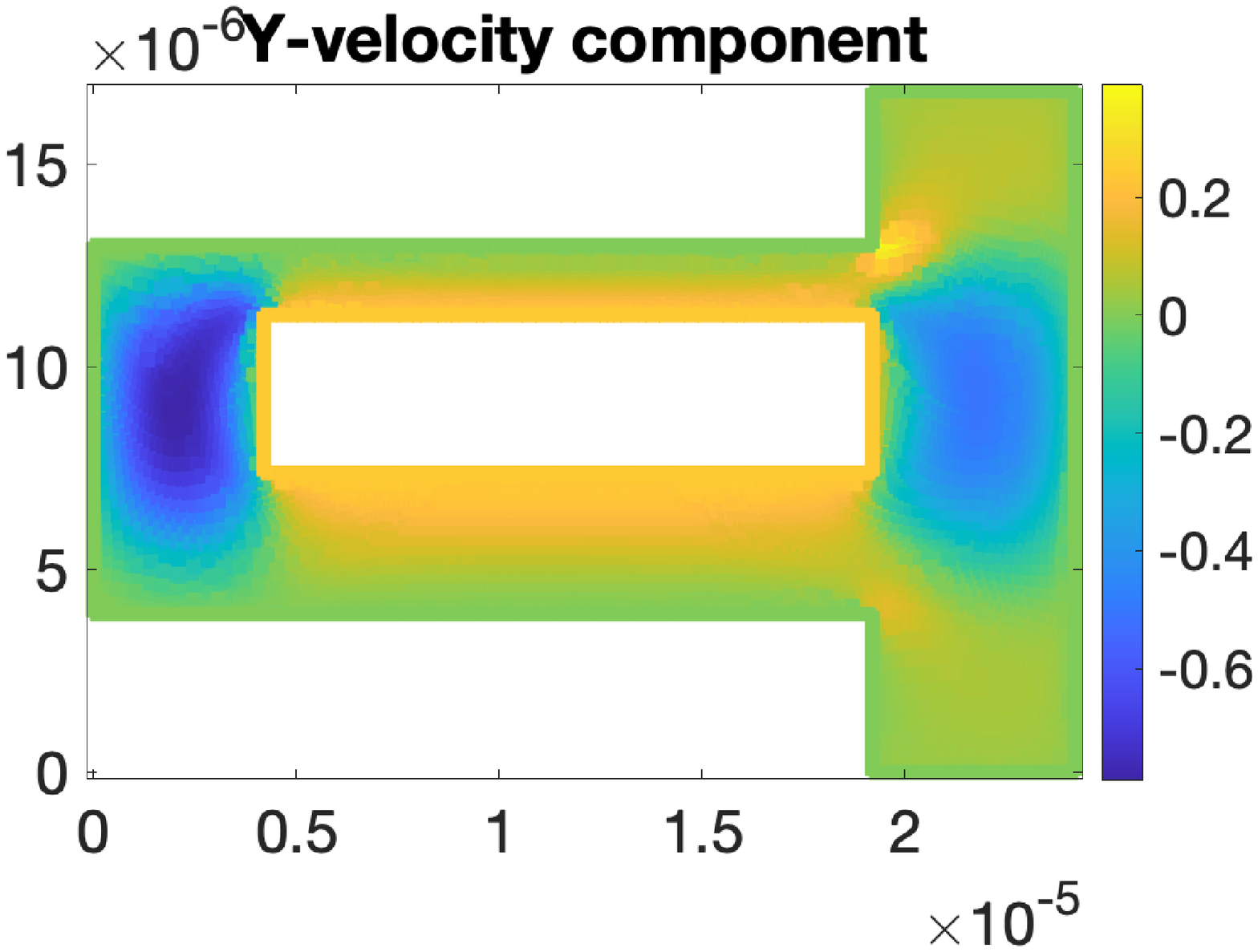}\\
\includegraphics[keepaspectratio=true, width=0.45\textwidth]{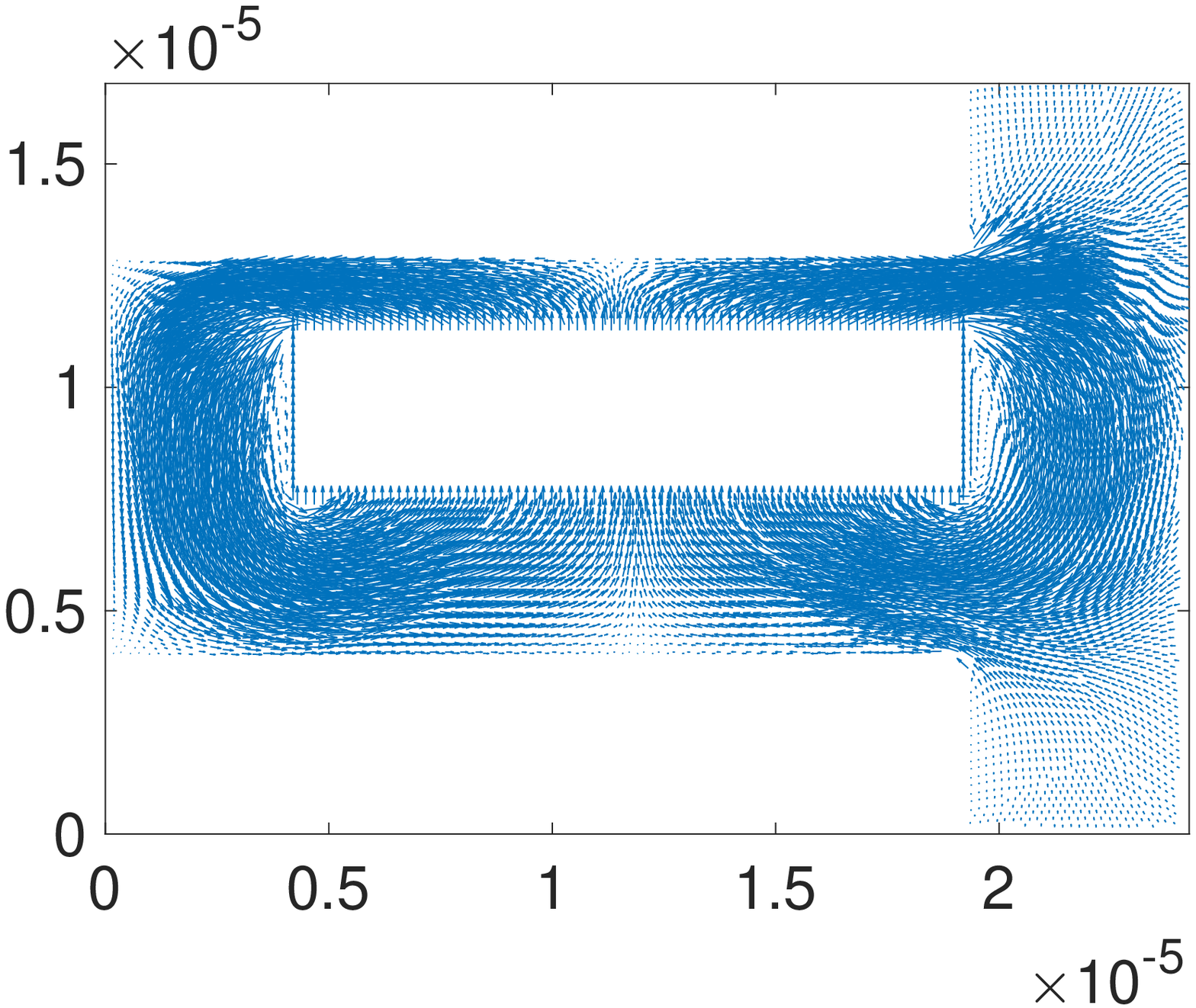} 
\caption{ Example 6: 2D shuttle.  First  row:  $x$- and $y$- velocity components at time $t = 1.2 \cdot 10^{-6}$. Second row: velocity fields at time $t=1.2\cdot10^{-6}$.  }
\label{shuttle_velocity}
\centering
\end{figure}

 \subsubsection{Convergence study}
 
 In Figure \ref{stress_on_shuttle} we have plotted the normal stress tensor on the top wall of the shuttle at time $t=1.2\cdot 10^{-6}$. As a reference solution we consider the one obtained at the  finest resolution with  $h = 1.2\cdot 10^{-7}$, which corresponds to approximately $111000$ grid points including boundary points. The finest time step is chosen as $8\cdot 10^{-11}$. 
 The results of the  convergence study are  presented in Table \ref{table_shuttle}. 
 
 Table \ref{table_shuttle} shows the results for the first order scheme in time and space. In order to estimate the error, we have generated a fixed number  $N=100$  of points in equal distance at the upper boundary of the shuttle.

 On these points we have interpolated the  stress tensors from different resolutions including the reference solutions and then computed the errors. In Table \ref{table_shuttle}  the $L^1$ and $L^2$ errors of the normal stress tensor $\varphi_{yy}$ are presented. The errors in the table show the first order convergence of the scheme. 
 
 Table \ref{table_shuttle2} shows the results for the ARS(2,2,1) scheme.
 We observe  an  improvement compared to the first order scheme, but we  obtain in this situation a rate of convergence still below $2$, which is expected due to the non-smooth geometry.
 
 \begin{figure}[!t]
 	\begin{center}
 	\includegraphics[keepaspectratio=true, width=0.45\textwidth]{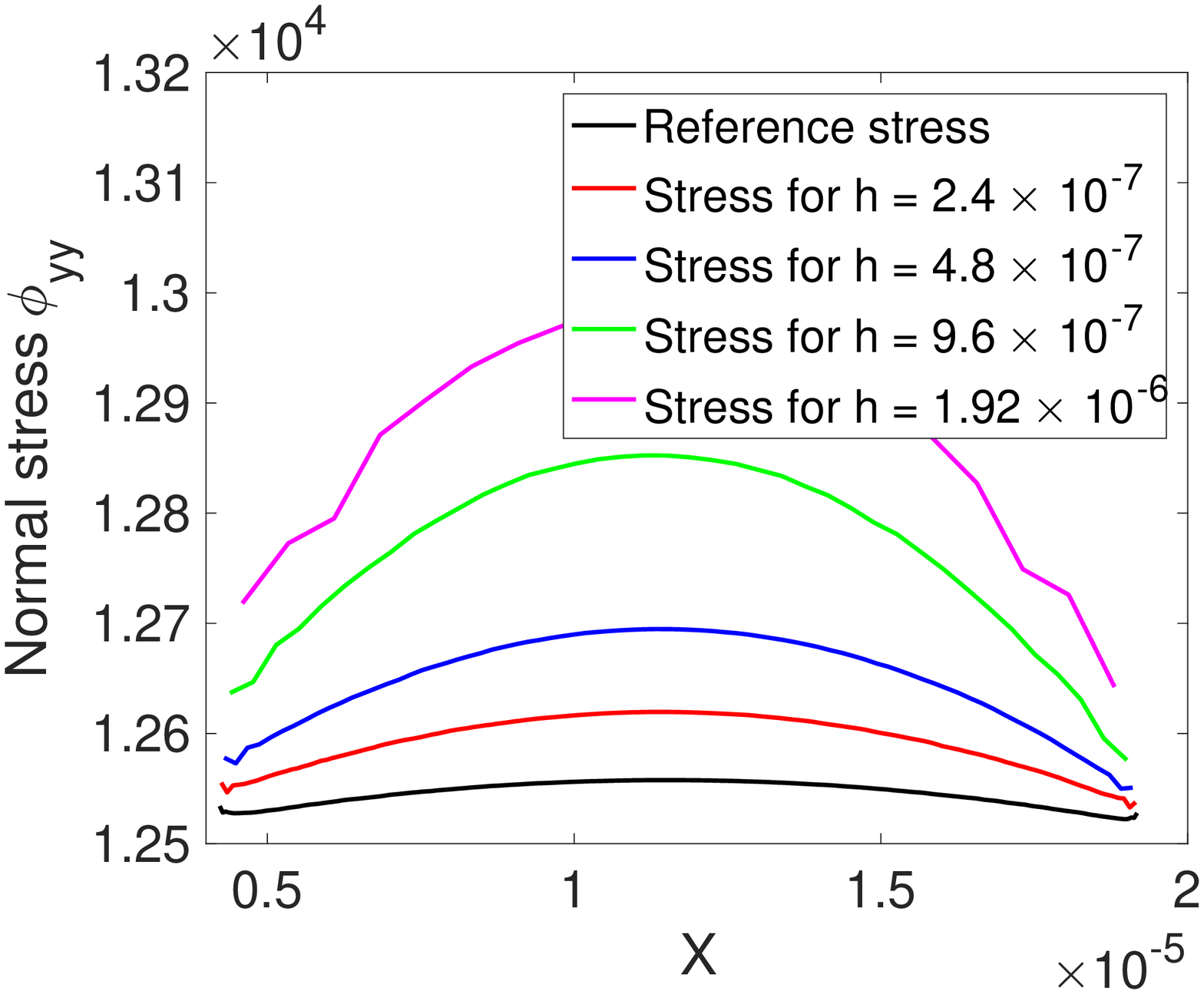}
	\includegraphics[keepaspectratio=true, width=0.45\textwidth]{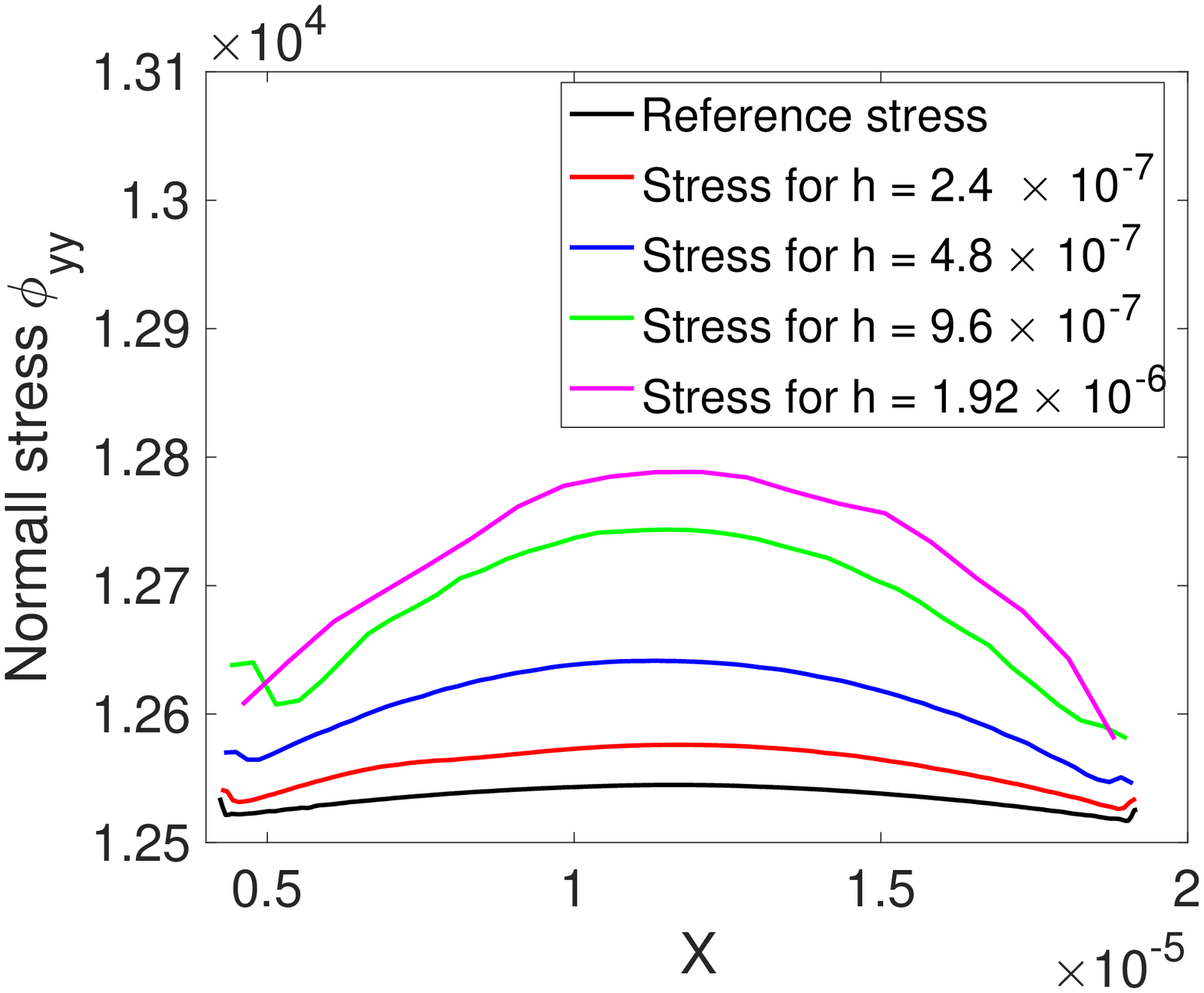}
 		\caption{Example 6: 2D shuttle. The normal stress tensor on the top wall of the shuttle at $t = 1.2\cdot 10^{-6}$ for different cell sizes. Left: First order in space and time. Right: Second order in space and time.  }
 		\label{stress_on_shuttle}  
 	\end{center}
 \end{figure}

 \begin{table}
 	\begin{center}
 		\begin{tabular}{|r|r|r|r|r|r|}
 			\hline
 		$\Delta t$ &	$h $  & $L^1$-error & Order \\ \hline 
 		$ 64\cdot 10^{-11}$ &	$1.92\cdot10^{-6}$ & $4.85\cdot10^{-3}$   &    $--$           \\  
 		$ 32\cdot 10^{-11}$ &	$9.6\cdot10^{-7}$ & $3.27\cdot 10^{-3}$  &   $0.57$ \\ 
 	    $ 16\cdot 10^{-11}$ &	$4.8\cdot10^{-7}$ & $1.54\cdot 10^{-3}$  &   $1.08$ \\ 
 		$ 8\cdot 10^{-11}$ &	$2.4\cdot10^{-7}$ & $7.04\cdot10^{-4}$ &   $1.13$ \\  
 			\hline
 		\end{tabular}
 		\caption{Example 6: 2D shuttle. Convergence of the normal stress tensor $\phi$ on the top wall of the shuttle at time  $t =1.2\cdot10{-6}$ from the first  order scheme in space and time.}
 		\label{table_shuttle}
 	\end{center}
 \end{table}

  \begin{table}
  	\begin{center}
  		\begin{tabular}{|r|r|r|r|r|r|}
  			\hline
  			$\Delta t$ &	$h $  & $L^1$-error & Order \\ \hline 
  			$ 96\cdot 10^{-11}$ &	$1.92\cdot10^{-6}$ & $2.73\cdot10^{-3}$   &    $--$ \\  
  			$ 48\cdot 10^{-11}$ &	$9.6\cdot10^{-7}$ & $2.22\cdot 10^{-3}$  &   $0.30$ \\ 
  			$ 24\cdot 10^{-11}$ &	$4.8\cdot10^{-7}$ & $1.09\cdot 10^{-3}$  &   $1.03$ \\ 
  			$ 12\cdot 10^{-11}$ &	$2.4\cdot10^{-7}$ & $3.55\cdot10^{-4}$ &   $1.62$ \\
  			\hline
  		\end{tabular}
  		\caption{Example 6: 2D shuttle. Convergence of the normal stress tensor $\phi$ on the top wall of the shuttle at time  $t =1.2\cdot10^{-6}$ from the ARS(2,2,1) scheme.
		%\giovanni{Something si not right: the velocities should be of order 1 m/s, not 
		%$10^{-6}$ m/s}
		}
  		\label{table_shuttle2}
  	\end{center}
  \end{table}

 %%%%%%%%%%%%%%%%%%%%%%%%%%%%%%%%
 
 %%%%%%%%%%%%%%%%%%%%%%%%%%%%%%%%
\subsection{Example 7: Transport of rigid particles}

The main aim of the following  tests is to demonstrate the ability  of the scheme to simulate arbitrary shapes of rigid body motion immersed in a rarefied
gas. We consider again two dimensional physical and velocity space.   In the previous 2-D test case a one-way coupling of rigid body motion and gas was investigated. In the present  example we consider a two-way coupling,
where the gas is also influencing the  motion of the rigid body. Using SI units, we consider the 
computational domain $\Omega = [0,2\cdot10^{-6}]\times[0,3\cdot10^{-6}]$. The initial density is $\rho_0 = 1$, the initial temperature $T_0 = 270$ and the initial mean velocity $U_0 = (0,0)$. 
These parameters yield the initial relaxation time $\tau = 3.71 \times 10^{-10}$ which is fixed for all times.
On the top we prescribe a  Maxwellian with  parameters $\rho=\rho_0, T=290, U=U_0$. On the  bottom  boundary 
%prescribe a   Maxwellian distribution with wall temperature $T_w = 290$, density $\rho_0$ and mean velocity $U_0$. On the bottom 
we use a  diffuse reflection boundary condition with wall temperature $T_w = T_0, U_w=U_0$. 

% \giovanni{Note that the motion is driven only by the temperature gradient, and we do not take into account the effect of gravity.}
% \axel{
% This is not really true, since by fixing the density of the Maxweillian we actually generate a pressure difference which is also driving the motion of the gas. Otherwise its a very slow motion}
On the left and right wall we apply  far field boundary conditions, that means, we prescribe a  Maxwellian with initial parameters $\rho_0, T_0, U_0$. On the rigid body we apply a diffuse reflection boundary condition with temperature $T(t,x)$ and velocity $U(t,x)$. We consider  circular as well as chiral particles. For the following simulations we use  the first order scheme in space and time. 

\subsubsection{Transportation of a circular particle}
First we consider a circular particle of radius $0.1\cdot10^{-6}$ and initial center of mass $(1.0\cdot10^{-6}, 2.5\cdot10^{-6})$. The grid points are generated equidistantly with $h=2.5 \Delta x=5.25\cdot10^{-8}$ which gives an  initial number of grid points  equal to $7273$.
% as plotted in Fig \ref{thermo_circle0} . 
 The time step is $\Delta t = 1\cdot10^{-11}$. 
            
In Figures \ref{thermo_circle_velo} and \ref{thermo_circle_temp} we have plotted the positions of the circular particle together with velocity fields and temperature fields, respectively, at  times $  1\cdot10^{-7}, 3\cdot10^{-7} $ and  $4.5\cdot10^{-7} $.

\begin{figure}[!t]
	\centering
   \includegraphics[keepaspectratio=true, angle=0, width=0.35\textwidth]{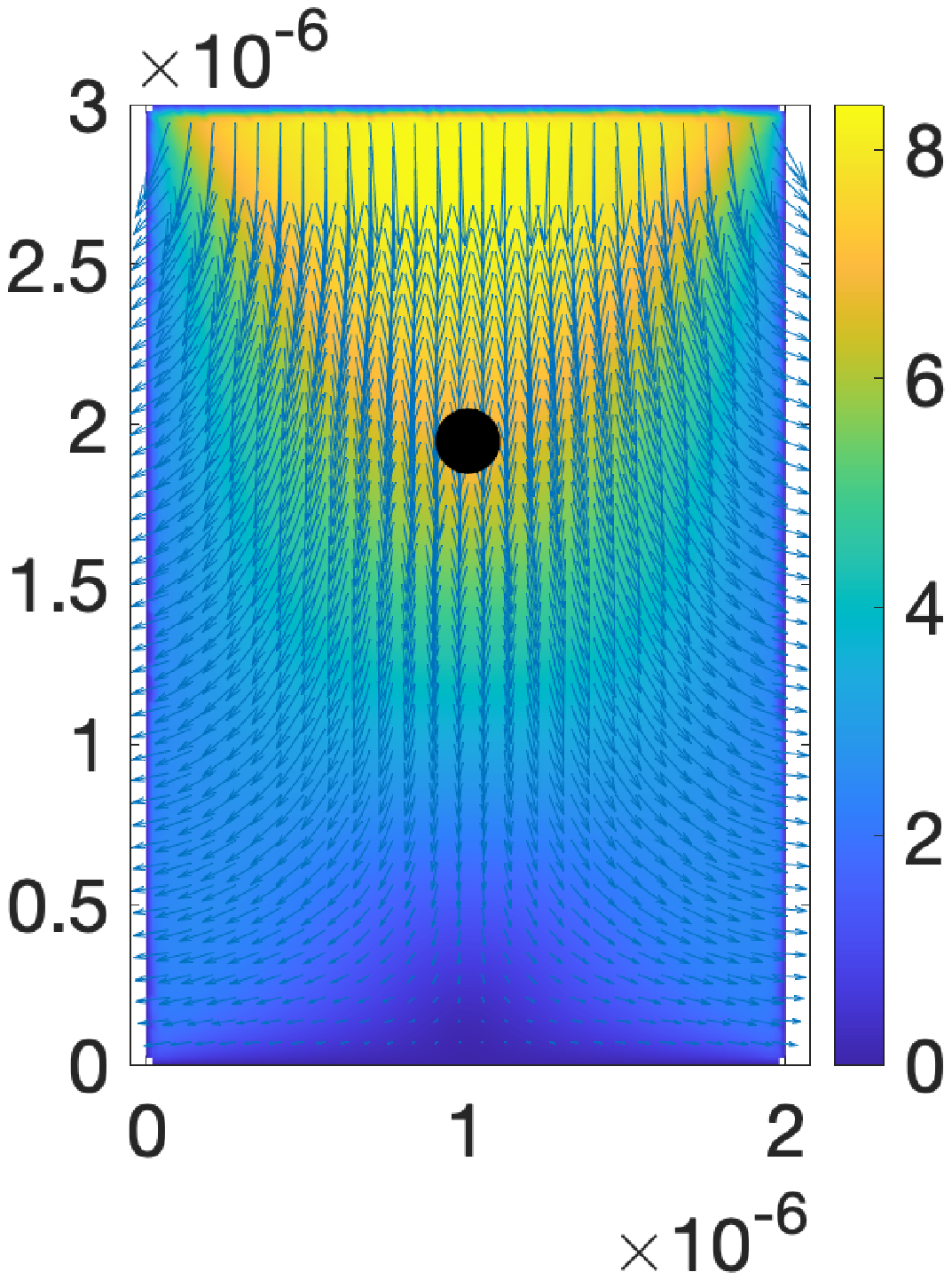}
   \hspace{-0.8cm}
     \includegraphics[keepaspectratio=true, angle=0, width=0.35\textwidth]{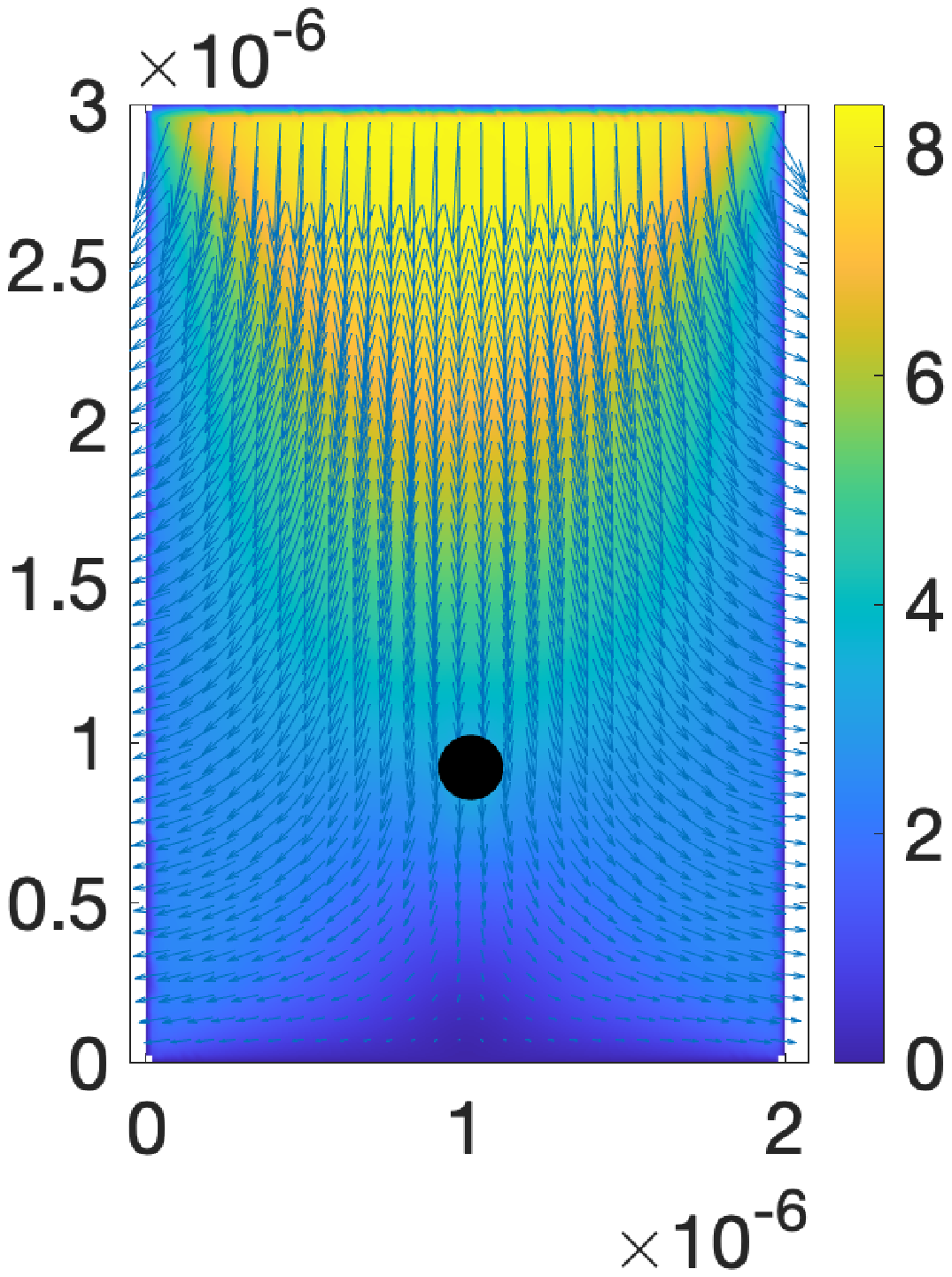}
    \hspace{-0.8cm} 
      \includegraphics[keepaspectratio=true, angle=0, width=0.35\textwidth]{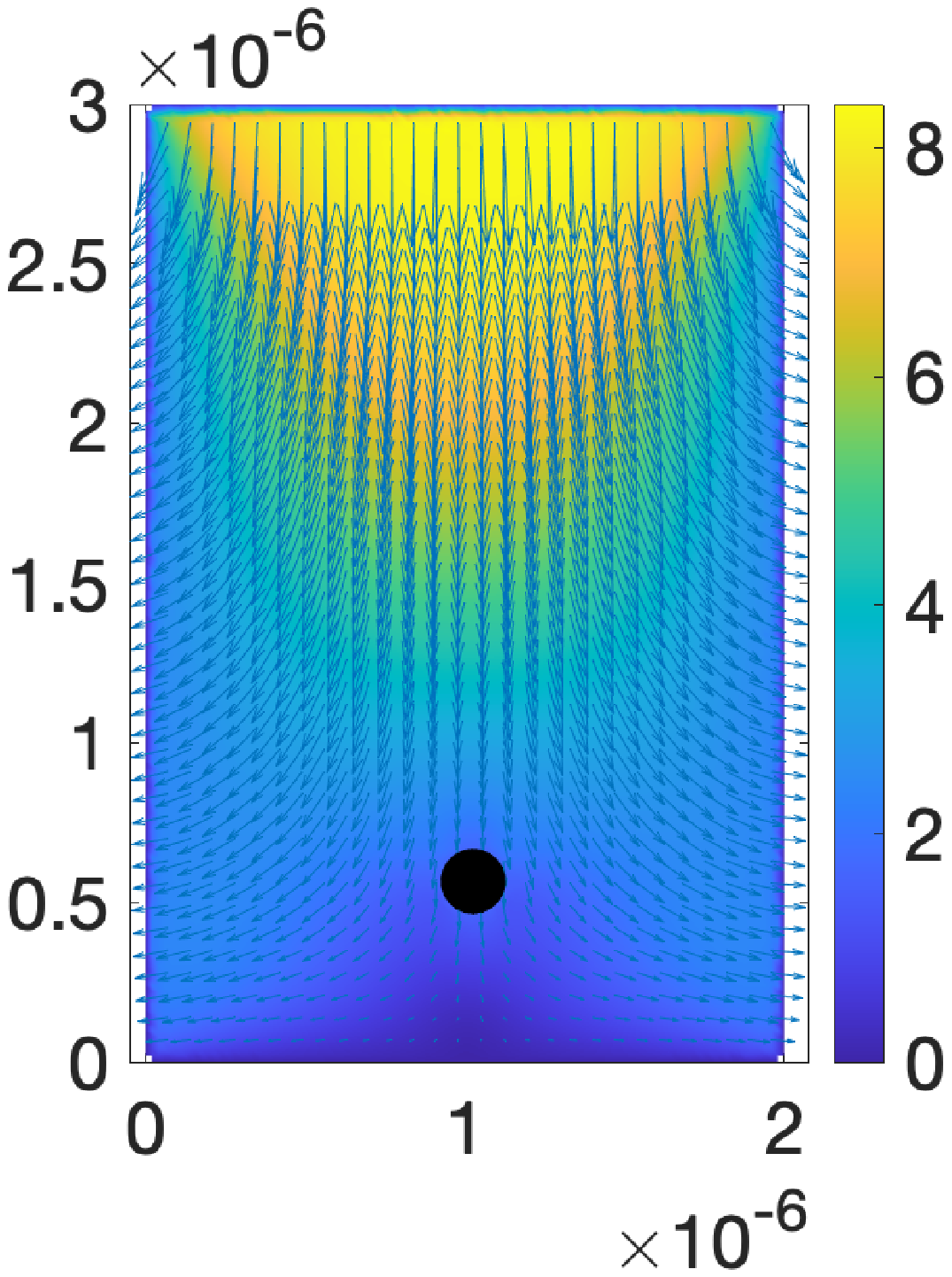}
% \vspace{-0.5cm}
	\caption{Example 7: Circular particles.  Particle positions and velocity field at %$ t =2\cdot10^{-7}$ (Left),
		$ t =1\cdot10^{-7}$,
		$ t =3\cdot10^{-7}$ and  $ t = 4.5\cdot10^{-7}$. }
	\label{thermo_circle_velo}
	\centering
\end{figure}   

%%%%%%%%%%%%%%%%%%
\begin{figure}[!t]
	\centering
     \includegraphics[keepaspectratio=true, angle=0, width=0.35\textwidth]{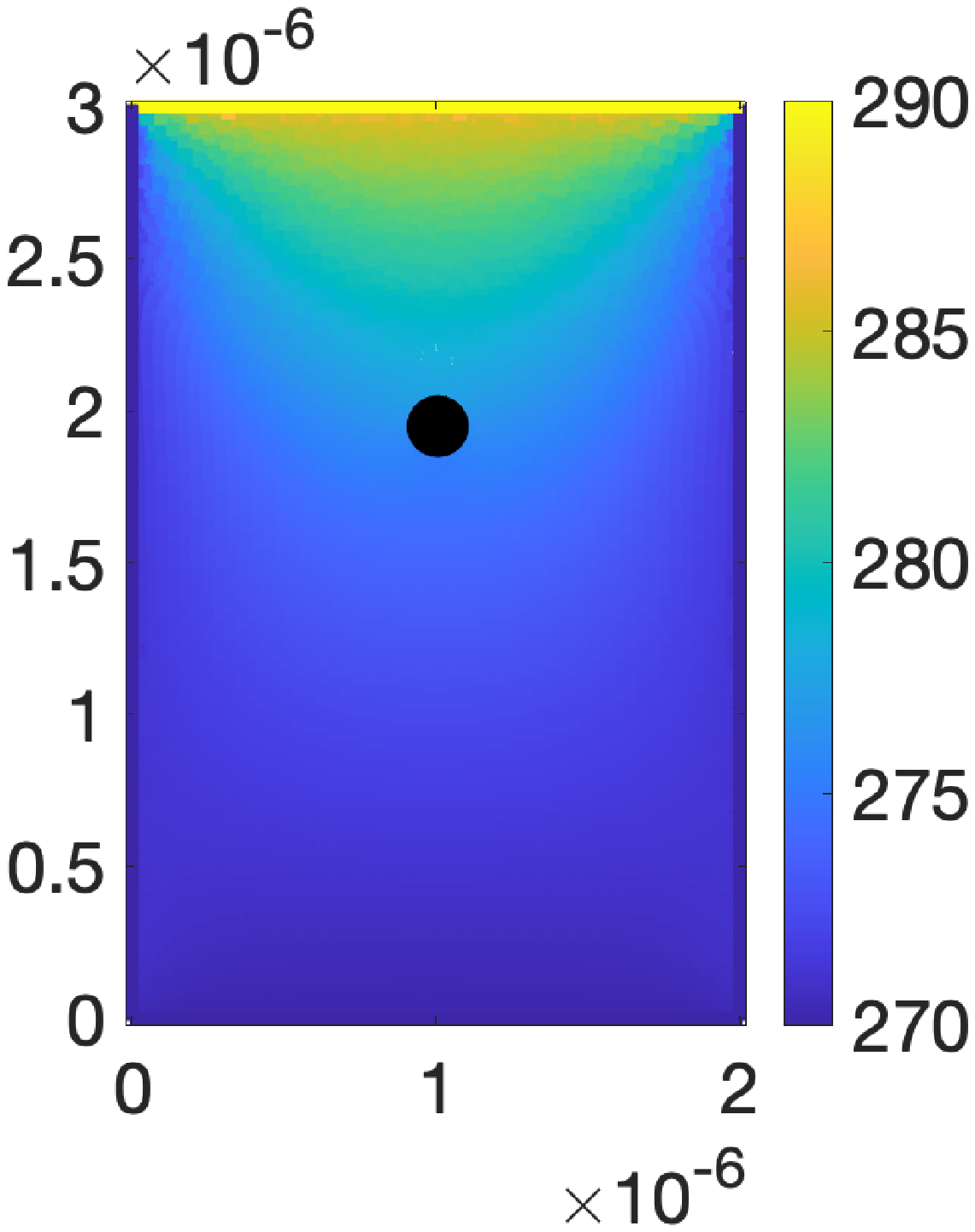}\hspace{-0.8cm} 
     \includegraphics[keepaspectratio=true, angle=0, width=0.35\textwidth]{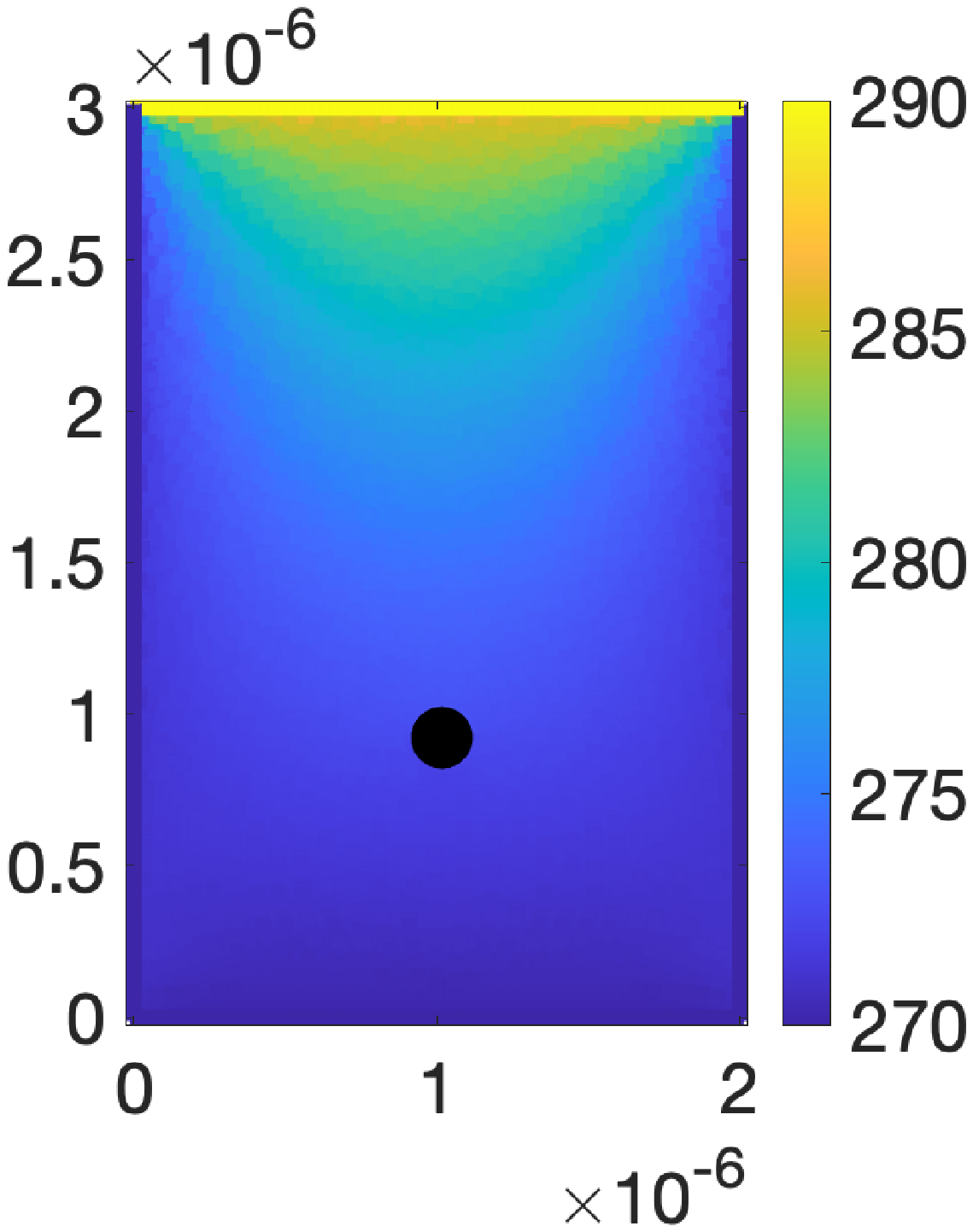}\hspace{-0.8cm} 
   	 \includegraphics[keepaspectratio=true, angle=0, width=0.35\textwidth]{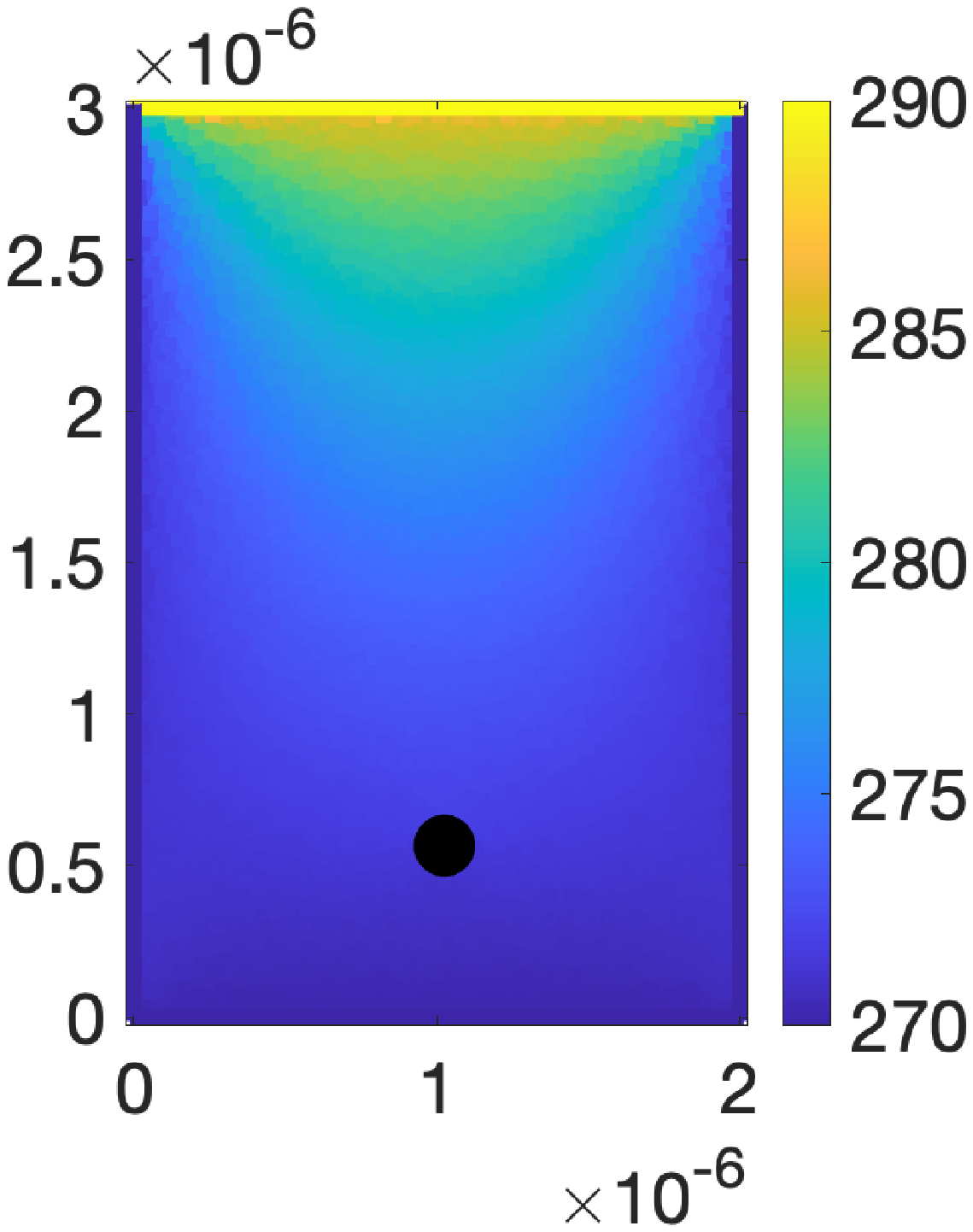}
 %  	  \vspace{-0.5cm}
	\caption{ Example 7: Circular particles.  Particle positions and temperature field at
		%$ t =2\cdot10^{-7}$ (Left),
		$ t =1\cdot10^{-7}$, $ t =3\cdot10^{-7}$  and $ t = 4.5\cdot10^{-7}$.  }
	\label{thermo_circle_temp}
	\centering
\end{figure}           
                  
%%%%%%%%%%%%%%%%%
\subsubsection{Transportation of a chiral particle}    
 
In this example, we consider a chiral particle with  initial center of mass $(1.0\cdot10^{-6}, 2.3\cdot10^{-6})$.
We have used a relatively fine grid with $h=6.25\cdot10^{-8}$, which gives   $9528$ particles and a time step $\Delta t = 1\cdot10^{-11}$. The boundary conditions are the same as in the case of the circular particle in the previous subsection.
                  
In Figures \ref{thermo_chiral_velo} and \ref{thermo_chiral_temp} we have plotted the positions of the chiral particle together with velocity fields and temperature fields, respectively, at  times $ 1\cdot10^{-7}, 3\cdot10^{-7} $ and  $4.5\cdot10^{-7} $. 
                  
 \begin{figure}[!t]
   \centering
    \includegraphics[keepaspectratio=true, angle=0, width=0.35\textwidth]{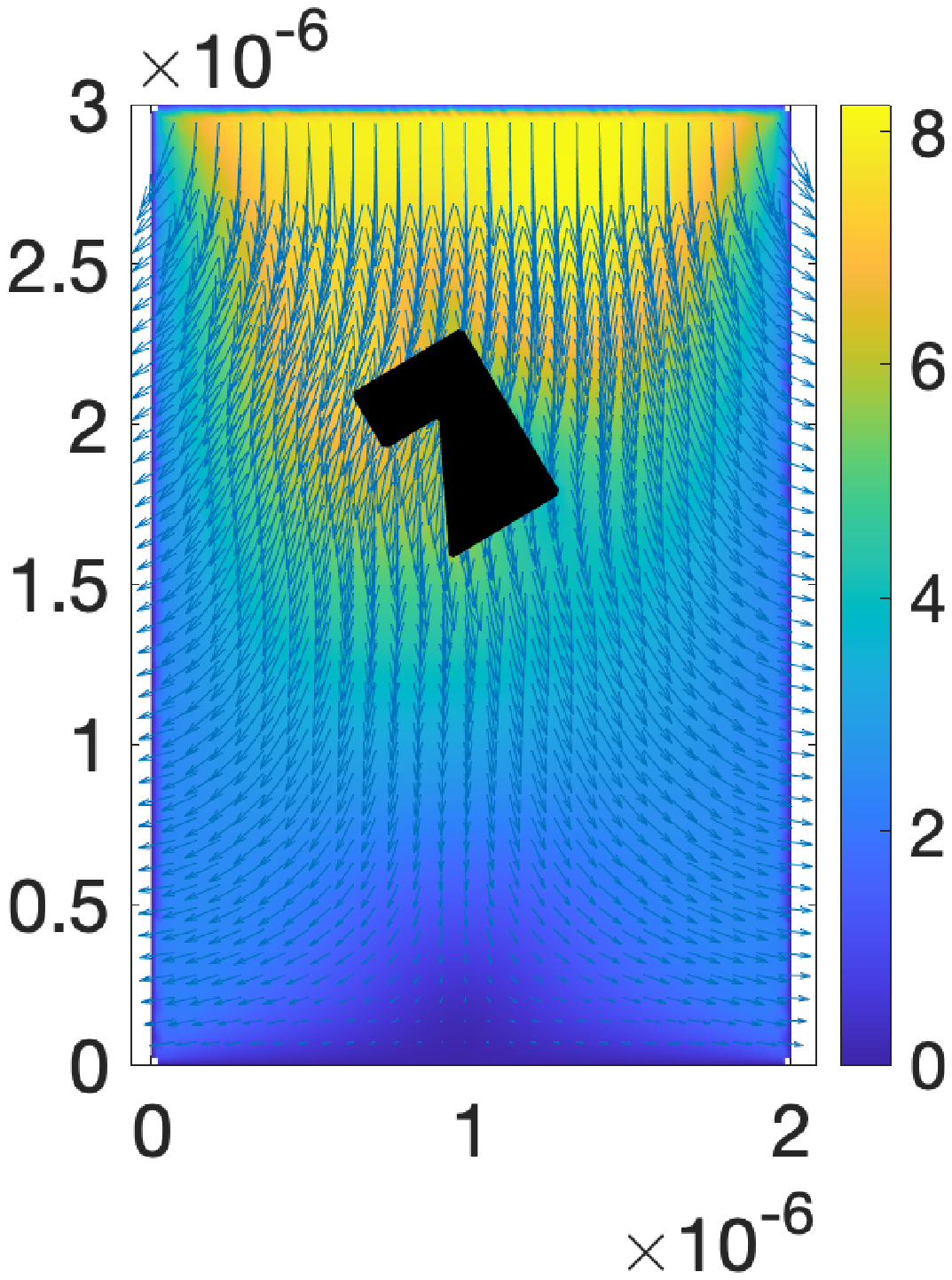}\hspace{-0.8cm} 
     \includegraphics[keepaspectratio=true, angle=0, width=0.35\textwidth]{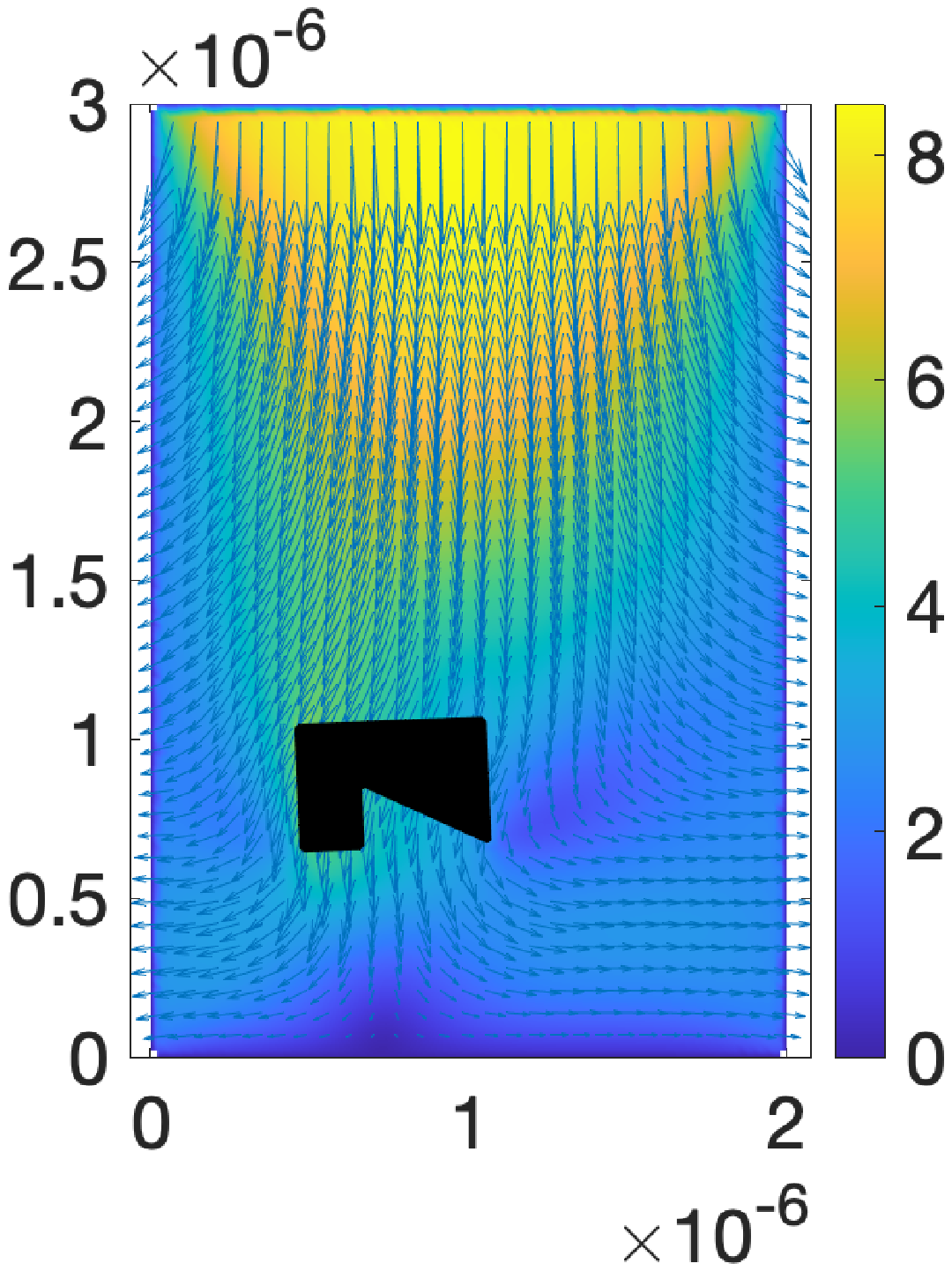}\hspace{-0.8cm} 
      \includegraphics[keepaspectratio=true, angle=0, width=0.35\textwidth]{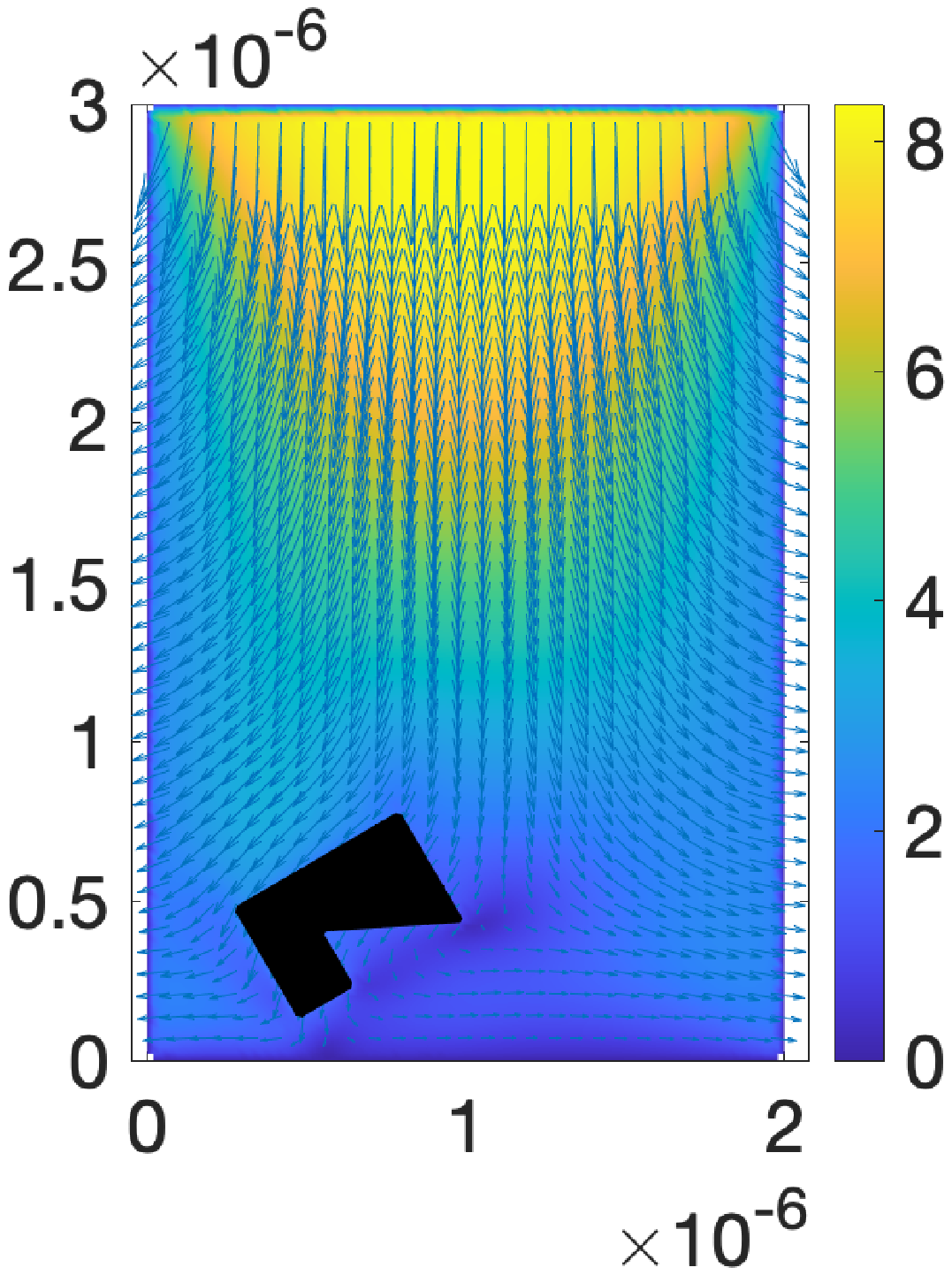}
   %    \vspace{-0.5cm}
      \caption{Example 7: Chiral particles.  Particle positions and velocity field at 
      	%$ t =2\cdot10^{-7}$ (Left),
      	$ t =1\cdot10^{-7}$, $ t =3\cdot10^{-7}$ and  $ t = 4.5\cdot10^{-7}$ }
   	\label{thermo_chiral_velo}
   \centering
  \end{figure}

 %%%%%%%%%%%%%%%%%%%
\begin{figure}[!t]
 \centering
      \includegraphics[keepaspectratio=true, angle=0, width=0.35\textwidth]{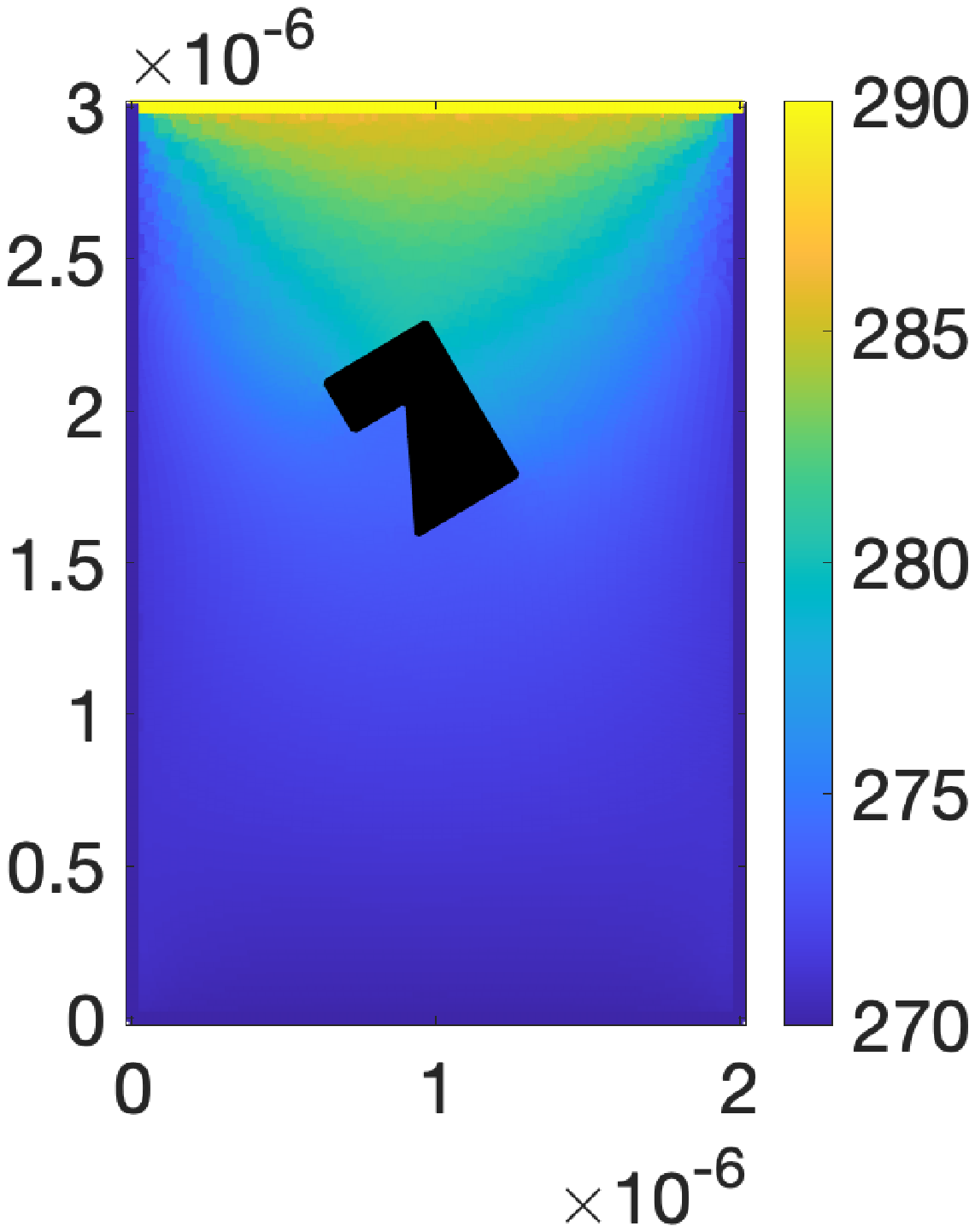}\hspace{-0.6cm} 
     \includegraphics[keepaspectratio=true, angle=0, width=0.35\textwidth]{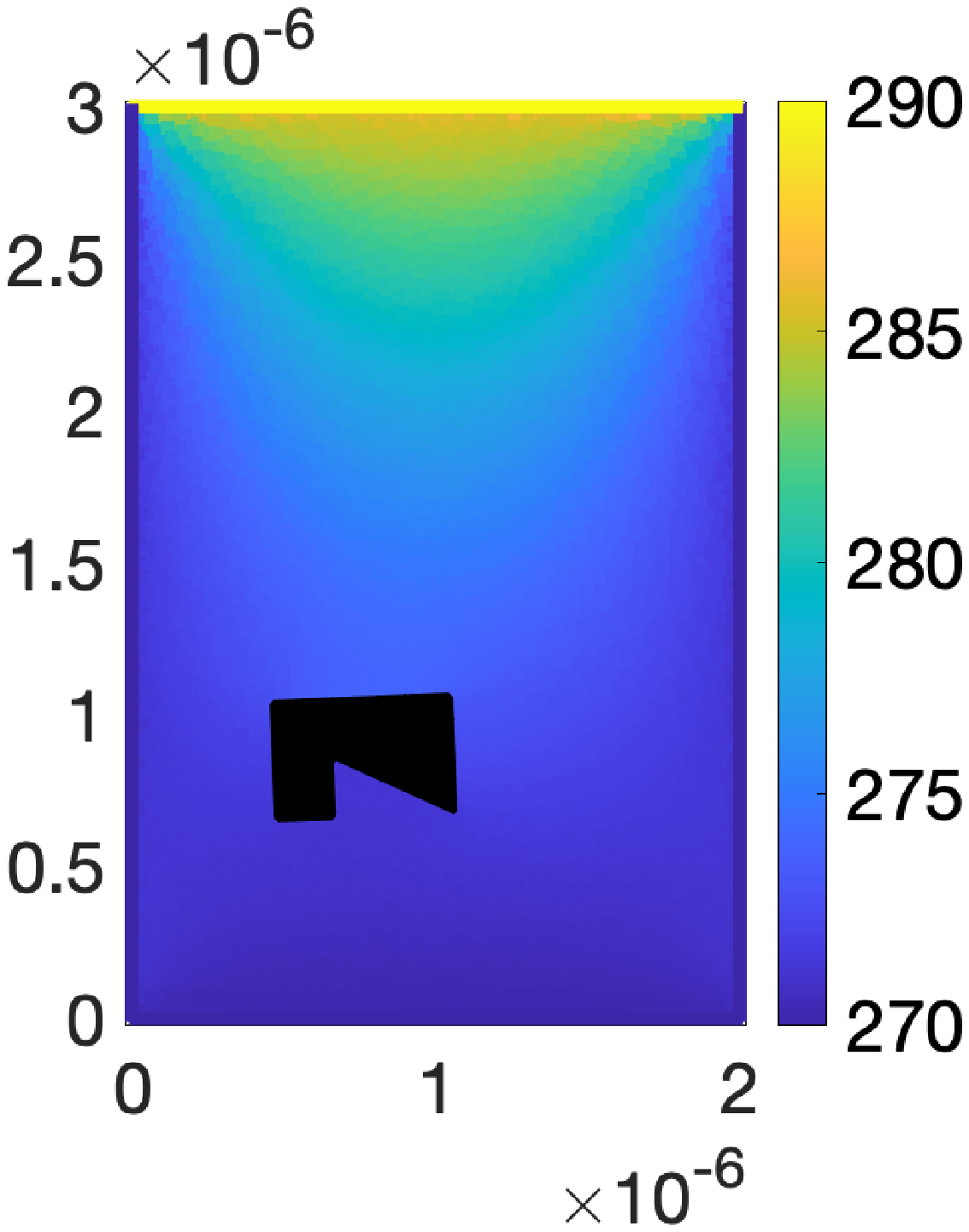}\hspace{-0.9cm} 
   	 \includegraphics[keepaspectratio=true, angle=0, width=0.35\textwidth]{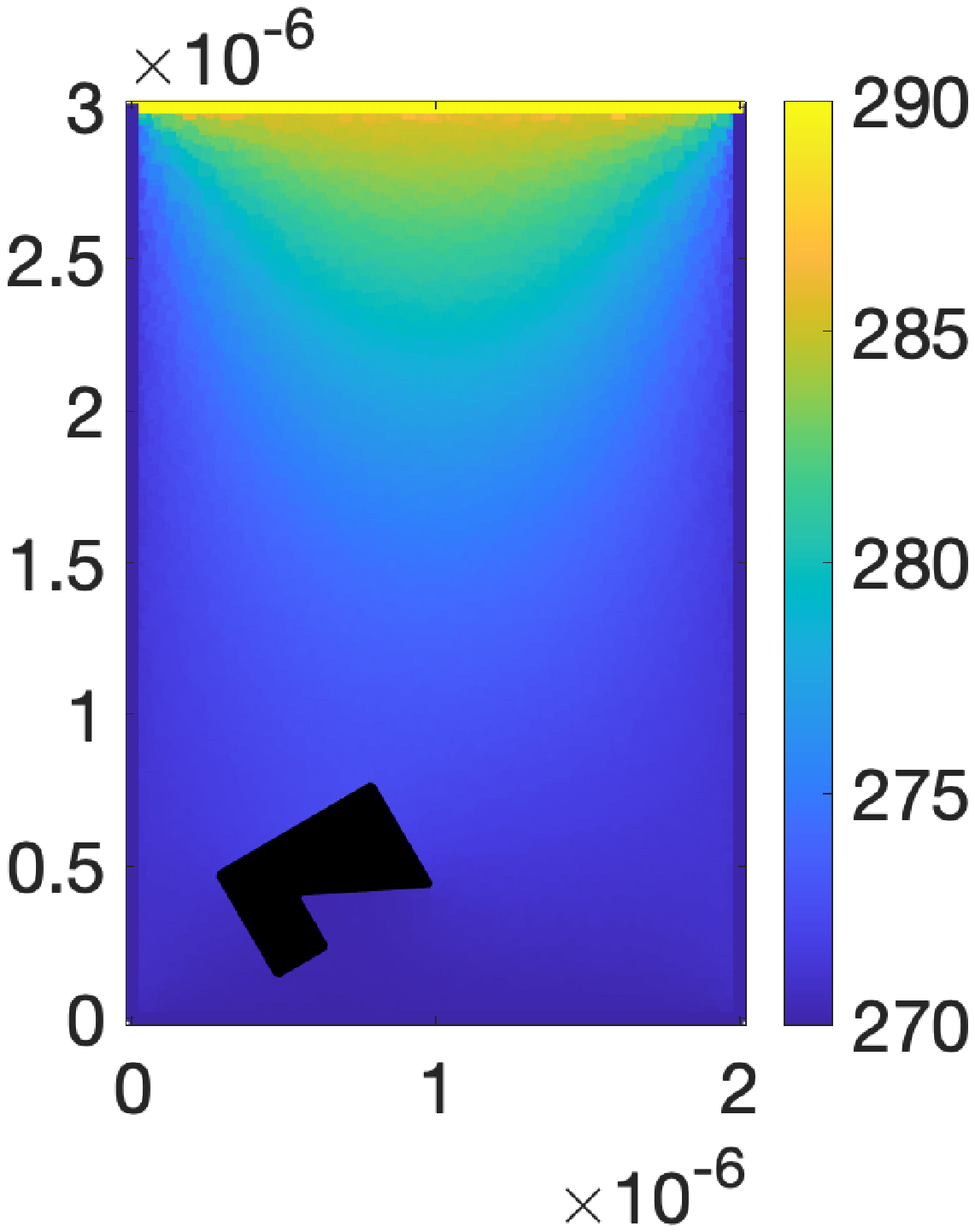}
%   	 \vspace{-0.5cm}
      \caption{Example 7: Chiral particles.   Particle positions and temperature field at 
      	%$ t =2\cdot10^{-7}$ (Left),
      	$ t = 1\cdot10^{-7}$, $ t =3\cdot10^{-7}$ and  $ t = 4.5\cdot10^{-7}$.  }
       \label{thermo_chiral_temp}
   	\centering
 \end{figure}           
                  
%%%%%%%%%%%%%%%%%%%%%%%%%%

\subsection{Multiple rigid particles in a driven cavity}

We consider a square cavity  $[0, L] \times [0, L] $ with  $L=1\cdot10^{-6}$.
The initial parameters of the Maxwellian are the same as in the previous test case. 
Diffuse reflection boundary conditions with  temperature  $T_0$ are applied on all boundaries as well as on the rigid particles. At the top wall we prescribe a  non-zero velocity in  $x$-direction given by 
\[
U_0^{(x)} =  10 \left(\frac{2x}{L}\right)^2 \left( 2-   \left(\frac{2x}{L}\right)^2\right).
\] 
This leads to a  maximum velocity equal to $10$ at  the center of the wall. The $y$-component of the top wall velocity is zero.  The velocities on all other walls and on the rigid particles are zero. We have generated 4 rigid particles of radius $0.075L$ with initial position as  in Figure \ref{circle_in_cavity}, first panel. The numerical particles are generated according to the parameter $h = 5.25\cdot10^{-8}$ which gives, initially,  a total number of $2313$ particles. The time step is chosen as $\Delta t = 1\cdot10^{-11}$. 

\begin{figure}[!t]
	\centering
	\includegraphics[keepaspectratio=true,  width=0.42\textwidth]{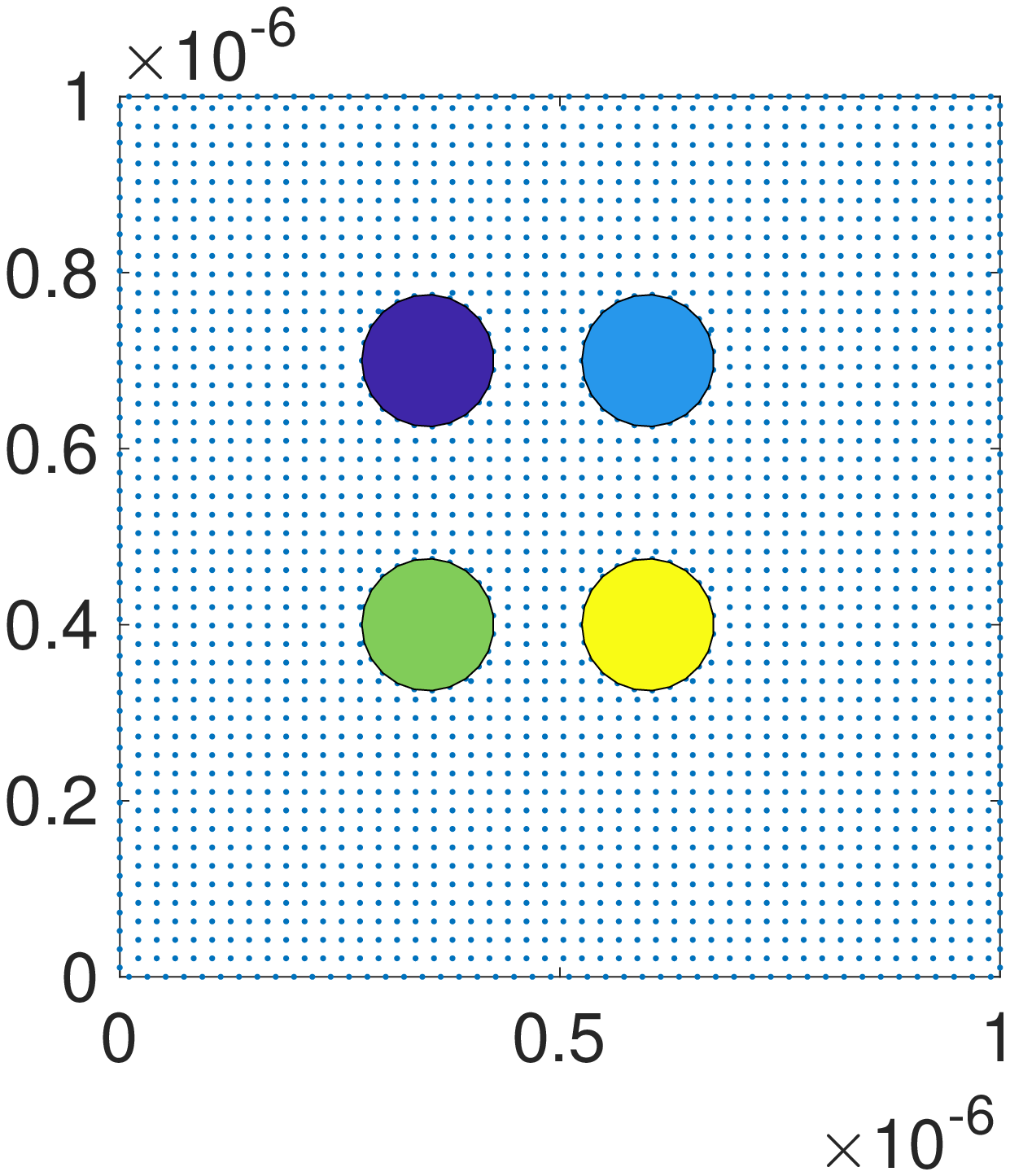}
	\includegraphics[keepaspectratio=true, width=0.42\textwidth]{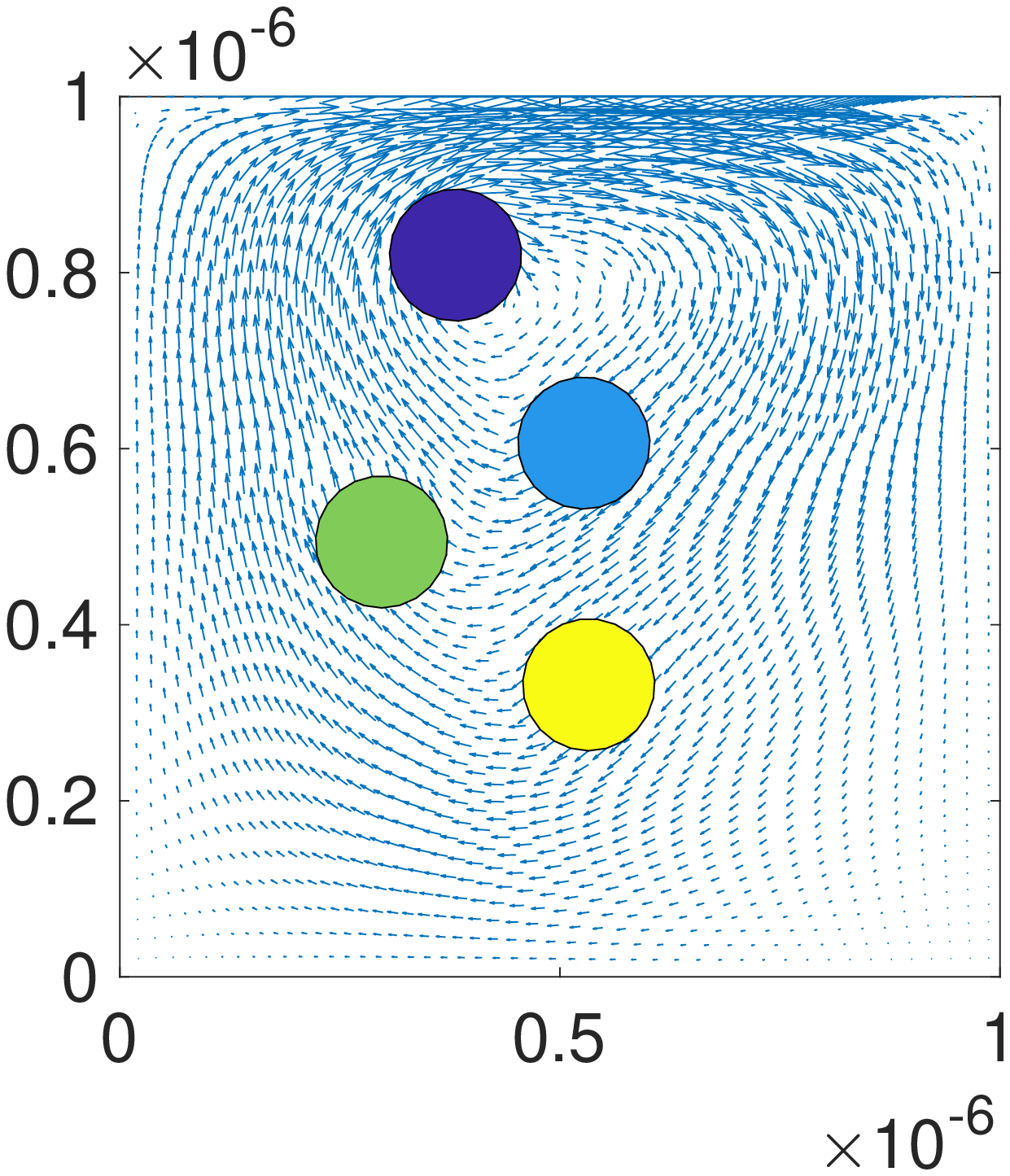}\\
	\includegraphics[keepaspectratio=true,  width=0.42\textwidth]{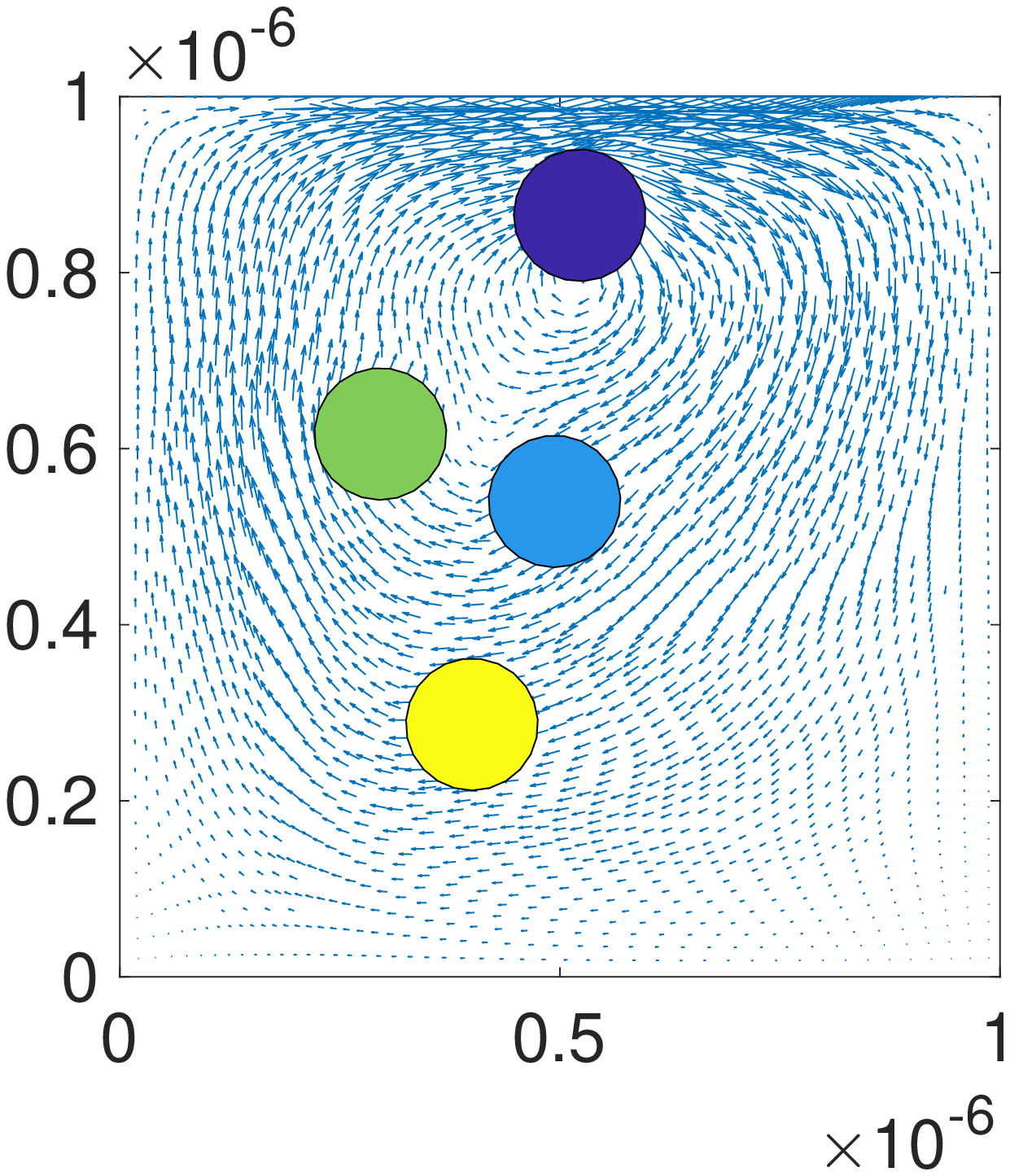}
	\includegraphics[keepaspectratio=true,  width=0.42\textwidth]{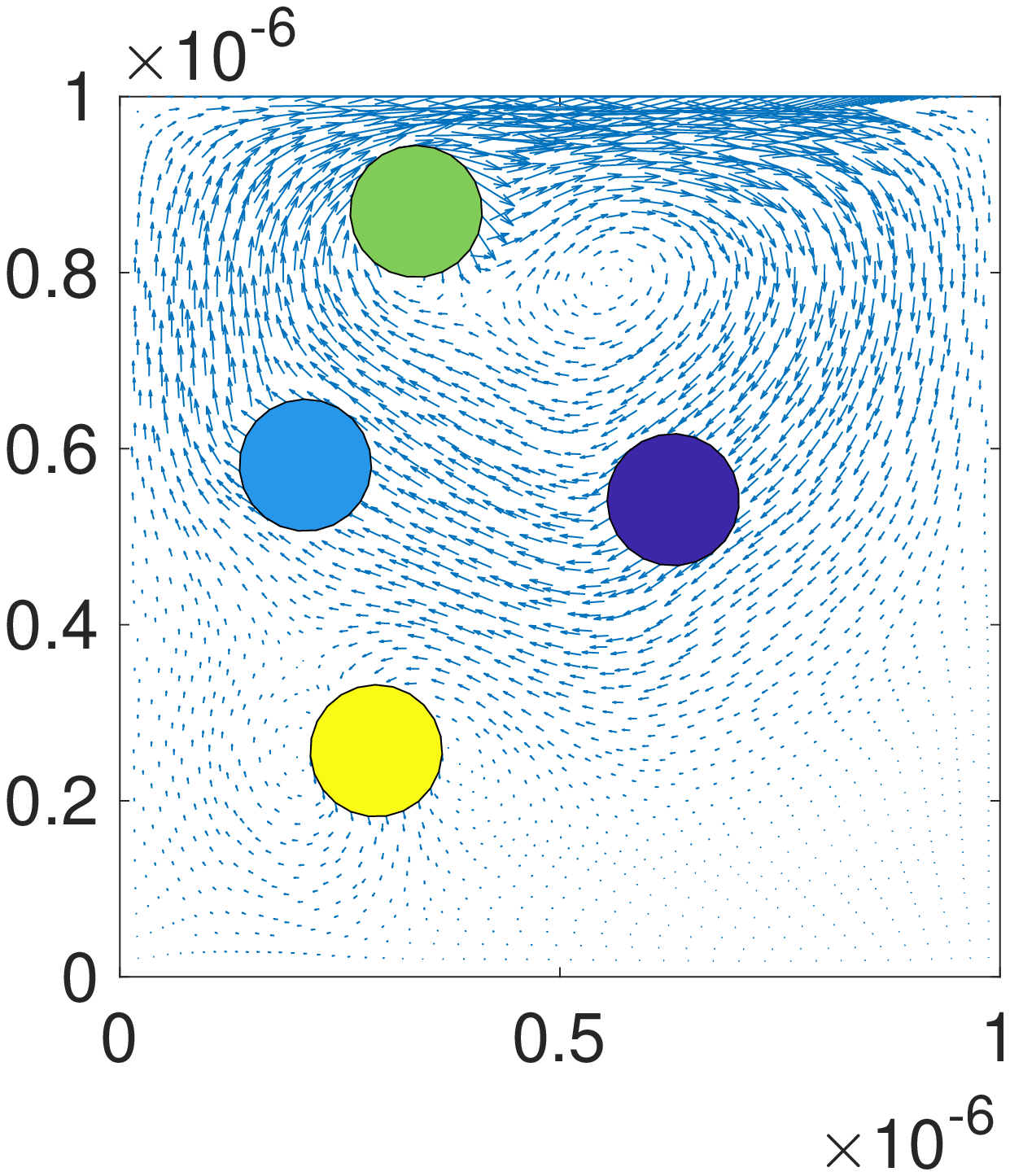}
	\caption{ Example 8: Multiple  particles in a driven cavity.  Particle positions and velocity  field. First row: $ t =0$ and $ t =1.5\cdot10^{-7}$. Second row $ t =3\cdot10^{-7}$ and  $ t = 6\cdot10^{-7}$.  }
	\label{circle_in_cavity}
	\centering
\end{figure}

\section{Conclusion and Outlook}
\label{sec:conclusion}
In this paper, we have presented an  Arbitrary Lagrangian-Eulerian method for the simulation of the BGK equation with moving boundaries.  
Besides the ALE approach, the method is based on   first and second order least squares approximations.
Several numerical tests are performed in order to validate the method, both in one and two space dimensions.
Moreover, we compared the results with those obtained by DSMC solution of the Boltzmann equation and by a higher order conservative semi-Lagrangian scheme.

In particular, in
 1D we consider the case of a moving plate immersed in a rarefied gas.
 In a first test we assume that the motion of the plate is prescribed (one way coupling), while in a second test the motion of the plate is computed from Newton's equations (two way coupling). 
% Notice that the DSMC results require to take the average of a lot of runs in order to decrease statistical fluctuations.  Thus, accurate DSMC solutions require a computational time which is several orders of magnitude higher than the one needed by the numerical solution of the BGK model.
% 
 In two space dimensions we considered several test problems. A first test case investigates  a situation where the motion of the object is prescribed
 (one -way coupling). We consider the motion of a  shuttle in a 2D model of a Micro Electro Mechanical System, see \cite{FFL}. Moreover, we considered  some tests with  rigid bodies/mesoscopic particles of arbitrary shape immersed in a gas and driven by either thermophoresis or driven cavity  flow (two way coupling).

In future work the scheme will be extended to the case of gas-mixtures
\cite{GR} and to three space dimensions.
Moreover, larger collections of mesoscopic particles dispersed in a rarefied gas will be considered, thus providing a quantitative tool that can be used to validate homogenised macroscopic models of suspensions.

%%%%%%%%%%%%%%%%%%%%%

%\begin{figure}[!t]
%\centering
%\includegraphics[scale=.5]{jcompfig01}
%\caption{Studio setup for capturing face images indoor. Three light
%sources L1, L2, L3 were used in conjunction with normal office lights.}
%\end{figure}

%%%%%%%%%%%%%%%%%%%%

\section*{Acknowledgments}
All authors would like to thank Dr.~Seung-Yong Cho for computing the numerical solution of the BGK model with a conservative semi-Lagrangian scheme.
 This  work is supported 
by the DFG (German research foundation) under Grant No. KL 1105/30-1 and by the ITN-ETN Marie-Curie Horizon 2020 program ModCompShock, Modeling and computation of shocks and interfaces, Project ID: 642768.
G.R.~would like to thank the Italian Ministry of
Instruction, University and Research (MIUR) to support this research with funds coming from PRIN Project 2017 (No.2017KKJP4X entitled Innovative numerical methods for evolutionary partial differential equations and applications).
G. Russo is a member of the INdAM Research group GNCS. 

%\bibliographystyle{plain}
%\bibliography{References}

\end{document}